\documentclass[aop,seceqn,citesort,MSNbibl,dvips]{arximspdf}
\usepackage{stfloats}
\usepackage{graphicx}

\doi{10.1214/09-AOP493}
\volume{38}
\issue{2}
\pubyear{2010}
\firstpage{714}
\lastpage{769}

\makeatletter

\fnbelowfloat

\newcommand{\Id}{\mathbf{1}}
\renewcommand{\widetilde}{\tilde}

\newtheorem{theorem}{Theorem}[section]

\newproclaim{rem}[theorem]{Remark}
\newtheorem{conjecture}[theorem]{Conjecture}
\newtheorem{proposition}[theorem]{Proposition}
\newtheorem{lemma}[theorem]{Lemma}
\newtheorem{corollary}[theorem]{Corollary}

\makeatother

\begin{document}
\begin{frontmatter}

\title{Airy processes with wanderers and new universality
classes\thanksref{T2}}
\runtitle{Airy processes with wanderers}

\thankstext{T2}{Supported by a European Science Foundation
grant (MISGAM), a Marie Curie Grant (ENIGMA), a FNRS grant and a
``Interuniversity Attraction Pole'' grant.}

\begin{aug}
\author[A]{\fnms{Mark} \snm{Adler}\thanksref{t1}\ead[label=e1]{adler@brandeis.edu}},
\author[B]{\fnms{Patrik L.} \snm{Ferrari}\ead[label=e2]{ferrari@uni-bonn.de}} and
\author[A,C]{\fnms{Pierre}
\snm{van Moerbeke}\corref{}\thanksref{t1}\ead[label=e3]{pierre.vanmoerbeke@uclouvain.be}\ead[label=e4]{vanmoerbeke@brandeis.edu}}
\runauthor{M. Adler, P. L. Ferrari and P. van Moerbeke}
\affiliation{Brandeis University, Bonn University and Universit\'{e} Catholique de
Louvain\break and Brandeis University}
\address[A]{M. Adler\\
P. van Moerbeke\\
Department of Mathematics\\
Brandeis University\\
Waltham, Massachusetts 02454\\
USA\\
\printead{e1}\\
\phantom{E-mail: }\printead*{e4}} 
\address[B]{P. L. Ferrari\\
Institute for Applied Mathematics\\
Bonn University\\
Endenicher Allee 60\\
53115 Bonn\\
Germany\\
\printead{e2}}
\address[C]{P. van Moerbeke\\
Department of Mathematics\\
Universit\'{e} Catholique de Louvain\\
1348 Louvain-la-Neuve\\
Belgium\\
\printead{e3}}
\end{aug}

\thankstext{t1}{Supported by NSF Grant DMS-07-04271.}

\received{\smonth{11} \syear{2008}}

%
\begin{abstract}
Consider $n+m$ nonintersecting Brownian bridges, with $n$ of them
leaving from $0$ at time $t=-1$ and returning to $0$ at time $t=1$,
while the $m$ remaining ones (wanderers) go from $m$ points $a_i$ to
$m$ points $b_i$. First, we keep $m$ fixed and we scale $a_i,b_i$
appropriately with $n$. In the large-$n$ limit, we obtain a new Airy
process with wanderers, in the neighborhood of $\sqrt{2n}$, the
approximate location of the rightmost particle in the absence of
wanderers. This new process is governed by an Airy-type kernel, with a
rational perturbation.

Letting the number $m$ of wanderers tend to infinity as well, leads to
two Pearcey processes about two cusps, a closing and an opening cusp,
the location of the tips being related by an elliptic curve. Upon
tuning the starting and target points, one can let the two tips of the
cusps grow very close; this leads to a new process, which might
be governed by a kernel, represented as a double integral involving
the exponential of a quintic polynomial in the integration variables.
\end{abstract}

%

\end{frontmatter}

\tableofcontents[alignleft,level=2]

\section{Introduction}\label{SectIntroduction}
Consider $n+m$ nonintersecting Brownian motions\break (Brownian bridges) on
${\mathbb R}$ depending on time $t\in[-1,1]$, with $n$ of them leaving
from and returning
to $0$, while the $m$ remaining ones leave from $a_m\leq
\cdots\leq a_1$ and are forced to end up at $b_m\leq\cdots\leq b_1$.
We denote by $x_i(t)$ the position at time $t$ of the $i$th largest
Brownian particle among the $n+m$ nonintersecting Brownian bridges.
Denote by $\mathcal{D}$ the conditioning event defined by the following
conditions:

\begin{longlist}
\item nonintersecting paths: $x_1(t)>x_2(t)>\cdots>x_{m+n}(t)$, $t\in
(-1,1)$,
\item
$n$ bridges from $0$ to $0$: $x_i(-1)=x_i(1)=0$ for $i=m+1,\ldots
,m+n$,
\item $m$ \textit{wanderers} from $a_i$ to $b_i$: $x_i(-1)=a_i$,
$x_i(1)=b_i$ for $i=1,\ldots,m$.
\end{longlist}

Then denote the conditional probability under $\mathcal D$ by ${\mathbb P}
_{\mathrm{ab}}$, that is,
%
\begin{equation}\label{1.1}
{\mathbb P}_{\mathrm{ab}}(\cdot)={\mathbb P}(\cdot| \mathcal{D}).
\end{equation}

The interest in nonintersecting Brownian motions stems from a paper by
Dyson~\cite{Dyson}, who made the important observation that
putting dynamics into the GUE-random matrix model (Ornstein--Uhlenbeck
Processes on the real and imaginary parts of the entries) leads to
finitely many
nonintersecting Brownian motions on ${\mathbb R}$ for the eigenvalues
(stationary process). A space--time transformation enables one to map
the above Dyson process into nonintersecting Brownian motions starting
from $0$ and returning to $0$; see formula (1.7) in \cite{AvMD08}. In
their work on coincidence probabilities, Karlin and McGregor \cite{Karlin}
found a determinantal formula for the transition probability of
nonintersecting Brownian motions. The relationship between
nonintersecting Brownian motions, matrix models and random matrix
theory has been developed starting with Johansson \cite{Johansson1} and
has led to the Airy and other processes
\cite{AptBleKui,AvM-Airy-Sine,DK,BleKui1,Brezin5,TW-Dyson,TW-Pearcey,Okounkov,Pastur,AvM-Pearcey,Zinn1,Zinn2},
when the number of particles tend to infinity, see also \cite{DelK}.

At first, consider the motion of the nonintersecting Brownian particles
above, but with $m=0$, and
let $n$ become very large. The \textit{Airy process} $\mathcal{A}(\tau)$
describes this cloud of particles (``infinite-dimensional diffusion''),
but viewed from any point on the ``edge'' $\mathcal{C}\dvtx x= \sqrt
{2n(1-t^2)}$ of the set of particles, with time and space properly
rescaled; the Airy process will be independent of the point chosen and
will be governed by the Airy kernel. This process was found by
Pr\"{a}hofer and Spohn \cite{Spohn} in the context of stochastic growth
models and further investigated in
\cite{Johansson2,Johansson3,TW-Dyson,AvM-Airy-Sine,FS02}.

Assume now a fixed and finite $m\geq1$ and all $a_i=0$, with the
target points all equal to $b$ scaled as
$b=\rho_0\sqrt{2n}>0$. Does it affect the Brownian
fluctuations along the curve $\mathcal{C}$ for large
$n$? No new process appears
as long as one considers points $(y,t)\in\mathcal{C}$,
below the point of tangency of the tangent to the curve
passing through $(\rho_0\sqrt{2n},1)$. At this tangency point
the fluctuations obey a new statistics, which we call
the \textit{Airy process with $m$ outliers} $\mathcal
{A}_{m}^b(\tau)$, governed by a rational perturbation of the Airy
kernel, see \cite{AvMD08}. This kernel was already considered by
Baik, Arous and P\'{e}ch\'{e} \cite{BBP,Baik} and P\'{e}ch\'{e} \cite
{Peche} in the context of multivariate statistics.

\begin{figure}

\includegraphics{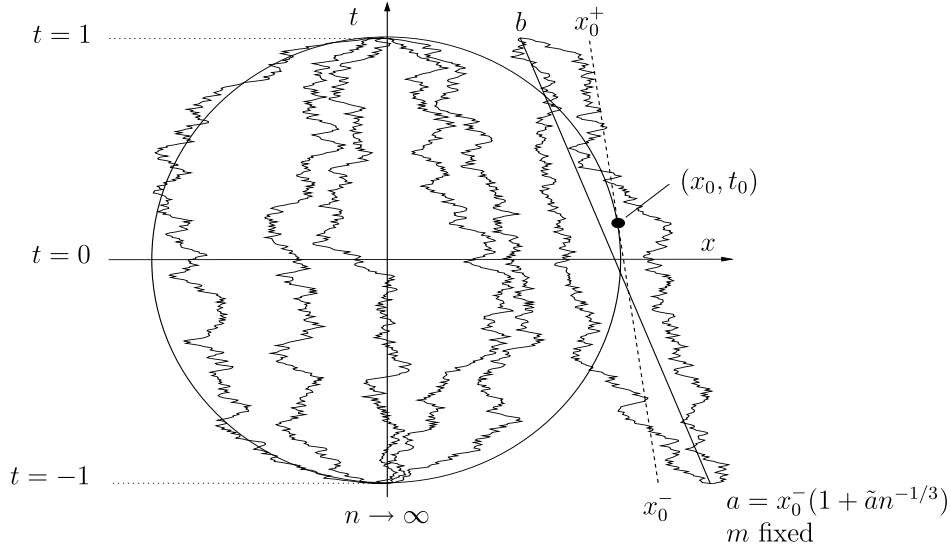}

\caption{Nonintersecting Brownian bridges with $m$
wanderers, leaving from $a=x_0^-(1+\tilde a n^{-1/3})$ and forced to
$b=x_0^+(1-\tilde b n^{-1/3})$, with $\tilde a<\tilde b$, where
$x_0^-=\sqrt{2n}
\sqrt{\frac{1+t_0}{1-t_0}}$, $x_0^+=\sqrt{2n}
\sqrt{\frac{1-t_0}{1+t_0}}$. The dotted line linking $(x_0^-,-1)$ to
$(x_0^+,1)$ is tangent to the curve $x= \sqrt{2n(1-t^2)}$ at the point
$(x_0,t_0)$.}
\label{FigAsymmetric}
\end{figure}

The first\vspace*{-2pt} result in this paper concerns the limiting process, described
in (\ref{1.1}), in the large-$n$ limit, while keeping $m$ fixed; this
process is denoted by $\mathcal{A}_m^{(\tilde a,\tilde b)}(\tau)$. This
paper deals with the statistical fluctuations of the edge of the cloud
of particles near any point on the curve $\mathcal{C}\dvtx x=\sqrt
{2n(1-t^2)}$, in the presence of wanderers.
To do so, consider the tangent line to the curve $\mathcal{C}$, with
point of tangency $(x_0,t_0)$, as in Figure \ref{FigAsymmetric}; this
tangent intersects the lines $t=-1$ and $t=1$ at the points
$x_0^-=\frac
{x_0}{1-t_0}=\sqrt{2n} \sqrt{\frac{1+t_0}{1-t_0}}$ and $x_0^+=\frac
{x_0}{1+t_0}=\sqrt{2n} \sqrt{\frac{1-t_0}{1+t_0}}$,
respectively. Consider now $m$ wanderers leaving from neighboring
points (when $n$ gets large) of the point $x_0^-$ at time $t=-1$ and
forced to neighboring points of $x_0^+$ at time $t=1$. The first part
of this paper is to show that the fluctuations near the edge of the
cloud and near the point $(x_0,t_0)$ obeys a new statistic, independent
of the point $(x_0,t_0)$ chosen on the curve above, showing \textit{universality} within that class.

%

\begin{figure}

\includegraphics{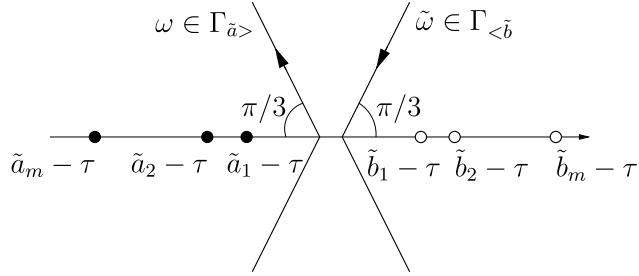}

\caption{Integration paths of the kernel $K_m^{\tilde a,\tilde b}$ of
(\protect\ref{eqKernel}).}\label{Fig1}
\vspace*{10pt}
\end{figure}

At the first stage (Theorems \ref{th:1} and \ref{th:2}), the result
will be shown for a vertical line tangent to $\mathcal{C}$ at the point
$(\sqrt{2n},0)$, whereas Theorem \ref{th:2'} deals with the
universality result.
The nonintersecting nature of the first $n$ bridges implies that the
largest one will again reach a height of about $\sqrt{2n}$. So, it is
natural to consider the following scaling of the starting and the
target points
%
\begin{equation}\label{scaling}
a_i=\sqrt{2n} \biggl(1+\frac{\tilde a_i}{n^{1/3}} \biggr)
\quad\mbox{and}\quad
b_i=\sqrt{2n} \biggl(1-\frac{\tilde b_i}{n^{1/3}} \biggr).
\end{equation}
With this scaling, the $m$ wanderers will interact with the bulk (of
$n$ particles, with $n$ very large) in a nontrivial way, upon
considering regions close to $x=\sqrt{2n}$ and $t=0$, namely at
space--time positions $(x,t)$ which scale like
%
\begin{equation}\label{scalingXT}
t=\tau n^{-1/3},\qquad x=\sqrt{2n}+\frac{\xi-\tau^2}{\sqrt{2}n^{1/6}}.
\end{equation}
This will only be so under some geometric condition: the lines
connecting the starting and target points in $(x,t)$-space must pass to
the left of $\sqrt{2n}$ at $t=0$; see Figure \ref{FigAsymmetric}.
Then the first result concerns the gap probability at a given time
$\tau
$ for very large $n$ and keeping $m$ finite and fixed, i.e., the
probability that a set is not visited by any of the $n+m$ Brownian
bridges at time $\tau$. Thus, in Theorems \ref{th:1} and \ref{th:2}, a
different (nontrivial) process $\mathcal{A}_m^{(\tilde a,\tilde
b)}(\tau)$ will appear due to the interaction of the $m$ wanderers with
the Airy field in the neighborhood of $(x,t)=(\sqrt{2n},0)$. Note that
in the absence of wanderers the particles must look, near the edge,
like the Airy process. This also explains why the kernel (\ref
{eqKernel}) obtained below is another rational extension of the Airy kernel.
\begin{theorem}\label{th:1}
Consider points $a_i$ and $b_i$, as in (\ref{scaling}), with $\tilde
a_m\leq\cdots\leq\tilde a_1<\tilde b_1\leq\cdots\leq\tilde b_m$
on the
real line.\footnote{The inequalities that all the $\tilde a_i$ be
smaller than all the $\tilde b_i$'s means geometrically that the lines
connecting corresponding points intersect the $x$-axis to the left of
$x=\sqrt{2n}$; see Figure \ref{FigAsymmetric}.\label{foot1}}
Given any compact set $E\subset{\mathbb R}$, the gap probability at rescaled
time--space (\ref{scalingXT}) is given, in the large-$n$ limit, by
%
%
\begin{eqnarray}\label{eq.1.4}
&&\lim_{n\to\infty} {\mathbb P}_{\mathrm{ab}} \biggl( \biggl\{\mbox{all } x_i \biggl(\frac
{\tau}{n^{1/3}} \biggr)\in\sqrt{2n}+\frac{E^c-\tau^2}{\sqrt
{2}n^{1/6}} \biggr\} \biggr) \nonumber\\
&&\qquad=\det(\Id-\chi_{_{E}} K_m^{\tilde a,\tilde b}\chi_{_{E}})_{L^2({\mathbb R})}\\
&&\qquad={\mathbb P}\bigl(\mathcal{A}_m^{(\tilde a,\tilde b)}(\tau)\cap
E=\varnothing\bigr),\nonumber
\end{eqnarray}
where $\chi_{_{E}}(\xi)=\Id(\xi\in E)$, where $\det$ denotes the Fredholm
determinant on $L^2({\mathbb R})$ and where the kernel $K_m^{\tilde
a,\tilde
b}$ is given by
%
%
\begin{eqnarray}\label{eqKernel}\hspace*{17pt}
&&K_m^{\tilde a,\tilde b}(\tau;\xi_1,\xi_2)
\nonumber\\
&&\qquad=\frac{1}{ (2\pi i )^2}
\int_{\Gamma_{\tilde a -\tau>}}d\omega\int_{\Gamma_{<\tilde
b-\tau
}}d\widetilde\omega\,\frac{e^{-{\omega^3}/3+\xi_2\omega
}}{e^{-
{\widetilde\omega^{3}}/3+\xi_1
\widetilde\omega}}\frac{1}{\omega-\widetilde\omega}
\\
&&\hspace*{156.8pt}{}\times\Biggl(\prod_{k=1}^m \biggl(\frac{\widetilde\omega-\tilde a_k+\tau
}{\omega-\tilde a_k+\tau}\biggr)
\biggl(\frac{\omega-\tilde b_k+\tau}{\widetilde\omega-\tilde b_k+\tau
}\biggr)\Biggr).\nonumber
\end{eqnarray}
The integration contours are as follows: $\Gamma_{\tilde a>}$ goes from
$e^{-2\pi i/3}\infty$ to $e^{2\pi i/3}\infty$, and passes on the right
of all the $\tilde a_i-\tau$, while $\Gamma_{<\tilde b}$ goes from
$e^{\pi i/3}\infty$ to $e^{-\pi i/3} \infty$, and passes to the left of
all $\tilde b_i-\tau$. Moreover, the two contours do not intersect; see
Figure \ref{Fig1} for an illustration.
\end{theorem}

This kernel has also appeared in recent work of Borodin and P\'{e}ch\'{e}
\cite{BP08}, as a limit of a directed percolation in a quadrant
with defective rows and columns, itself a generalization of a kernel of
Baik, Arous and P\'{e}ch\'{e} \cite{BBP,Baik} and P\'{e}ch\'{e}
\cite{Peche} and considered in \cite{AvMD08} in the context of
nonintersecting Brownian motions. The same limit process occurs in the
asymmetric exclusion process, see \cite{IS07,BFS09}.
The proof of Theorem \ref{th:1} will be given in Section \ref
{SectThm1}, when the points $a_i$ and the points $b_i$ are all
different. When the $a_i$'s all coincide, and similarly the $b_i$'s,
the proof of Theorem \ref{th:1} breaks down and must be replaced by
another one; two approaches are being discussed here (see Section \ref
{SectTwoPackets}): (1) using a certain moment matrix, (2)~using
biorthogonal polynomials.

%

In Theorem \ref{th:2} (see Section \ref{SectThm2}), the first result
will be extended to the joint gap probabilities at different (rescaled)
times $\tau_1,\ldots,\tau_\ell$. Obviously, Theorem \ref{th:1} is the
specialization of Theorem \ref{th:2} to the one-time case.
\begin{theorem}\label{th:2}
Consider $\ell$ distinct times $\tau_1,\tau_2,\ldots,\tau_\ell$ and
compact sets $E_1,\ldots,E_\ell\subset{\mathbb R}$. Then,
%
%
\begin{eqnarray}\label{0.Prob-extA}
&&\lim_{n\to\infty} {\mathbb P}_{\mathrm{ab}} \Biggl(\bigcap
_{k=1}^\ell\biggl\{
\mbox{all } x_i \biggl(\frac{\tau_k}{n^{1/3}} \biggr)\in\sqrt{2n}+\frac{E_k^c-\tau_k^2}{\sqrt{2}n^{1/6}} \biggr\} \Biggr)\nonumber\\
&&\qquad=\det(\Id-\chi_{_{E}} K_m^{\tilde a,\tilde b} \chi_{_{E}})\\
&&\qquad
={\mathbb P}\Biggl(\bigcap_{k=1}^\ell\bigl\{
\mathcal{A}_m^{(\tilde a,\tilde b)}(\tau_k)\cap
E_k=\varnothing\bigr\}\Biggr),\nonumber
\end{eqnarray}
where $\chi_{_{E}}(\tau_k,\xi):=\Id(\xi\in E_k)$. Here, $\det$
denotes the (matrix) Fredholm determinant
on the space $L^2(\{\tau_1,\ldots,\tau
_\ell
\}\times{\mathbb R})$ and the extended kernel $K_m^{\tilde a,\tilde
b}$ is
given~by\looseness=1
%
%
\begin{eqnarray}\label{eqKernelExtended}
&&K_m^{\tilde a,\tilde b}(\tau_1,\xi_1;\tau_2,\xi_2)\nonumber\\
&&\qquad=-\frac{\Id(\tau_2>\tau_1)}{\sqrt{4\pi(\tau_2-\tau_1)}}\nonumber\\
&&\qquad\quad\hspace*{6.2pt}{}\times e^{-{(\xi_2-\xi_1)^2}/({4(\tau_2-\tau_1)})-({1}/{2})
(\tau_2-\tau_1)(\xi_2+\xi_1)+({1}/{12}) (\tau_2-\tau_1)^3}\\
&&\qquad\quad{}+ \frac{1}{ (2\pi i )^2}\int_{\Gamma_{\tilde a-\tau_2>}}d\omega\int_{\Gamma_{<\tilde b-\tau_1}}
d\widetilde\omega\frac{e^{-\omega^3/3+\xi_2\omega}}{e^{-\widetilde\omega^{3}/3+\xi_1 \widetilde\omega}}
\frac{1}{(\omega+\tau_2)-(\widetilde\omega+\tau_1)}\nonumber\\
&&\qquad\quad\hspace*{140.5pt}{}\times
\Biggl(\prod_{k=1}^m \biggl(\frac{\widetilde\omega-\tilde a_k+\tau_1}{\omega-\tilde a_k+\tau_2} \biggr)
\biggl(\frac{\omega-\tilde b_k+\tau_2}{\widetilde\omega-\tilde b_k+\tau_1} \biggr)\Biggr).\nonumber
\end{eqnarray}
The integration contours are as in Figure \ref{Fig1}, but with $\tilde
a_k-\tau$ replaced by $\tilde a_k-\tau_2$ and $\tilde b_k-\tau$
replaced by $\tilde b_k-\tau_1$.
\end{theorem}

A similar statement can then be made along any point $(x_0,t_0)$ of the
curve $x=\sqrt{2n(1-t^2)}$, with tangent intersecting the lines $t=-1$
and $t=1$ at the points
%
\begin{equation}\label{points}
\quad x_0^-=\frac{x_0}{1-t_0}=\sqrt{2n} \sqrt{\frac{1+t_0}{1-t_0}}
\quad\mbox{and}\quad x_0^+=\frac{x_0}{1+t_0}=\sqrt{2n} \sqrt{\frac
{1-t_0}{1+t_0}},
\end{equation}
respectively. This is done in Theorem~\ref{th:2'} below.
\begin{theorem}[(Universality statement)]\label{th:2'}
As before, consider $\ell$ distinct times $\tau_1,\tau_2,\ldots,\tau
_\ell$ and compact sets $E_1,\ldots,E_\ell\subset{\mathbb R}$.
Also, consider
$m$ Brownian wanderers, now leaving from the points $a_\ell
=x_0^-(1+\frac{\tilde a_\ell}{n^{1/3}})$ and forced to $b_\ell
=x_0^+(1-\frac{\tilde b_\ell}{n^{1/3}})$, with the condition\footnote
{Here also, the inequalities that all the $\tilde a_i$ be smaller than
all the $\tilde b_i$'s means geometrically that the lines connecting
corresponding points intersect the horizontal line through $(x_0,t_0) $
to the left of $(x_0,t_0) $; see Figure \ref{FigAsymmetric}.\label
{foot1'}} $\tilde a_m\leq\cdots\leq\tilde a_1<\tilde b_1\leq\cdots
\leq
\tilde b_m$.
For $n$ large, pick $\ell$ points in a $n^{-1/3}$-neighborhood of
$(x_0,t_0)$, lying on the curve $x=\sqrt{2n(1-t^2)}$,
%
%
\begin{equation}
x_k:=\sqrt{2n(1-t_k^2)}\qquad \mbox{with } t_k:= t_0+\frac{(1-t_0^2) \tau
_k}{n^{1/3}}, 1\leq k\leq\ell.
\end{equation}
Then the following limit holds\footnote{Expanded out, $ (1+\frac
{E_k^c}{2n^{2/3}} )x_k$ reads
\[
\sqrt{1-t_0^2}\sqrt{2n} \biggl(1-\frac{\tau_k t_0}{n^{1/3}}+\frac
{E_k^c-\tau_k^2}{2n^{2/3}}-t_0\tau\frac{ E_k^c+\tau_k^2 }{2n}
\biggr)+O \biggl(\frac{1}{{n}^{5/6}} \biggr).
\]
\label{foot3''}}:
%
%
\begin{eqnarray}\label{0.Prob-ext}
&&\lim_{n\to\infty} {\mathbb P}_{\mathrm{ab}} \Biggl(\bigcap_{k=1}^\ell\biggl\{
\mbox{all } x_i (t_k )\in
\biggl(1+\frac{E_k^c}{2n^{2/3}} \biggr)x_k \biggr\}
\Biggr)\nonumber\\[-8pt]\\[-8pt]
&&\qquad={\mathbb P}\Biggl(\bigcap_{k=1}^\ell\bigl\{\mathcal{A}_m^{(\tilde a,\tilde
b)}(\tau_k)\cap E_k=\varnothing\bigr\}
\Biggr).\nonumber
\end{eqnarray}
\end{theorem}
\begin{rem}
For $(x_0,t_0)=(\sqrt{2n},0)$, this statement reduces to Theorem \ref
{th:2}, as can be seen from footnote \ref{foot3''}.
\end{rem}

\begin{figure}

\includegraphics{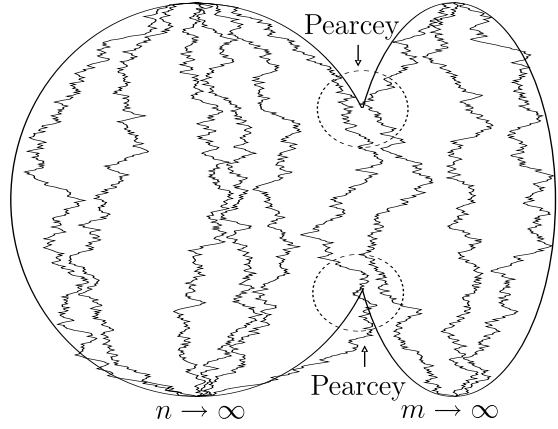}

\caption{Illustration of the two Pearcey processes, arising around the
two cusps.}
\label{FigPearcey}
\end{figure}

In view of the new process $\mathcal{A}_m^{(\tilde a,\tilde b)}(\tau)$,
it seems natural to let the number of wanderers $m$ to go infinity. For
simplicity, consider the case where the $m$ wanderers all start from
the same point $\tilde a$, and end up at the same point $\tilde b$,
with $\tilde a<\tilde b$, with the scaling
%
%
\begin{equation}\label{scalingAB}
\tilde a=\alpha m^{1/3},\qquad \tilde b=\beta m^{1/3}\qquad \mbox{with
}\alpha<\beta.
\end{equation}
Under this scaling, the set of $m$ wanderers itself produces an Airy
field which then interacts with the one already present after the $n\to
\infty$ limit. Thus, we might expect that there will be two regions
where the Pearcey process arises. Indeed, the first Pearcey process
occurs when the ``Pearcey cusp'' closes, while the second does when the
cusp opens, as illustrated in Figure \ref{FigPearcey}.

%

%
\begin{figure}[b]

\includegraphics{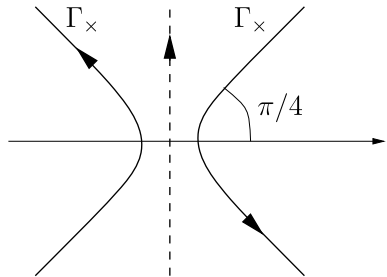}

\caption{Integration paths of the Pearcey kernel $K^\mathcal{P}_{u}$
defined in (\protect\ref{eqPearceyKernel}). The two solid lines form together
$\Gamma_\times$, the dashed line is the $z$-integration path.}
\label{Fig2}
\end{figure}

The reader is reminded of the extended Pearcey kernel $K^\mathcal
{P}(\theta_1,v_1;\theta_2,v_2)$ with space--time parameters $(\theta
_i,v_i)$, which is given by
%
%
\begin{eqnarray}\label{eqPearceyKernel}\qquad\quad
K^\mathcal{P}(\theta_1,v_1;\theta_2,v_2) &=& -\frac{\Id(\theta
_2>\theta
_1)}{\sqrt{2\pi(\theta_2-\theta_1)}} e^{-(v_2-v_1)^2/2(\theta
_2-\theta
_1)} \nonumber\\[-8pt]\\[-8pt]
&&{} + \frac{1}{(2\pi i)^2}\int^{i\infty}_{-i\infty}dz \int_{\Gamma
_\times
}d\tilde z\frac{1}{z-\tilde z} \frac{e^{-z^4/4+\theta_2 z^2/2-v_2
z}}{e^{-\tilde z^4/4+\theta_1\tilde z^2/2-v_1\tilde z}},\nonumber
\end{eqnarray}
where the path $\Gamma_\times$ is illustrated in Figure \ref{Fig2}; see
Tracy and Widom \cite{TW-Pearcey}. This leads to Theorem \ref{th:Pearcey},
established in Section \ref{SectPearcey}.
\begin{theorem}\label{th:Pearcey}
Let the starting point $\tilde a$ and the target point $\tilde b$ of
the $m$ wanderers for the Airy process with wanderers (\ref
{0.Prob-ext}) grow with m, as $\tilde a=\alpha m^{1/3}$ and $\tilde b=
\beta m^{1/3}$ with arbitrary $\alpha<\beta$. Given $\alpha<\beta$, the
following equations:
%
\begin{eqnarray}\label{ell_curve}
\beta-\alpha={\frac{4{\sigma}^{4}{x}^{3}}{2-x}}\hspace*{80pt}
\nonumber\\[-8pt]\\[-8pt]
\eqntext{\mbox{with } (x,\sigma)\in\mathcal{E}\dvtx 4\sigma^6 x^4-2 x+ 3=0
\mbox
{ (elliptic curve)}}
\end{eqnarray}
have a
\begin{eqnarray*}
\mbox{unique solution }(x,\sigma)
&:=&
(x,\sigma_+)\in\bigl(\bigl(\tfrac32,2\bigr)\times\bigl(-\tfrac12,0\bigr)\bigr)\mbox{ (\textup{opening cusp}),}
\\
\mbox{unique solution }(x,\sigma)
&:=&
(x,\sigma_-)\in\bigl(\bigl(\tfrac32,2\bigr)\times\bigl(0,\tfrac12\bigr)\bigr)\mbox{ (\textup{closing cusp})}.
\end{eqnarray*}
Then, the Airy process with $m$ wanderers $\mathcal{A}_m^{(\tilde
a,\tilde b)}(\tau)$ properly rescaled as \mbox{$m\to\infty$}, converges
to two (identical) Pearcey processes $\mathcal{P}(\theta)$ about two
cusps, one \textup{opening cusp} ($T_+$) and one \textup{closing cusp}
($T_-$) about \footnote{$T_+$ corresponds to $\sigma_+<0$ and $T_-$
corresponds to $\sigma_->0$, with obviously $\sigma_+=-\sigma_-$.}
%
%
\begin{eqnarray}\label{5.XT}
\tau\sim T_{\pm} m^{1/3},\qquad \xi\sim X m^{2/3}\nonumber\\[-8pt]\\[-8pt]
\eqntext{\mbox{with
} \displaystyle T_{\pm}:=
\frac{\alpha+\beta}2- \frac{2\sigma_\pm}{2-x}, X:=\sigma_\pm^2(1-2x),}
\end{eqnarray}
with $T_-<\frac{\alpha+\beta}2<T_+$. To be precise, upon using the two
different scalings (\ref{1.scaling}) below, depending on the opening or
closing cusp, one has, for any $\ell=1,2,\ldots,$ that the limit of the
gap probability of the sets $\widetilde E_1,\ldots,\widetilde E_\ell$
at times $\tau_1,\ldots,\tau_\ell$ is given by the same (matrix)
Fredholm determinant,
%
%
\begin{eqnarray}\label{1}\hspace*{27pt}
&&\lim_{m\to\infty}
{\mathbb P}\Biggl(\bigcap_{k=1}^\ell\bigl\{\mathcal{A}_m^{(\alpha m^{1/3},\beta
m^{1/3})}(\tau_k)
\cap\widetilde E_k=\varnothing\bigr\}\Biggr)\nonumber\\
&&\qquad=\det(\Id-\chi_{_{E}} K^\mathcal{P}\chi_{_{E}} )_{L^2((\theta
_1,\ldots,\theta_\ell)\times{\mathbb R})}\\
&&\qquad=:
{\mathbb P}\Biggl(\bigcap_{k=1}^\ell\{\mathcal{P}(\theta_k)\cap
E_k=\varnothing\} \Biggr),\nonumber
\end{eqnarray}
where the rescaling from the space--time variables $(\tilde E_i,\tau_i)$
to the new space--time variables $( E_i,\theta_i)$ is imposed by the
initial scaling (\ref{5.XT}), to yield
%
%
\begin{eqnarray}\label{1.scaling}
\tau_i&=& T_\pm m^{1/3} + \frac12 \kappa^2 \theta_i m^{-1/6},
\qquad \kappa := \biggl(\frac{2(x-1)}{|\sigma_\pm| x^2}
\biggr)^{1/4},\nonumber\\[-8pt]\\[-8pt]
\widetilde E_i &=& X m^{2/3}-\kappa^2\sigma_\pm\theta_i
m^{1/6}-\kappa
E_i m^{-1/12}.\nonumber
\end{eqnarray}
\end{theorem}
\begin{rem}
Note that the involution: $v_1\leftrightarrow v_2$, $\theta_1
\leftrightarrow-\theta_2$, $T_+\leftrightarrow T_-$, $\sigma_+
\leftrightarrow\sigma_-=-\sigma_+$, where $v_k\in E_k$, maps the
opening cusp into the closing cusp and, in particular, acts on the
kernel (\ref{eqPearceyKernel}) to produce the kernel going with the
closing cusp.
\end{rem}

\begin{figure}

\includegraphics{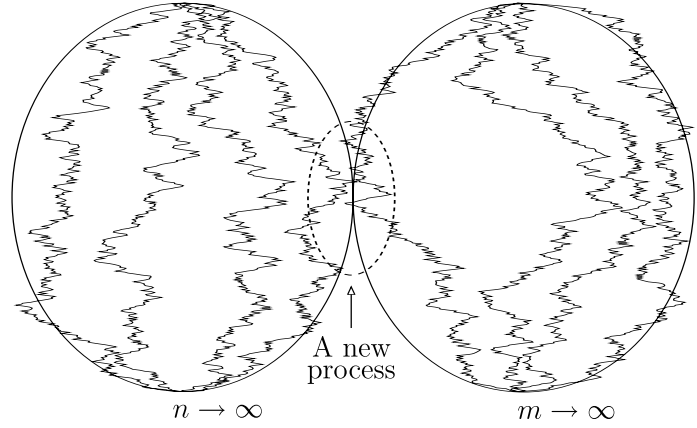}

\caption{When two Pearcey cusps touch, there will be a new process.}
\label{FigNewProcess}
\end{figure}

The tips of the two cusps in Theorem \ref{th:Pearcey} come together,
when $\alpha, \beta\to0$, and hence
$x\to3/2$, $\sigma_\pm\to0$ and $T_\pm\to0$; this is not the
only way for this to happen, as will be mentioned below.
At the very point where the two cusps meet, a new process will emerge
(as in Figure \ref{FigNewProcess}), which might be
%
governed by a ``\textit{quintic kernel}.''
\begin{conjecture}\label{conjecture}
The gap probability for the new process appearing in Figure \ref
{FigNewProcess} is given by the Fredholm determinant of the following
\textit{quintic kernel}:
%
\begin{equation}\label{Quintic}
K^\mathcal{Q}(\theta,\eta;x,y) =\frac{1}{(2\pi i)^2}\int_{
\mathcal
{C}}dz\int_{\widetilde{\mathcal{C}}}d\tilde z \frac{1}{z-\tilde z}
\frac{e^{2z^5/5-\theta z^3/3-\eta z^2+z x}}{e^{2\tilde z^5/5-\theta
\tilde z^3/3-\eta\tilde z^2+\tilde z y}},
\end{equation}
where the $z$ and $\tilde z$-integration paths are given by appropriate subpaths of the $z$ and
%
\begin{figure}[b]

\includegraphics{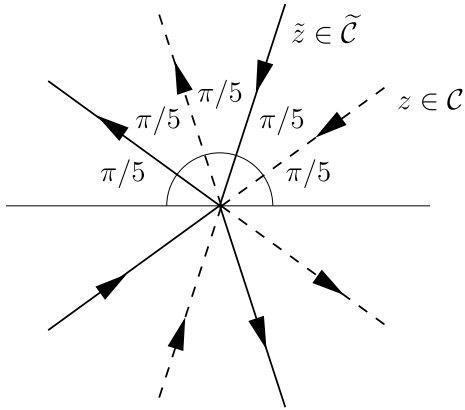}

\caption{Integration path $\mathcal{C}$ (dashed line) and $\widetilde
{\mathcal{C}}$ (solid line) of the quintic kernel.}\label{QuinticProcess}
\end{figure}
$\tilde z$-paths in Figure \ref{QuinticProcess}, with the orientation
indicated.
\end{conjecture}

To explain this attempt, we first notice that the
curve $\mathcal E$ [introduced in (\ref{ell_curve})] contains another
real point [besides the real segments introduced just after (\ref
{ell_curve})] appropriate subpaths of namely at $(x,\sigma)= (\infty,0)$, for which $(\alpha
,\beta)=
(2^{1/3},-2^{1/3})$; there the critical point $w_c$ of the associated
steepest-descent \mbox{$F$-function} becomes order $5$, with $(X,T)=
(-2^{2/3},0)$, rather than order $4$ as in the Pearcey case; this
expresses the fact that the two tips come together. For this choice of
$(\alpha,\beta)= (2^{1/3},-2^{1/3})$, the source and the target points\vspace*{2pt}
\begin{eqnarray}
a&=&\sqrt{2n} \biggl(1+\frac{\tilde a}{n^{1/3}}
\biggr),\qquad
b=\sqrt{2n} \biggl(1-\frac{\tilde b}{n^{1/3}} \biggr)\nonumber\\[2pt]
\eqntext{\mbox{with } \tilde a=\alpha m^{1/3} \mbox{ and } \tilde b= \beta m^{1/3},}
\end{eqnarray}
do not, of course, satisfy the inequality $\tilde a<\tilde b$, but
rather the opposite inequality.
We then perform an analytic continuation of the (one-time)\vspace*{2pt}
kernel\footnote{In the one-time case, one can just absorb the time
$\tau
$ in the $\tilde a$ and $\tilde b$.}
%
%
\begin{eqnarray}\label{eqKernelExtended1}\qquad
&&K_m^{\tilde a,\tilde b}( \xi_1; \xi_2)
\nonumber\\[2pt]
&&\qquad= \frac{1}{ (2\pi i )^2}
\int_{\Gamma_{\tilde a >}}d\omega\int_{\Gamma_{<\tilde
b}}d\widetilde
\omega\frac{e^{-\omega^3/3+\xi_2\omega}}{e^{-\widetilde\omega
^{3}/3+\xi_1
\widetilde\omega}}\\[2pt]
&&\hspace*{137.7pt}{}\times
\frac{ (({\widetilde\omega-\tilde a})/({\omega-\tilde a}))^m
(({\omega-\tilde b})/({\widetilde\omega-\tilde b} ))^m}{ \omega
-\widetilde\omega}\nonumber
\end{eqnarray}
by moving\vspace*{2pt} $\tilde a$ and $\tilde b$ in the complex plane from their
original position $\tilde a<\tilde b$ to a new position $\tilde
b<\tilde a$ on the real line. Then by picking $\tilde a=\alpha m^{1/3}$ and
$\tilde b= \beta m^{1/3}$, with $(\alpha,\beta)= (2^{1/3},-2^{1/3})$ and
letting $m\rightarrow\infty$, we show the kernel (\ref{eqKernelExtended1}) tends to the quintic kernel (\ref{Quintic}) with
the precise contour of integration in the figure above. Some evidence
in favor of this guess is given in Section \ref{SectConjecture},
which contains two rigorous statements, with proofs. However, this does
not suffice to prove the conjecture; e.g., it is still unknown whether
the Fredholm determinant of the quintic kernel (\ref{Quintic})
determines a probability. For numerical methods, see, for instance,
Bornemann \cite{Bornemann}. Folkmar Bornemann and Georg
Wechslberger developed a Mathematica-program to numerically compute
the kernels obtained above. The full paths as in Figure~\ref{QuinticProcess} did
not pass the positivity test. However, there are many other
possibilities of selecting the paths and/or their orientations, some of
which have positive density.
\begin{rem}
It is interesting to put the three kernels in parallel,
%
%
\begin{figure}

\includegraphics{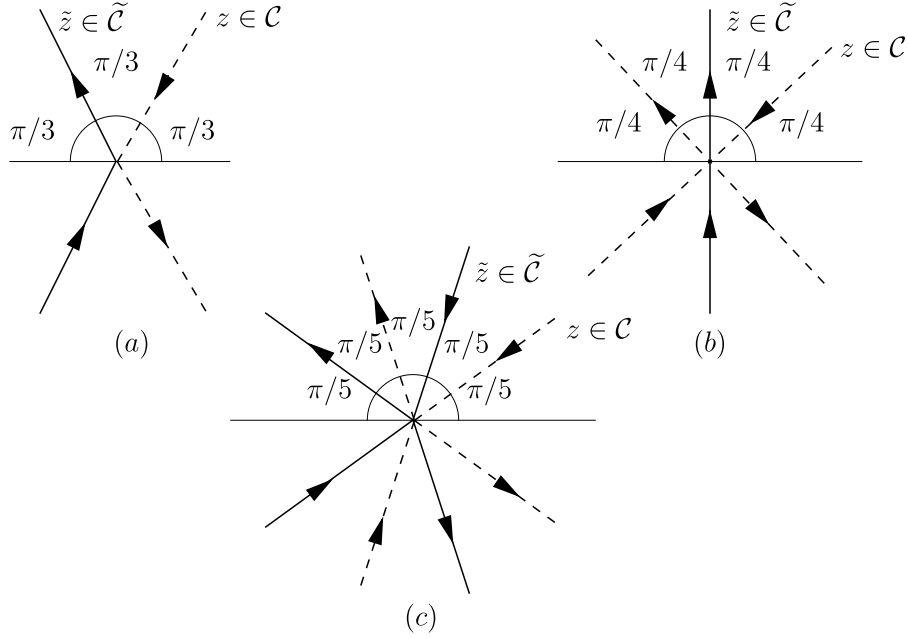}

\caption{Integration paths $\mathcal{C}$ (dashed line) and
$\widetilde
{\mathcal{C}}$ (solid line) of the \textup{(a)} Airy kernel $K^\mathcal{A}$,
\textup{(b)}
Pearcey kernel $K^\mathcal{P}$ and \textup{(c)} Quintic kernel $K^{\mathcal{Q}}$.}
\label{345Kernels}
\end{figure}
Airy, Pearcey and quintic together with their appropriate contours, as
in Figure \ref{345Kernels}:\vspace*{2pt}
\begin{eqnarray*}
&&K^\mathcal{A}(\tau_1,\xi_1;\tau_2,\xi_2)\\[2pt]
&&\qquad= \frac{1}{ (2\pi i )^2}
\int_{\tilde\mathcal{C}}d\tilde z \int_{\mathcal{C}}d z %
\frac{1}
{(\tilde z+\tau_2)-( z+\tau_1)} \frac{e^{ z^{3}/3-\xi_1
z}}{e^{\tilde
z^3/3-\xi_2\tilde z}}
\\[2pt]
&&\qquad\quad{} -\frac{\Id(\tau_2>\tau_1)}{\sqrt{4\pi(\tau_2-\tau_1)}}
e^{-{(\xi_2-\xi_1)^2}/({4(\tau_2-\tau_1)})-({1}/{2})(\tau
_2-\tau
_1)(\xi_2+\xi_1)+({1}/{12}) (\tau_2-\tau_1)^3},
\\
&&K^\mathcal{P}(\tau_1,\xi_1;\tau_2,\xi_2) \\
&&\qquad= \frac{1}{(2\pi
i)^2}\int
_{\tilde\mathcal{C}}d\tilde z \int_{\mathcal{C}}d z\frac{1}{\tilde
z- z}
\frac{e^{ z^4/4-\tau_1 z^2/2+\xi_1 z}}{e^{\tilde z^4/4-\tau_2
\tilde
z^2/2+\xi_2 \tilde z}}
\\
&&\qquad\quad{} -\frac{\Id(\tau_2>\tau_1)}{\sqrt{2\pi(\tau_2-\tau_1)}}
e^{-(\xi
_2-\xi_1)^2/(2(\tau_2-\tau_1))},
\\
&&K^\mathcal{Q}(\tau,\eta,\xi_1;\tau,\eta;\xi_2)\\
&&\qquad=\frac
{1}{(2\pi
i)^2}\int_{\widetilde{\mathcal{C}}}d\tilde z\int_{ \mathcal{C}}dz
\frac
{1}{z-\tilde z}
\frac{e^{2z^5/5-\tau z^3/3-\eta z^2+\xi_1 z }}{e^{2\tilde z^5/5-\tau
\tilde z^3/3-\eta\tilde z^2+\xi_2 \tilde z }}.\hspace*{36.4pt}\hspace*{36.4pt}
\end{eqnarray*}
\end{rem}

\section{Airy process with wanderers leaving from and going to distinct
points}\label{SectThm1}
The aim of this section is to prove Theorem \ref{th:1} in the case that
all points $\tilde a_i$ are distinct and all $\tilde b_i$ as well. We
first present the case $\tau=0$. The multi-time joint gap probabilities
will be discussed in the next section, implying the case of the
one-time process; i.e., general $\tau$ beyond $\tau=0$. However, first
presenting the one-time case will prove useful for understanding the
basic structure.

Denote by $p(x,y;t)$ the one-particle Brownian motion transition from
$x$ to $y$ during a time interval $t$, namely
%
\begin{equation}
p(x,y;t)=\frac{1}{\sqrt{2\pi t}}e^{-{(x-y)^2}/({2t})}.
\end{equation}
Let us consider $n+m$ Brownian bridges leaving at $t=-1$ from
$a_{m+n}<\cdots<a_{m+1}<a_m<\cdots<a_1$ and ending at $t=1$ at
positions $b_{m+n}<\cdots<b_{m+1}<b_m<\cdots<b_1$. The positions of
these particles at time $t$ are denoted by $\mathbf{x}(t)=\{
x_1(t),\ldots,x_{m+n}(t)\}$. Then, the probability density that
$\mathbf{x}(t)=\mathbf{x}$, conditioned that the Brownian bridges do not
intersect in $t\in(-1,1)$, is given by the Karlin and McGregor formula
\cite
{Karlin}, namely
%
%
\begin{eqnarray}\label{eq2.1}
P\bigl(\mathbf{x}(t)=\mathbf{x}\bigr)&=&\frac{1}{Z}\det\bigl(p(a_i,x_j;1+t)\bigr)_{1\leq
i,j\leq m+n}\nonumber\\[-8pt]\\[-8pt]
&&{}\times\det\bigl(p(x_i,b_j;1-t)\bigr)_{1\leq i,j\leq m+n},\nonumber
\end{eqnarray}
with $Z$ the normalization constant, which is equal to the probability
that the $m+n$ paths do not intersect, given the initial and final
conditions at $t=\pm1$.

It is known that a measure on $x=(x_1,\ldots,x_{m+n})$ of the form
(\ref
{eq2.1}) has determinantal correlation functions (see, e.g.,
Proposition 2.2 of \cite{Bor98}, or for information on determinantal
processes \cite{BKPV05,Jo05,Lyo03,Sos00,Spo05}).

As mentioned before, we restrict the discussion in this section to the
case $t=0$. Then the $k$-point correlation functions $\rho^{(k)}$ are
given by
%
%
\begin{equation}
\rho^{(k)}(x_1,\ldots,x_k)=\det(K(x_i,x_j) )_{1\leq i,j\leq k}
\end{equation}
with the kernel $K$ explicitly given by
%
%
\begin{equation}
K(x,y)=\sum_{i,j=1}^{m+n}p(x,b_i;1) [B^{-1}]_{i,j} p(a_j,y;1),
\end{equation}
where
%
%
\begin{equation}
B=[B_{i,j}]_{1\leq i,j\leq m+n},\qquad B_{i,j}=\int_{{\mathbb R}}dx\, p(a_i,x;1)
p(x,b_j;1).
\end{equation}
In particular, the gap probability of a set $E$, that is, the probability
that none of the $x_1,\ldots,x_{n+m}$ belongs to the set $E$, is given
in terms of a Fredholm determinant,
%
%
\begin{equation}\label{eq2.6}\qquad
{\mathbb P}(\mbox{none of the }x_i\in E)=\det(\Id-\chi_{_{E}} K
\chi_{_{E}})_{L^2({\mathbb R})},\qquad \chi_{_{E}}(x)=\Id(x\in E).
\end{equation}

The structure of the measure does not change when taking the limit of
one of more of the Brownian bridges starting and/or leaving from the
same position. Thus, the determinantal structure of correlation still
holds, yielding the
following proposition.
\begin{proposition} \label{Prop:1.1}
Consider $a_{m+1}=\cdots=a_{m+n}=0$ and $b_{m+1}=\cdots=b_{m+n}=0$ and
the other $m$ Brownian bridges from $a_i$ to $b_i$, with $0<a_m<\cdots
<a_1$ and $0<b_m<\cdots<b_1$. Then
%
\begin{equation}\label{1.Prob}
{\mathbb P}\bigl(\mathbf{x}(0)\notin E\bigr)=\det(\Id-\chi_{_{E}} K_{n,m} \chi
_{_{E}} ),
\end{equation}
where the kernel $K_{n,m}$ is given by
%
\begin{equation}\label{finite_kernel}
K_{n,m}(x,y)=K_n^{\mathrm{Hermite}}(x,y)+\sum^m_{i,j=1}\psi
^{(n)}_i(x)(\mu
^{-1})_{ij}\varphi^{(n)}_j(y).
\end{equation}
The Hermite kernel $K_n^{\mathrm{Hermite}}$ is defined by the classical
Hermite polynomials and their $L^2$-norms\footnote{$ \int_{{\mathbb R}}dx H_k(x)H_{\ell}(x)e^{-x^2}=\delta
_{k,\ell
} c_k^2
= \delta_{k,\ell}{2^{k}k!}\sqrt{\pi}$.}
%
\begin{equation}
K_n^{\mathrm{Hermite}}(x,y)=e^{-(x^2+y^2)/2}\sum^{ n-1}_{i=0} \frac
{1}{c_i^2} H_{i}(x)H_{i} (y);
\end{equation}
the functions $\psi_k^{(n)}$ and $\varphi_k^{(n)}$ are defined as follows
for $1\leq k\leq m$:
%
\begin{eqnarray}\label{1.phi}
\varphi_k^{(n)}(x) &=& {\frac{e^{-x^2/2}}{2\pi i} \oint_{\Gamma
_{0,a/2}}
dz\frac{e^{-z^2+2xz}}{z^n(z-{a_k}/{2})}},\nonumber\\[-8pt]\\[-8pt]
\psi_k^{(n)}(x) &=& {\frac{e^{-x^2/2}}{2\pi i} \oint
_{\Gamma_{0,b/2}}
dz\frac{e^{-z^2+2xz}}{z^n(z-{b_k}/{2})}},\nonumber
\end{eqnarray}
where $\Gamma_{0,a/2}$ denotes any contour containing the points
$z=0,a_1/2,\ldots,a_m/2$, and similarly for $\Gamma_{0,b/2}$. Finally,
the entries of the matrix of inner products
%
\begin{equation}\label{1.mu}\qquad\quad
\mu=(\mu_{k\ell})_{1\leq k,\ell\leq m}\qquad\mbox{with } \mu_{k\ell
}=\bigl\langle\varphi
_k^{(n)},\psi_{\ell}^{(n)}\bigr\rangle\equiv\int_{\mathbb R}dx\, \varphi
_k^{(n)}(x)\psi
_\ell^{(n)}(x)
\end{equation}
can be written\footnote{Similarly $\Gamma_{0,({a_kb_{\ell}})/{2}}$ denotes a
contour containing $0$ and $a_kb_{\ell}/2$. Note that $\frac{1}{2\pi
i}\oint_{\Gamma_{0,u}} \frac{e^z\,dz}{z^n (z-u)}=\frac{1}{u^n}(\sum
_{k=n}^\infty\frac{u^k}{k!})$.}
%
\begin{equation}\label{1.mukl}
\mu_{k\ell}=\frac{\sqrt{\pi} 2^n}{2\pi
i}\oint_{\Gamma_{0,a_k b_{\ell}/2}}dz \frac{e^z}{z^n (z-a_k b_{\ell
}/2 )}.
\end{equation}
\end{proposition}
\begin{pf}
We start from the setting (\ref{eq2.1})--(\ref{eq2.6}) and take
the limit when the $2n$ points $a_{m+n},\ldots,a_{m+1} \to0$ and
$b_{m+n},\ldots,b_{m+1} \to0$, and leaving the $2m$ points
$a_m<\cdots
<a_1$ and $b_m<\cdots<b_1$ fixed. Then the probability density on the
$x_i$'s becomes
%
%
\begin{eqnarray}\label{eq2.9}
P\bigl(\mathbf{x}(0) = \mathbf{x}\bigr)&=&\frac{1}{Z'}\det\pmatrix{
(e^{a_ix_j-x_j^2/2})_{\fontsize{8.36}{10}\selectfont{\begin{array}{l}
1\leq i\leq m\\ 1\leq
j\leq m+n
\end{array}}} \cr
(x_j^{i-1}e^{-x_j^2/2})_{\fontsize{8.36}{10}\selectfont{\begin{array}{l}
1\leq i\leq n\\ 1\leq j\leq m+n
\end{array}}}}
\nonumber\\[-8pt]\\[-8pt]
&&{}\times
\det\pmatrix{(e^{b_ix_j-x_j^2/2})_{\fontsize{8.36}{10}\selectfont{\begin{array}{l}
1\leq i\leq m\\ 1\leq j\leq m+n
\end{array}}}
\cr
(x_j^{i-1}e^{-x_j^2/2})_{\fontsize{8.36}{10}\selectfont{\begin{array}{l}
1\leq i\leq n\\ 1\leq j\leq m+n
\end{array}}}},\nonumber
\end{eqnarray}
where $Z'$ is a normalization constant. Consider any set of functions
$\{\varphi_k^{(n)}(x),k=1,\ldots,n+m\}$ spanning the vector space
%
%
\begin{equation}\label{vspace}\qquad
V(a_1,\ldots,a_m)=\operatorname{span}\{e^{a_ix-x^2/2},1\leq i \leq m,
x^{j-1}e^{-x^2/2}, 1\leq j \leq n\},
\end{equation}
and similarly a set of functions $\{\psi_k^{(n)}(x),k=1,\ldots,n+m\}$
spanning $V(b_1,\ldots,\break b_m)$. Then,
%
%
\begin{equation}\qquad
P\bigl(\mathbf{x}(0) = \mathbf{x}\bigr)=\frac{1}{Z''} \det\bigl(\varphi
_i^{(n)}(x_j) \bigr)_{1\leq i,j\leq n+m} \det\bigl(\psi
_i^{(n)}(x_j) \bigr)_{1\leq i,j\leq n+m}.
\end{equation}
As mentioned above, this measure defines a determinantal point process
with defining kernel
%
%
\begin{equation}
K(x,y)=\sum_{i,j=1}^{n+m} \psi_i^{(n)}(x) [B^{-1}]_{i,j}\varphi_j^{(n)}(y),
\end{equation}
where $B=[B_{i,j}]_{1\leq i,j\leq n+m}$ has entries $B_{i,j}=\langle
\varphi
_i^{(n)},\psi_j^{(n)}\rangle$. Thus, the goal is to find nice functions
$\psi_k^{(n)}$ and $\varphi_k^{(n)}$ such that the inverse of the matrix
$B$ is manageable; usually one looks for a set of functions such that
$B$ becomes the identity matrix (biorthogonalization). In this
instance, it is more convenient for doing asymptotics to find functions
such that the matrix $B$ has the form
%
%
\begin{equation}\label{eq2.13}
B=\pmatrix{\mu& 0 \cr0 & \Id_n}.
\end{equation}
As will be shown below, the choice of functions for which this is the
case is as follows:
%
%
\begin{eqnarray}\label{phi}\quad
\varphi_k^{(n)}(x)&=&\frac{e^{-x^2/2}}{2\pi i}\oint_{\Gamma
_{0,a_k/2}}dz\frac
{e^{-z^2+2xz}}{z^n(z-a_k/2)},\qquad 1\leq k \leq m,\nonumber\\
\varphi_{m+k}^{(n)}(x)&=&\frac{(k - 1)!}{c_{k-1}}\frac
{e^{-x^2/2}}{2\pi i}
\oint_{\Gamma_0}dz\frac{e^{-z^2+2xz}}{z^{k}}
=\frac{H_{k-1}(x)}{c_{k-1}}e^{-x^2/2},\\
\eqntext{1\leq k \leq n.}
\end{eqnarray}
The $H_k(x)$ are the classical Hermite polynomials, with generating function
%
\begin{equation}\label{H-integr}
e^{-z^2+2xz}=\sum^{\infty}_{j=0}\frac{z^j}{j!}H_j(x)
\quad\mbox{and thus}\quad \frac{1}{2\pi i}\oint_{\Gamma
_0}e^{-z^2+2xz}\,\frac{dz}{z^{j+1}}=\frac{H_{j}(x)}{j!},\hspace*{-33pt}
\end{equation}
and with orthogonality relations
%
\begin{equation}\label{Herm-ortho}
\int_{{\mathbb R}}dx\, H_k(x)H_{\ell}(x)e^{-x^2}=\delta_{k,\ell}
c_k^2\qquad\mbox{with }
c_k=\sqrt{2^{k}k!}\sqrt[4]{\pi}.
\end{equation}
By the residue theorem, it follows that
%
\begin{equation}\quad
\frac{e^{-x^2/2}}{2\pi i}\oint_{\Gamma_{0,a/2}}dz\frac
{e^{-z^2+2xz}}{z^n(z-a_k/2)}\in
\operatorname{span}(e^{a_kx},H_0,\ldots,H_{n-1})e^{-x^2/2}.
\end{equation}
Similarly, one defines the functions $\psi_k^{(n)}(x)$ upon replacing
$a_k$ by $b_k$ in (\ref{phi}). Thus, the set of functions $\{\varphi
_k^{(n)}(x),k=1,\ldots,n+m\}$ spans the vector space $ V(a_1,\ldots
,a_m)$, and the set $\{\psi_k^{(n)}(x),k=1,\ldots,n+m\}$ the vector
space $ V(b_1,\ldots, b_m)$, as defined in (\ref{vspace}).

The last step is to show that with our choice we actually obtain (\ref
{eq2.13}). From the representation (\ref{phi}) of the $\varphi
_k^{(n)}(x),\psi_{\ell}^{(n)}(y)$ in terms of Hermite polynomials, it
follows immediately that
%
\begin{equation}
\mu_{k\ell}=\bigl\langle\varphi_k^{(n)},\psi_{\ell}^{(n)}\bigr\rangle
=\delta_{k\ell} \quad\mbox{for } m+1\leq k,\ell\leq m+n.
\end{equation}
Next, we show that
\begin{eqnarray*}
\bigl\langle\varphi_k^{(n)},\psi_{\ell}^{(n)}\bigr\rangle&=&0
\qquad\mbox{for }1\leq k\leq m,m+1\leq\ell\leq m+n\\
&&\hspace*{24.77pt}\mbox{ and }m+1\leq k\leq m+n,1\leq\ell\leq m.
\end{eqnarray*}
Indeed for $1\leq k\leq m$ and $m+1\leq\ell\leq m+n$, we have
%
%
\begin{eqnarray}
\bigl\langle\varphi_k^{(n)},\psi_{\ell}^{(n)}\bigr\rangle&=&\frac{
\mathrm{const}}{(2\pi i)^2}
\oint_{\Gamma_{0,a/2}}\frac{dz\, e^{-z^2}}{z^n(z-a_k/2)}\oint_{\Gamma
_0}dw\frac{e^{-w^2}}{w^{\ell-m}} \int_{-\infty}^{\infty
}dx\, e^{-x^2+2x(w+z)} \nonumber\hspace*{-20pt}\\
&=&\frac{\operatorname{const}\sqrt{\pi}}{(2\pi i)^2} \oint_{\Gamma
_{0,a/2}}\frac
{dz}{z^n(z-a_k/2)}\oint_{\Gamma_0}\frac{dw}{w^{\ell-m}}e^{2zw}
\hspace*{-20pt}\\
&=&\frac{\operatorname{const}\sqrt{\pi}}{2\pi i} \oint_{\Gamma
_{0,a/2}}\frac
{dz\, P_{\ell-m-1}(z)}{z^n(z-a_k/2)}=0,\nonumber\hspace*{-20pt}
\end{eqnarray}
where $P_i(x)$ is a polynomial of degree $i$. The result is zero
because for $\ell-m-1\leq n-1$ the residue at infinity is zero.

Finally, for $1\leq k,\ell\leq m$, by the same argument one gets
%
%
\begin{equation}\label{eq2.21}
\mu_{k,\ell}=\bigl\langle\varphi_k^{(n)},\psi_{\ell}^{(n)}\bigr\rangle
=\frac{\sqrt{\pi}}{(2\pi i)^2}
\oint_{\Gamma_{0,a_k/2}}\frac{dz}{z^n(z-a_k/2)}\oint_{\Gamma
_{0,b_{\ell}/2}}
\frac{dw\, e^{2zw}}{w^n(w-b_{\ell}/2)}.\hspace*{-38pt}
\end{equation}
By the residue theorem, the contribution of the pole at $w=0$ is a
polynomial of degree $n-1$ in $z$. Thus, the integral over $z$ is zero,
because the residue at infinity is zero. Thus, it remains to compute
the contribution of the pole at $w=b_\ell/2$, namely
%
%
\begin{eqnarray}\label{2.25}
\mu_{k,\ell}&=&\frac{\sqrt{\pi}}{(2\pi i)^2}
\oint_{\Gamma_{0,a_k/2}}dz \frac{e^{z b_\ell
}}{z^n(z-a_k/2)}\frac{2^n}{b_\ell^n}
\nonumber\\[-8pt]\\[-8pt]
&=&\frac{\sqrt{\pi}2^n}{(2\pi i)^2}
\oint_{\Gamma_{0,a_k b_\ell/2}}dz \frac
{e^{z}}{z^n(z-a_kb_\ell/2)}.\nonumber
\end{eqnarray}
This ends the proof of Proposition \ref{Prop:1.1}.
\end{pf}

The next step in showing Theorem \ref{th:1} for $\tau=0$ is to
determine the $n\to\infty$ limit of the kernel under the space scaling
%
%
\begin{equation}\label{eq2.28}
x=\sqrt{2n}+\frac{\xi_1}{\sqrt{2}n^{1/6}},\qquad y=\sqrt{2n}+\frac{\xi
_2}{\sqrt{2}n^{1/6}}
\end{equation}
with $a_i, b_i$ scaled as in (\ref{scaling}),
%
\begin{equation}\label{2.ab-scaling}
a_i=\sqrt{2n} \biggl(1+\frac{\tilde a_i}{n^{1/3}} \biggr)
\quad\mbox{and}\quad
b_i=\sqrt{2n} \biggl(1-\frac{\tilde b_i}{n^{1/3}} \biggr)
\end{equation}
and with the assumption
%
%
\begin{equation}\label{Assumption}
\tilde a_i<\tilde b_j,\qquad 1\leq i,j\leq m.
\end{equation}
Thus, we have to show that for $\xi_1,\xi_2$ in a bounded set,
%
%
\begin{equation}\label{eqKK}
\lim_{n\to\infty} \frac{1}{\sqrt{2}n^{1/6}}K_{n,m}(x,y)
=K_m^{\tilde a,\tilde b}(0;\xi_1,\xi_2).
\end{equation}

It is well known that the Hermite kernel under the above scaling, for
$\xi_1,\xi_2$ in a bounded set, converges to the Airy kernel
$K_\mathcal
{A}$ (see, e.g., Appendix A.7 of \cite{Ferrari})
%
%
\begin{eqnarray}\label{eq2.31}
\lim_{n\to\infty} \frac{1}{\sqrt{2}n^{1/6}} K_n^{\mathrm{Hermite}}(x,y)
&=&\frac{1}{ (2\pi i )^2}\int_{\Gamma_>} d\omega\int_{\Gamma_<}
d\widetilde\omega\frac{1}{\omega-\widetilde\omega}\frac
{e^{-\omega
^3/3+\xi
_2\omega}}{e^{-\widetilde\omega^{3}/3+\xi_1\widetilde\omega
}}\hspace*{-38pt}\nonumber\\[-8pt]\\[-8pt]
&=&\!: K_\mathcal{A}(\xi_1,\xi_2),\hspace*{-38pt}\nonumber
\end{eqnarray}
where the path $\Gamma_>$ goes from $e^{-2\pi i/3}\infty$ to $e^{2\pi
i/3}\infty$, the path $\Gamma_<$ from $e^{\pi i/3}\infty$ to
$e^{-\pi
i/3}\infty$, with $\Gamma_>$ and $\Gamma_<$ not intersecting each other.

What remains is to compute the limit of the last term in (\ref
{finite_kernel}). Since $m$ remains finite, one can take the $n\to
\infty
$ limit inside the sum. Below we compute the asymptotics for $\psi
_i^{(n)}$, $\varphi_j^{(n)}$, and $\mu^{-1}_{i,j}$ separately.
Let us start with the matrix $\mu$, as defined in (\ref{eq2.13}).
\begin{lemma} \label{lemma_mubar}
The following asymptotics holds for the inverse of the $m\times m$ matrix:
%
%
\begin{equation}\qquad\quad
\lim_{n\to\infty} \frac{1}{\sqrt{2} n^{1/6}} \biggl(\frac{2e}{n} \biggr)^n
\mu^{-1}=-A^{-1}\qquad \mbox{where }
A= \biggl(\frac{1}{\tilde a_k-\tilde b_{\ell}} \biggr)_{1\leq k,\ell\leq m}.
\end{equation}
\end{lemma}
\begin{pf}
Using the scaling (\ref{scaling}), the quantity
%
\begin{equation}\label{ab}
\frac{a_kb_{\ell}}{2}=n \biggl(1+\frac{\tilde a_k-\tilde b_{\ell
}}{n^{1/3}}+\mathcal{O}\biggl(\frac{1}{n^{2/3}} \biggr) \biggr)
\end{equation}
is, for $n$ large enough, strictly less than $n$ by assumption (\ref
{Assumption}). We use (\ref{1.mukl}) and make the change of variable
$z=u n$
%
%
\begin{eqnarray}\label{eq2.35}\quad
\mu_{k\ell}&=&\frac{\sqrt{\pi} 2^n}{2\pi i}
\oint_{\Gamma_{0,{(a_kb_{\ell})}/{2}}}\frac{dz\, e^z}{z^n (z-
{a_kb_{\ell}}/{2} )} \nonumber\\[-8pt]\\[-8pt]
&=& \frac{\sqrt{\pi} (2/n)^n}{2\pi i}
\oint_{|u|=1}du \frac{e^{n F(u)}}{u-1-(\tilde a_k-\tilde b_\ell
)n^{-1/3}+\mathcal{O}(n^{-2/3})}, \nonumber
\end{eqnarray}
where
%
%
\begin{equation}
F(u):=u-\ln u=1+\tfrac12(u-1)^2+\mathcal{O}\bigl((u-1)^3 \bigr),
\end{equation}
with
%
\begin{equation}\label{eqdescent}
\operatorname{Re}(F(u))=\operatorname{Re}(u)-\ln(|u|).
\end{equation}

Thus, we can deform the path $|u|=1$ into $\gamma_\delta=\{
1+iy,-\delta
\leq y \leq\delta\}$ plus a circle segment $\gamma'$ centered at zero
joining the extremities of $\gamma_\delta$. By (\ref{eqdescent}), the
path $\gamma_\delta\vee\gamma'$ is a steepest descent path for $F$ with
maximum at $u=1$, $F(1)=1$. We choose $\delta=n^{-2/5}$, then, the
contribution of the integral in (\ref{eq2.35}) from $\gamma'$ is of
order $\mathcal{O}(e^{-c n^{1/5}})$ smaller than the main
contribution, coming
from $\gamma_\delta$, for some $c>0$. Thus, continuing (\ref{eq2.35}),
%
%
\begin{eqnarray}\qquad
\mu_{k\ell} &=& \biggl(\frac{2e}{n} \biggr)^{ n} \frac{\sqrt{\pi}}{2\pi i}
\int_{1-i n^{-2/5}}^{1+i n^{-2/5}} du \frac{e^{n(u-1)^2/2+n\mathcal{O}
((u-1)^3 )}}{u - 1-(\tilde a_k - \tilde b_\ell)n^{-
{1}/{3}}+\mathcal{O}(n^{-{2}/{3}})}\\
&&\hspace*{68.5pt}\hspace*{40.6pt}{}\times \bigl(1 + \mathcal{O}(e^{-c
n^{1/5}}) \bigr).\nonumber
\end{eqnarray}
By the change of variable $\omega=(u-1)\sqrt{n}$, the last integral becomes
%
%
\begin{equation}\label{eq2.39}
\frac{\sqrt{\pi}}{2\pi i}
\int_{-i n^{1/10}}^{i n^{1/10}} d\omega\frac{e^{(1/2)\omega
^2(1+\mathcal{O}
(n^{-2/5}))}}{\omega-(\tilde a_k-\tilde b_\ell)n^{1/6}+\mathcal
{O}(n^{-1/6})}.
\end{equation}
In the $n\to\infty$ limit, we finally have
%
%
\begin{equation}
\lim_{n\to\infty}\sqrt{2}n^{1/6}(\ref{eq2.39})=\frac{-1}{\tilde
a_k-\tilde b_\ell}.
\end{equation}
Thus, we have shown that
%
%
\begin{equation}
\lim_{n\to\infty}\sqrt{2}n^{1/6} \biggl(\frac{n}{2e} \biggr)^n \mu
_{k,l}=\frac{-1}{\tilde a_k-\tilde b_\ell}=-A_{k,l}.
\end{equation}
This suffices to prove Lemma \ref{lemma_mubar}, since the dimension of
the matrix does not depend on $n$.
\end{pf}

The next item is to determine the asymptotics of $\varphi_k^{(n)}$ and
$\psi
_k^{(n)}$.
\begin{lemma}\label{Lemma1.4}
Consider the scaling (\ref{eq2.28}) and (\ref{scaling}), with $\xi
_1,\xi
_2$ in a bounded set. Then
%
%
\begin{eqnarray}\label{eq2.42}
\varphi_k(\xi_2):\!&=& \lim_{n\to\infty} \biggl(\frac{n}{2e} \biggr)^{n/2}\varphi
_k^{(n)} \biggl(\sqrt{2n}+\frac{\xi_2}{\sqrt{2} n^{1/6}}
\biggr)\nonumber\\[-8pt]\\[-8pt]
&=&\frac{1}{2\pi i}\int_{\Gamma_{\tilde a_k>}}d\omega\frac
{e^{-\omega
^3/3+\xi
_2\omega}}{\omega-\tilde a_k},\nonumber
\end{eqnarray}
where $\Gamma_{\tilde a_k>}$ is a simple path from $e^{-2\pi
i/3}\infty
$ to $e^{2\pi i/3}\infty$ and passing onto the right of $\tilde a_k$.
Similarly,
%
%
\begin{eqnarray}
\psi_k(\xi_1):\!&=&\lim_{n\to\infty} \biggl(\frac{n}{2e} \biggr)^{n/2}\psi
_k^{(n)} \biggl(\sqrt{2n}+\frac{\xi_1}{\sqrt{2} n^{1/6}}
\biggr)\nonumber\\[-8pt]\\[-8pt]
&=&\frac{1}{2\pi i}\int_{\Gamma_{<\tilde b_k}}d\widetilde\omega
\frac
{e^{\widetilde\omega^{3}/3-\xi_1\widetilde\omega}}{\widetilde
\omega
-\tilde b_k},\nonumber
\end{eqnarray}
where $\Gamma_{<\tilde b_k}$ is a simple path from $e^{\pi i/3}\infty$
to $e^{-\pi i/3}\infty$ and passing onto the left of $\tilde b_k$
(similar to Figure \ref{Fig1}).
\end{lemma}

\begin{pf}
The plan is to compute the large $n$ behavior of
%
\begin{equation}
\varphi_k^{(n)}(x)=\frac{e^{-x^2/2}}{2\pi i}\oint_{\Gamma
_{0,a_k/2}}dz\frac
{e^{-z^2+2xz}}{z^n(z-a_k/2)},
\end{equation}
with
%
\begin{equation}
a_k=\sqrt{2n} \biggl(1+\frac{\tilde a_k}{n^{1/3}} \biggr),\qquad x=\sqrt
{2n}+\frac{\xi}{\sqrt{2} n^{1/6}}.
\end{equation}
Rescaling the integration variable $z=u\sqrt{n/2}$, one gets
%
\begin{equation}\hspace*{28pt}
\varphi_k^{(n)}(x) = \biggl(\frac{2}{n} \biggr)^{n/2} \frac{e^{-n-\xi
n^{1/3}+\mathcal{O}(n^{-1/3})}}{2\pi i} \oint_{\Gamma_{0,1+\tilde
a_k/n^{1/3}}}
du \frac{e^{n F(u)+u\xi n^{1/3}}}{u-1-\tilde a_k/n^{1/3}},
\end{equation}
where $F(u)=-u^2/2+2u-\ln(u)$. The leading contribution comes from the
%
\begin{figure}

\includegraphics{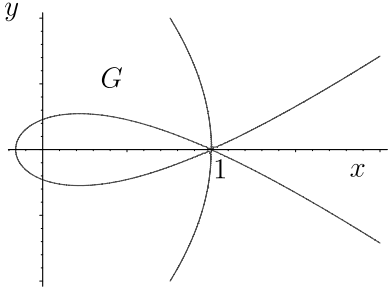}

\caption{Lines with $\operatorname{Re}(F(x+iy))=\operatorname{Re}(F(1))$.}
\label{FigSteepPhi}
\end{figure}
neighborhood of the double critical point of $F(u)$ at $u=1$, where we have
%
\begin{equation}
F(u)=\tfrac32-\tfrac13 (u-1)^3+\mathcal{O}\bigl((u-1)^4\bigr).
\end{equation}
As integration path one can choose any path passing through $u=1+\alpha
n^{-1/3}$, with $\tilde a_k<\alpha$, locally following the directions
$e^{\pm2\pi i/3}$, and which remain inside the region $G$ of
Figure \ref{FigSteepPhi}.
Then the integration away from a $\delta$-neighborhood of $u=1+\alpha
n^{-1/3}$ (where $\delta=n^{-\varepsilon}$, with $0<\varepsilon<1/3$)
will be of order
$\mathcal{O}
(e^{-c n})$ smaller than the leading term, with $0<c\sim\delta^3$ for
small $\delta$. Then, in a $\delta$-neighborhood of $u=1$, one can use
series expansions and after the change of variable $\omega
=n^{1/3}(u-1)$, one finds
\begin{eqnarray*}
{\varphi_k^{(n)} \biggl(\sqrt{2n}+\frac{\xi}{\sqrt{2} n^{1/6}}
\biggr)} &=& \biggl( \frac{2e}{n} \biggr)^{n/2} \biggl[ \frac{1}{2\pi i}
\biggl(1+\mathcal{O}\biggl(\frac{1}{n^{1/3}} \biggr) \biggr)\int
d\omega\frac{e^{-\omega
^3/3+\xi\omega}}{\omega-\tilde a_k} \biggr] \\
&&{}\times\bigl(1+\mathcal{O}(e^{-c n}) \bigr),
\end{eqnarray*}
where the integral goes from $e^{-2\pi i/3}\delta n^{1/3}$ to $e^{2\pi
i/3}\delta n^{1/3}$, and passing to the right of $\tilde a_k$. From
this, the $n\to\infty$ limit in (\ref{eq2.42}) holds.

The asymptotic for $\psi_k^{(n)}(x)$ is essentially the same, except that
\mbox{$\omega\mapsto-\omega$} and $\tilde a_k\mapsto-\tilde b_k$,
ending the proof of Lemma \ref{Lemma1.4}.
\end{pf}

We shall also need the following lemma.
\begin{lemma}\label{LemmaIdentity}
Given the matrix
%
\begin{equation}
A:= \biggl( \biggl(\frac{1}{\tilde a_i-\tilde b_j} \biggr)_{1\leq i,j\leq
m} \biggr),
\end{equation}
the following identity holds:
%
\begin{equation}\label{identity}\qquad
\sum_{1\leq i,j\leq m}\frac{(A^{\top-1})_{ij}}{(z-\tilde
a_i)(w-\tilde
b_j)} =\frac{1}{w-z} \Biggl(\prod_{k=1}^m \biggl(\frac{w-\tilde
a_k}{z-\tilde a_k} \biggr)
\biggl(\frac{z-\tilde b_k}{w-\tilde b_k} \biggr)-1 \Biggr).
\end{equation}
\end{lemma}
\begin{pf}
Since
%
\begin{equation}\label{det}
\det A=\frac{\Delta(\tilde a)\Delta(\tilde b)}{ {\prod_{1\leq
i,j\leq m}}(\tilde a_i-\tilde b_j)},
\end{equation}
one checks\vspace*{1pt} the identity (\ref{identity}), by computing the residue on
the right-hand
side at the points $z=\tilde a_i$, $w=\tilde b_j$ and identifying with
$(A^{\top-1})_{ij}$ using Cram\'{e}r's
rule and repeatedly using (\ref{det}).
\end{pf}
\begin{pf*}{Proof of Theorem \protect\ref{th:1}}
Assembling the asymptotic result
(\ref{eq2.31}), Proposition \ref{Prop:1.1}, Lemmas \ref{lemma_mubar} and
\ref{Lemma1.4}, one obtains Theorem \ref{th:1} in the special
case $\tau=0$, with distinct $\tilde a_i$, $\tilde b_i$, under the
condition $\tilde a_i<\tilde b_j$. Upon using the scaling (\ref
{eq2.28}) and (\ref{2.ab-scaling}), the limit kernel is thus given by
the limit of the sum of the kernels in (\ref{finite_kernel}); i.e., the
sum of the Airy kernel $K_\mathcal{A}$, defined in (\ref{eq2.31}), and
a new kernel:\looseness=1
%
%
\begin{eqnarray}\label{eq2.49}\qquad\quad
&&\lim_{n\to\infty} \frac{1}{\sqrt{2}n^{1/6}}K_{n,m}(x,y)
\nonumber\\
&&\qquad= K_\mathcal{A}(\xi_1,\xi_2)-\sum_{i,j=1}^m \psi_i(\xi_1)
[A^{-1}]_{i,j} \varphi_j(\xi_2)\nonumber\\[-8pt]\\[-8pt]
&&\qquad=\frac{1}{ (2\pi i )^2}\int_{\Gamma_>} d\omega\int
_{\Gamma_<} d\widetilde\omega\frac{e^{-\omega^3/3+\xi_2\omega
}}{e^{-\widetilde\omega^{3}/3+\xi_1\widetilde\omega}} \frac
{1}{\omega
-\widetilde\omega} \nonumber\\
&&\qquad\quad{}-\frac{1}{ (2\pi i )^2}\int_{\Gamma_{\tilde a>}}
d\omega\int_{\Gamma_{<\tilde b}} d\widetilde\omega\frac
{e^{-\omega
^3/3+\xi_2\omega}}{e^{-\widetilde\omega^{3}/3+\xi_1\widetilde
\omega}}
\sum_{i,j=1}^m \frac{[A^{-1}]_{i,j}}{(\widetilde\omega-\tilde
b_i)(\omega-\tilde a_j)}. \nonumber
\end{eqnarray}
The fact that this expression actually equals the kernel $K_m^{\tilde
a,\tilde b}(0;\xi_1,\xi_2)$, as defined in Theorem \ref{th:1}, follows
from Lemma \ref{LemmaIdentity}.
\end{pf*}

\section{Extended kernel for the Airy process with wanderers}\label{SectThm2}

In this section, we will prove Theorem \ref{th:2}. For this purpose, we
need to know the measure, defined on the positions of the Brownian
bridges at different times $-1<T_1<T_2<\cdots<T_\ell<1$. Set
$\mathbf{x}(T_i):=(x_1(T_i),\ldots,x_{m+n}(T_i))$. Then, by Karlin--McGregor
applied to these different times, the measure obtained by the
nonintersecting condition on the Brownian bridges is given by
%
%
\begin{eqnarray}\label{3.KM}
&&P\bigl(\mathbf{x}(T_1)=\mathbf{x}^1,\ldots,\mathbf{x}(T_\ell)=\mathbf{x}^\ell
\bigr)\nonumber\\
&&\qquad=\frac{1}{Z}\det\bigl(p(a_i,x^1_j,T_1+1)\bigr)_{1\leq i,j\leq n+m}
\nonumber\\[-8pt]\\[-8pt]
&&\qquad\quad{}\times
\Biggl(\prod_{k=1}^{\ell-1} \det\bigl(p(x_i^k,x^{k+1}_j,T_{k+1}-T_k)\bigr)_{1\leq
i,j\leq n+m} \Biggr)\nonumber\\
&&
\qquad\quad{}\times
\det\bigl(p(x_i^\ell,b_j,1-T_\ell)\bigr)_{1\leq i,j\leq n+m}.\nonumber
\end{eqnarray}
It is well known that this measure, a generalization of (\ref{eq2.1})
to multi-times, or any measure of this form has determinantal
correlations in space--time \cite{EM97,TW98,FN98,Jo03b,RB04} (even in
cases when the size of the determinant is increasing \cite{BFPS06,BF07}).\looseness=1
\begin{proposition}\label{prop:det-corr} Any measure on $\{
x_i^{(n)},1\leq i \leq N,1\leq n \leq\ell\}$ of the form\footnote{The
functions $\phi(T_n,x ;T_{n+1},y)$ themselves may in fact vary with
$n$ above.}
%
%
\begin{eqnarray}\label{eq3.1}
&&\frac{1}{Z}\det\bigl(\phi\bigl(T_0,a_i;T_1,x^{(1)}_j\bigr)\bigr)_{1\leq i,j\leq N}
\nonumber\\
&&\qquad
{}\times\Biggl(\prod_{n=1}^{\ell-1}\det\bigl(\phi
\bigl(T_n,x^{(n)}_i;T_{n+1},x^{(n+1)}_j\bigr)\bigr)_{1\leq i,j\leq N}
\Biggr)\hspace*{-12pt}\\
&&\qquad{}\times\det\bigl(\phi\bigl(T_\ell,x^{(\ell)}_i;T_{\ell+1},b_j\bigr)\bigr)_{1\leq
i,j\leq N},\hspace*{-12pt}\nonumber
\end{eqnarray}
has, assuming $Z\neq0$, the following $k$-point correlation functions
for $t_1,\ldots,t_k\in\{T_1,\ldots,T_\ell\}$:
%
%
\begin{equation}
\rho^{(k)}(t_1,x_1,\ldots,t_k,x_k)=\det(K(t_i,x_i;t_j,x_j)
)_{1\leq i,j\leq k},
\end{equation}
where the space--time kernel (often called extended kernel) is given by
%
%
\begin{eqnarray}\label{eqExtKernelGeneral}
&&K(t_1,x_1;t_2,x_2)\nonumber\\[-2pt]
&&\qquad=-\phi(t_1,x_1;t_2,x_2)\Id(t_2>t_1)\\[-2pt]
&&\qquad\quad{} + \sum_{i,j=1}^N\phi(t_1,x_1;T_{\ell+1},b_i) [B^{-1}]_{i,j} \phi
(T_0,a_j;t_2,x_2) \nonumber
\end{eqnarray}
with ($*$ means integration with regard to the consecutive dots)
%
%
\begin{eqnarray}\label{3.5}
&&\phi(T_r,x;T_s,y)
\nonumber\\[-10pt]\\[-10pt]
&&\qquad=
\cases{\phi(T_r,x;T_{r+1},\cdot)*\cdots*\phi(T_{s-1},\cdot
;T_s,y), &\quad
if $T_r<T_s$,\cr
0,&\quad if $T_r\geq T_s$,}\nonumber
\end{eqnarray}
and with the $N\times N$ matrix $B$ having entries $B_{i,j}=\phi
(T_0,a_i;T_{\ell+1},b_j)$. Remark that $(N!)^\ell\det(B)=Z$, so that
$B^{-1}$ exists as soon as $Z\neq0$.
\end{proposition}

We now apply this general fact to the nonintersecting Brownian motion
formula (\ref{3.KM}): here $x_i^{(n)}$ denotes the position $x_i(T_n)$
of the $i$th Brownian motion at time $T_n$, while one sets $T_0=-1$,
$T_{\ell+1}=1$, and one sets
%
\begin{equation}\label{3.def1}
\phi(t,x;t',x'):=p(x,x',t'-t).
\end{equation}

As for the one-time case, the structure is unchanged, even after
letting $a_{n+m},\ldots,a_{m+1} \to0$ and $b_{n+m},\ldots,b_{m+1}
\to0$, keeping $a_m<\cdots<a_1$ and $b_m<\cdots<b_1$ fixed. The only
difference is that the entries on the first and last determinants in
(\ref{eq3.1}) will be different (together with a different
normalization constant $Z$). Indeed, the first determinant in (\ref
{eq3.1}) is just replaced by
%
%
\begin{equation}\label{eq3.6}
\det\pmatrix{\bigl(e^{a_i x_j^{(1)}/(1+T_1)}p\bigl(0,x_j^{(1)},T_1+1\bigr)\bigr)
_{\fontsize{8.36}{10}\selectfont{\begin{array}{l}
1\leq i\leq m\\ 1\leq j\leq m+n
\end{array}}} \cr
\biggl(\biggl(\dfrac{x_j^{(1)}}{1+T_1}\biggr)^{i-1}p\bigl(0,x_j^{(1)},T_1+1\bigr)\biggr)
_{\fontsize{8.36}{10}\selectfont{\begin{array}{l}
1\leq i\leq n\\ 1\leq j\leq m+n
\end{array}}}},
\end{equation}
while the last determinant is replaced by
%
%
\begin{equation}\label{eq3.7}
\det\pmatrix{
\bigl(e^{b_i x_j^{(\ell)}/(1-T_\ell)}p\bigl(x_j^{(\ell)},0,1-T_\ell\bigr) \bigr)_{
{\fontsize{8.36}{10}\selectfont{\begin{array}{l} 1\leq i\leq m\\ 1\leq j\leq m+n
\end{array}}}}
\cr \biggl(\biggl(\dfrac{x_j^{(\ell)}}{1-T_\ell}\biggr)^{i-1}p\bigl(x_j^{(\ell
)},0,1-T_\ell\bigr) \biggr)_{\fontsize{8.36}{10}\selectfont{\begin{array}{l}
1\leq i\leq n\\ 1\leq j\leq m+n
\end{array}}}}.
\end{equation}

As for the one-time situation, one looks for sets of functions
generating the same vector spaces as the functions in (\ref{eq3.6}) and
(\ref{eq3.7}), namely one searches for functions $\varphi^{(n)}_k(T_1,x)$
and $\psi^{(n)}_k(T_\ell,x)$, such that
%
%
\begin{eqnarray}
(\ref{eq3.6})&=&\operatorname{const} \times\det\bigl(\varphi
^{(n)}_i(T_1,x_j^1)\bigr)_{1\leq
i,j\leq n+m}, \nonumber\\[-10pt]\\[-10pt]
(\ref{eq3.7})&=&\operatorname{const} \times\det\bigl(\psi^{(n)}_i(T_\ell
,x_j^\ell
)\bigr)_{1\leq i,j\leq n+m},\nonumber
\end{eqnarray}
and such that the matrix $B$ has the same form (\ref{eq2.13}) as
before. Setting
%
%
\begin{equation}
\gamma(t):=\sqrt{\frac{1-t}{1+t}},\qquad \sigma(t):=\sqrt{1+t},
\end{equation}
one picks, for $1\leq k\leq m$,
%
%
\begin{eqnarray}\label{eqphipsim}
\varphi_k^{(n)}(t,x) &:=& \frac{e^{-x^2/2\sigma(t)^2}}{\sigma(t)}
\frac
{1}{2\pi i}\oint_{\Gamma_{0,a_k/2}} dz\frac{e^{-z^2\gamma
(t)^2+2xz/\sigma^2(t)}}{z^n (z-a_k/2 )},\nonumber\\[-8pt]\\[-8pt]
\psi_k^{(n)}(t,x) &:=& \frac{e^{-x^2/2\sigma(-t)^2}}{\sigma(-t)}
\frac
{1}{2\pi i}\oint_{\Gamma_{0,b_k/2}} dz\frac{e^{-z^2\gamma
(-t)^2+2xz/\sigma^2(-t)}}{z^n (z-b_k/2 )},\nonumber
\end{eqnarray}
and for $1\leq k \leq n$,
%
%
\begin{eqnarray}\label{eqphipsin}\hspace*{33pt}
\varphi_{m+k}^{(n)}(t,x) &:=& \frac{(4\pi)^{1/4}}{\sqrt
{(k-1)!2^{k-1}}}\gamma
(t)^{k-1}
\nonumber\\
&&{}\times H_{k-1}\biggl(\frac{x}{\gamma(t)\sigma^2(t)} \biggr) p(0,x;t+1),
\nonumber\\[-8pt]\\[-8pt]
\psi_{m+k}^{(n)}(t,x) &:=& \frac{(4\pi)^{1/4}}{\sqrt
{(k-1)!2^{k-1}}}\gamma(-t)^{k-1}
\nonumber\\
&&{}\times H_{k-1}\biggl(\frac{x}{\gamma(-t)\sigma^2(-t)} \biggr) p(x,0;1-t).\nonumber
\end{eqnarray}
Remark that, using the integral representation of the Hermite
polynomials, an equivalent expression for (\ref{eqphipsin}) is
%
%
\begin{eqnarray}\label{eqphipsinB}
\varphi_{m+k}^{(n)}(t,x) &=& \frac{(4\pi)^{1/4}}{\sqrt{(k-1)!2^{k-1}}}
\frac
{e^{-x^2/2\sigma(t)^2}}{\sigma(t)}\nonumber\\
&&{}\times \frac{(k-1)!}{2\pi i}\oint
_{\Gamma
_{0}}
dz\frac{e^{-z^2\gamma(t)^2+2xz/\sigma^2(t)}}{z^k},\nonumber\\[-8pt]\\[-8pt]
\psi_{m+k}^{(n)}(t,x) &=& \frac{(4\pi)^{1/4}}{\sqrt{(k-1)!2^{k-1}}}
\frac
{e^{-x^2/2\sigma(-t)^2}}{\sigma(-t)}\nonumber\\
&&{}\times \frac{(k-1)!}{2\pi i}\oint
_{\Gamma
_{0}} dz\frac{e^{-z^2\gamma(-t)^2+2xz/\sigma^2(-t)}}{z^k}.\nonumber
\end{eqnarray}

It is immediate to verify that these functions generate at $t=T_1$,
resp., $t=T_\ell$, the same space as the function in (\ref{eq3.6}),
resp., (\ref{eq3.7}). So, one defines the functions appearing in the
first and last determinant of (\ref{eq3.1}) by
%
\begin{eqnarray}\label{3.def2}
\phi\bigl(T_0,a_i;T_1,x^{(1)}\bigr)&:=&\varphi_i^{(n)}\bigl(T_1,x^{(1)}\bigr)
\quad\mbox{and}\nonumber\\[-8pt]\\[-8pt]
\phi\bigl(T_\ell,x^{(\ell)};T_{\ell+1},b_i\bigr)&:=&\psi_i^{(n)}\bigl(T_\ell
,x^{(\ell)}\bigr),\nonumber
\end{eqnarray}
for which we show the following property:
\begin{lemma}\label{tr-pr}
For any $t_1<t_2$ and $1\leq k\leq n+m$, one has
%
%
\begin{eqnarray}
\int_{{\mathbb R}} dx\, \varphi_k^{(n)}(t_1,x) p(x,y;t_2-t_1)&=& \varphi
_k^{(n)}(t_2,y),\nonumber\\[-8pt]\\[-8pt]
\int_{{\mathbb R}} dy\, p(x,y;t_2-t_1) \psi_k^{(n)}(t_2,y)&=& \psi
_k^{(n)}(t_1,x).\nonumber
\end{eqnarray}
\end{lemma}
\begin{pf}
Since $\psi_k^{(n)}$ is obtained from $\varphi_k^{(n)}$ by the map
$t\mapsto-t$ and $a\mapsto b$, it suffices to present the proof for
$\varphi_k^{(n)}$.
At first, for $1\leq k\leq m$, one has
%
%
\begin{eqnarray}\label{eq3.16}
&&\int_{{\mathbb R}} dx\, \varphi_k^{(n)}(t_1,x) p(x,y;t_2-t_1)
\nonumber\\
&&\qquad
=\frac{1}{2\pi i} \oint_{\Gamma_{0,a_k/2}} dz\frac{e^{-z^2\gamma
(t_1)^2}}{z^n (z-a_k/2 )} \\
&&\qquad\quad{}\times\int_{{\mathbb R}} dx \frac{e^{-x^2/2(1+t_1)}}{\sqrt
{1+t_1}}e^{2xz/(1+t_1)}\frac{e^{-(x-y)^2/2(t_2-t_1)}
}{\sqrt{2\pi(t_2-t_1)}}\nonumber
\end{eqnarray}
and, after performing the Gaussian integration, one has
%
%
\begin{equation}\qquad\quad
\mbox{(\ref{eq3.16})}=\frac{e^{-y^2/2\sigma(t_2)^2}}{\sigma(t_2)} \frac
{1}{2\pi i}\oint_{\Gamma_{0,a_k/2}} dz\frac{e^{-z^2\gamma
(t_2)^2+2yz/\sigma^2(t_2)}}{z^n (z-a_k/2 )}=\varphi_k^{(n)}(t_2,y).
\end{equation}
Second, consider $1\leq k \leq n$. Comparing the representations (\ref
{eqphipsim}) and (\ref{eqphipsinB}), we see immediately that the
computations are exactly the same. Indeed, the only difference is a
$k$-dependent prefactor and the denominator in the integrand over $z$.
However, there are not affected by the computations above; thus
%
%
\begin{equation}\label{eq3.16a}
\int_{{\mathbb R}} dx\, \varphi_{m+k}^{(n)}(t_1,x) p(x,y;t_2-t_1)
=\varphi_{m+k}^{(n)}(t_2,y)
\end{equation}
holds, ending the proof of Lemma \ref{tr-pr}.
\end{pf}
\begin{proposition} \label{prop:Kext}The extended kernel is given by
%
%
\begin{eqnarray}\label{2.Kext}\quad
&&K_{n,m}(t_1,x_1;t_2,x_2)\nonumber\\
&&\qquad=-p(x_1,x_2; t_2-t_1)\Id(t_2>t_1)
+\sum_{i=1}^{n} \psi^{(n)}_{m+i}(t_1,x_1)\varphi
^{(n)}_{m+i}(t_2,x_2)\\
&&\qquad\quad{} + \sum_{i,j=1}^{m} \psi^{(n)}_i(t_1,x_1)[\mu^{-1}]_{i,j} \varphi
^{(n)}_j(t_2,x_2)\nonumber
\end{eqnarray}
with $\varphi^{(n)}_j(t,x)$ and $\psi^{(n)}_{m+i}(t,x)$ given by
(\ref
{eqphipsim}) and (\ref{eqphipsinB}) and with $\mu$ given by (\ref
{eq2.21}), the same as in the $1$-time case.
\end{proposition}
\begin{pf}
Given the definitions (\ref{3.def1}) and (\ref{3.5}), the first
term in the kernel (\ref{eqExtKernelGeneral}) is simply
%
%
\begin{equation}
-\phi(t_1,x_1;t_2,x_2)\Id(t_2>t_1)=-p(x_1,x_2,t_2-t_1)\Id(t_2>t_1).
\end{equation}
It remains to be shown that $B$ has the form
%
%
\begin{equation}\label{eq2.13a}
B=\pmatrix{\mu& 0 \cr0 & \Id_n}
\end{equation}
as in the 1-time case, with $\mu$ given in (\ref{eq2.21}).

Indeed, for any choice of $1\leq k\leq\ell-2$, and for $t_i=T_i$ with
$1\leq i\leq\ell$, one has, using the convolution property of the
Brownian transition probability and the convolution property in
Lemma \ref{tr-pr}, the property that ($*$ means integration with regard
to the common variable)
\begin{eqnarray*}
B_{i,j}&=&\phi(T_0,a_i;T_{\ell+1},b_j) \\
&=&\varphi_i^{(n)}\bigl(t_1,x^{(1)}\bigr)\ast p\bigl(x^{(1)},x^{(2)};t_2-t_1\bigr)\ast
\cdots\\
&&{}\ast
p\bigl(x^{(k)},x^{(k+1)};t_{k+1}-t_{k}\bigr)\ast
p\bigl(x^{(k+1)},x^{(k+2)};t_{k+2}-t_{k+1}\bigr)\ast\cdots
\\
&&{}
\ast p\bigl(x^{(\ell-1)},x^{(\ell)};t_{\ell}-t_{\ell-1}\bigr)
\ast\psi_j^{(n)}\bigl(t_\ell,x^{(\ell)}\bigr)\\
&=& \bigl(\varphi_i^{(n)}\bigl(t_1,x^{(1)}\bigr)\ast
p\bigl(x^{(1)},x^{(k+1)};t_{k+1}-t_1\bigr) \bigr)\\
&&{}\ast
\bigl(p\bigl(x^{(k+1)},x^{(\ell)};t_\ell-t_{k+1}\bigr)\ast\psi_j^{(n)}\bigl(t_\ell
,x^{(\ell) }\bigr) \bigr)\\
&=&\varphi_i^{(n)}\bigl(t_{k+1},x^{(k+1)}\bigr)\ast\psi_j^{(n)}\bigl(t_{k+1},x^{(k+1)}\bigr)
= \bigl\langle\varphi_i^{(n)}(t_{k+1},\cdot),
\psi_j^{(n)}(t_{k+1},\cdot)\bigr\rangle
\end{eqnarray*}
is independent of $t_{k+1}$; therefore, by setting $t_{k+1}=0$, it is,
in particular, equal to the value $\mu_{ij}$ obtained in (\ref{2.25})
and (\ref{eq2.21}). This establishes Proposition \ref{prop:Kext}.
\end{pf}

In order to prove Theorem \ref{th:2} (and thus also Theorem \ref{th:1}
for generic $\tau$), one needs to compute the $n\to\infty$ asymptotics
of the kernel. For convenience, recall the scaling for the starting and
ending points of the top $m$ Brownian bridges (\ref{scaling}) and of
the subsequent scaling (\ref{scalingXT}) of the space--time region one
focuses on
%
\begin{eqnarray}\label{scaling3}
a_i &=& \sqrt{2n}+\sqrt{2} \tilde a_i n^{1/6},\qquad b_i = \sqrt{2n}-\sqrt
{2} \tilde b_i n^{1/6},\nonumber\\[-8pt]\\[-8pt]
t_i &=& \tau_i n^{-1/3},\qquad {x_i= \sqrt{2n} +\frac{\xi_i-\tau
_i^2}{\sqrt{2} n^{1/6}}}\nonumber
\end{eqnarray}
with $\tilde a_i<\tilde b_j$, $1\leq i,j\leq m$. Below we prove that,
given the scaling (\ref{scaling3}) and for $\xi_1,\xi_2$ in a
bounded set,
%
%
\begin{equation}
\lim_{n\to\infty} \frac{1}{\sqrt{2}n^{1/6}}K_{n,m}(t_1,x_1;t_2,x_2)
\equiv K_m^{\tilde a,\tilde b}(\tau_1,\xi_1;\tau_2,\xi_2),
\end{equation}
where $\equiv$ we means an equivalent kernel.\footnote{Two kernels are
equivalent if they define the same determinantal point process. Namely,
if there exists some function $f(x)\neq0$ such that $K(x,y)=\widetilde
K(x,y) f(x)/f(y)$, then $K$ and $\widetilde K$ are equivalent, since
all the correlation functions are given by determinants in which the
functions $f$ cancel exactly.}
\begin{proposition}\label{Prop:kernel-ext}
With the\vspace*{2pt} above scaling, for $\xi_1,\xi_2$ in a bounded set (and $\tau
_1,\tau_2$ fixed), in the case where all the $\tilde a_i$ (and $\tilde
b_i$) are distinct, one has
%
%
\begin{equation}\qquad
\lim_{n\to\infty} \frac{1}{\sqrt
{2}n^{1/6}}K_{n,m}(t_1,x_1;t_2,x_2) =
K_m^{\tilde a,\tilde b}(\tau_1,\xi_1;\tau_2,\xi_2) \frac{f(\tau
_1,\xi
_1)}{f(\tau_2,\xi_2)},
\end{equation}
where $f(\tau,\xi)=\exp(\tau^3/3-\xi\tau)$.
\end{proposition}
\begin{pf}
Consider the first two terms in the kernel (\ref{2.Kext}). These
terms are independent of the $\tilde a_i,\tilde b_i$ and of $m$.
Indeed, it corresponds exactly to the kernel of the system without
wanderers, which can be denoted by $K_{n,0}$. Indeed,
%
%
\begin{eqnarray}\quad
&&\sum_{i=1}^{n} \psi^{(n)}_{m+i}(t_1,x_1)\varphi^{(n)}_{m+i}(t_2,x_2)
\nonumber\\
&&\qquad=
\sum
_{k=0}^{n-1} \frac{(4\pi)^{1/2}}{k!2^k}\gamma(t_2)^k\gamma(-t_1)^k
p(x_1,0;1-t_1) p(0,x_2;t_2+1)\\
&&\qquad\quad\hspace*{14.1pt}{}\times
H_{k} \biggl(\frac{x_1}{\gamma(-t_1)\sigma^2(-t_1)} \biggr)
H_{k} \biggl(\frac{x_2}{\gamma( t_2)\sigma^2( t_2)} \biggr).\nonumber
\end{eqnarray}
For fixed $\tau_1,\tau_2$, we show below that
%
%
\begin{equation}\label{eq3.26}
\lim_{n\to\infty} \frac{1}{\sqrt
{2}n^{1/6}}K_{n,0}(t_1,x_1;t_2,x_2) =
K_\mathcal{A}(\tau_1,\xi_1;\tau_2,\xi_2)\frac{f(\tau_1,\xi
_1)}{f(\tau
_2,\xi_2)},
\end{equation}
uniformly for $\xi_1,\xi_2$ in a bounded set, with $K_\mathcal{A}$ the
extended Airy kernel given by
%
%
\begin{eqnarray}\label{ExtAiry}
&&K_\mathcal{A}(\tau_1,\xi_1;\tau_2,\xi_2)\nonumber\\[-8pt]\\[-8pt]
&&\qquad =\cases{
\displaystyle\int_{{\mathbb R}_+} d\lambda\, e^{\lambda(\tau_2-\tau
_1)}\operatorname{Ai}(\xi
_1+\lambda)\operatorname{Ai}
(\xi_2+\lambda), &\quad $\tau_1\geq\tau_2$,\cr
\displaystyle-\int_{{\mathbb R}_-} d\lambda\, e^{\lambda(\tau_2-\tau
_1)}\operatorname{Ai}(\xi
_1+\lambda
)\operatorname{Ai}(\xi_2+\lambda), &\quad $\tau_1<\tau_2$.}\nonumber
\end{eqnarray}
To obtain this result, for $\xi_1,\xi_2$ in a bounded set, one can just
use the asymptotics of the classical Hermite polynomials (see, e.g., Appendix 7 of \cite{Ferrari}). Another, better, way is to first
perform the sum over $k$ using two different integral representations
for Hermite polynomials, a first one is (\ref{eqphipsinB}) and a second
one is an integral over $L+i{\mathbb R}$ (see, e.g., Section 2.2 of
\cite{Jo3})
for $L>0$; namely:
%
%
\begin{equation}
H_n(x)=\frac{n!}{2\pi i}\oint_\gamma e^{-z^2+2xz}\,\frac{dz}{z^{n+1}}
=\frac{2^ne^{x^2}}{i\sqrt{\pi}}\int_{L+i{\mathbb R}}e^{w^2-2xw}w^n\,dw.
\end{equation}
Then
%
%
\begin{eqnarray}\qquad\quad
&&K_{n,0}(x_1,t_1;x_2,t_2) \nonumber\\
&&\qquad=
-p(x_1,x_2;t_2-t_1)\Id(t_2>t_1)\nonumber\\[-8pt]\\[-8pt]
&&\qquad\quad{} + \frac{2}{ (2\pi i )^2}
\frac{e^{{x_1^2}/({2(1+t_1)})-{x_2^2}/({2(1+t_2)})}}
{\sqrt{(1+t_1)(1+t_2)}}
\nonumber\\
&&\qquad\quad\hspace*{10.5pt}{}\times
\int_{L+i{\mathbb R}}dU\oint_{\gamma}dV\frac{ (U/V)^n-1}{U-V}
\frac{e^{U^2({1-t_1})/({1+t_1}) -({2x_1U})/({1+t_1})}}{e^{V^2
({1-t_2})/({1+t_2})-({2x_2V})/({1+t_2})}}.\nonumber
\end{eqnarray}
Note that the $-1$ in $(\frac UV)^n-1$ (appearing in the integral
above) can actually be omitted, because there is no residue at $V=U$.
One then makes the substitution to new integration variables $\tilde U$
and $\tilde V$,
%
%
\begin{equation}
U\sqrt{\frac{1-t_1}{1+t_1}}=\tilde U\sqrt{\frac{n}{2}},\qquad
V\sqrt{\frac{1-t_2}{1+t_2}}=\tilde V\sqrt{\frac{n}{2}}
\end{equation}
and uses steepest descent in the integral to get the extended Airy
kernel (\ref{ExtAiry}),
which is just (\ref{eqKernelExtended}) in which one replaces $\tilde
a_k=\tilde b_k=0$ and $m=0$, namely
%
%
\begin{eqnarray}\hspace*{25pt}
&&K_\mathcal{A}(\tau_1,\xi_1;\tau_2,\xi_2)\nonumber\\
&&\qquad =-\frac{\Id(\tau
_2>\tau
_1)}{\sqrt{4\pi(\tau_2-\tau_1)}}
e^{-{(\xi_2-\xi_1)^2}/({4(\tau_2-\tau_1)})-({1}/{2})(\tau
_2-\tau
_1)(\xi_2+\xi_1)+({1}/{12}) (\tau_2-\tau_1)^3}\\
&&\qquad\quad{} + \frac{1}{ (2\pi i )^2}
\int_{\Gamma_{>}}d\omega\int_{\Gamma_{<}}d\widetilde\omega\frac
{e^{-\omega^3/3+\xi_2\omega}}{e^{-\widetilde\omega^{3}/3+\xi_1
\widetilde
\omega}}
\frac{1}{(\omega+\tau_2)-(\widetilde\omega+\tau_1)}.\nonumber
\end{eqnarray}

What remains is to compute the limit of the third term in (\ref
{2.Kext}), namely
%
%
\begin{equation}
\lim_{n\to\infty}\frac{1}{\sqrt{2}n^{1/6}} \sum_{i,j=1}^m \psi
^{(n)}_i(t_1,x_1)[\mu^{-1}]_{i,j} \varphi^{(n)}_j(t_2,x_2).
\end{equation}
Since $m$ remains finite, we can take the $n\to\infty$ limit inside the
sum. Also, the limit of $\mu^{-1}$, taking into account the prefactor,
has already been computed in Lem\-ma~\ref{lemma_mubar}. It remains to
determine the asymptotics of $\psi^{(n)}_k(t_1,x_1)$ and $\varphi
^{(n)}_k(t_2,x_2)$ (for $1\leq k \leq m$) under the above scaling.

As will be seen, the computations are very close to the ones for $t=0$
in Lemma \ref{Lemma1.4}.
For convenience, recall the notation $\gamma(t)=\sqrt{(1-t)/(1+t)}$
and $\sigma(t)=\sqrt{1+t}$. From (\ref{eqphipsim}), after the change of
variable $z=w/\gamma(t)$, one gets
%
%
\begin{equation}
\varphi^{(n)}_k(t,x) = \frac{e^{-x^2/2\sigma(t)^2}}{\sigma
(t)}\gamma(t)^n
\frac{1}{2\pi i} \oint_{\Gamma_{0,a_k'/2}} dw \frac
{e^{-w^2+2wx'}}{w^n(w-a_k'/2)},
\end{equation}
where
$ x' $ and $a_k'$ are defined below, together with their asymptotics:
%
\begin{eqnarray}\label{x,a}
x'&:=&\frac{x}{\sigma(t)^2\gamma(t)}=\frac{x}{\sqrt{1-t^2}}=\sqrt
{2n}+\frac{\xi}{\sqrt{2}
n^{1/6}}+\mathcal{O}(n^{-5/6}),\nonumber\\[-8pt]\\[-8pt]
a_k'&:=&a_k\gamma(t)=\sqrt{2n}+\sqrt{2}(\tilde a_k-\tau)
n^{1/6}+\mathcal{O}(n^{-1/6}).\nonumber
\end{eqnarray}
Now, we benefit from the computation made in the $\tau=0$ case. Indeed,
we showed that
%
%
\begin{eqnarray}
&&\frac{1}{2\pi i}\oint_{\Gamma_{0,a/2}} dz \frac
{e^{-z^2+2zy}}{z^n(z-a/2)}
\nonumber\\[-8pt]\\[-8pt]
&&\qquad= e^{y^2/2} \biggl(\frac{2e}{n}
\biggr)^{n/2}\frac{1}{2\pi i} \int_{\Gamma_{\tilde a>}}d\omega\frac
{e^{-\omega^3/3+\xi\omega}}{\omega-\tilde a} \bigl(1+o(1)\bigr),\nonumber
\end{eqnarray}
if $y$ and $a$ are scaled as
%
%
\begin{equation}
y=\sqrt{2n}+\frac{\xi}{\sqrt{2} n^{1/6}},\qquad
a=\sqrt{2n}+\sqrt{2}
\tilde a n^{1/6}.
\end{equation}
This is exactly our situation with $\tilde a=\tilde a_k-\tau$. Thus,
we get
%
%
\begin{eqnarray}
\varphi_k^{(n)}(t,x)&=&\frac{e^{-x^2/2\sigma(t)^2}}{\sigma(t)}\gamma(t)^n
e^{{x'}^2/2} \biggl(\frac{2e}{n} \biggr)^{n/2} \nonumber\\[-8pt]\\[-8pt]
&&{}\times \frac{1}{2\pi i}\int_{\Gamma_{\tilde a_k-\tau>}}d\omega\frac{e^{-\omega
^3/3+\xi\omega}}{\omega-\tilde a_k+\tau} \bigl(1+o(1)\bigr).\nonumber
\end{eqnarray}
Moreover, the asymptotics of the prefactor reads
%
%
\begin{equation}\label{3.35}
\frac{e^{-x^2/2\sigma(t)^2}}{\sigma(t)}\gamma(t)^n e^{{x'}^2/2} =
e^{-\tau^3/3+\xi\tau+\mathcal{O}(n^{-1/3})}.
\end{equation}
Thus, we have showed that
%
%
\begin{eqnarray}\label{eq3.37}
\varphi_k(\tau_2,\xi_2):\!&=&\lim_{n\to\infty}\varphi_k^{(n)}(t_2,x_2)
\biggl(\frac
{n}{2e} \biggr)^{n/2}\nonumber\\[-8pt]\\[-8pt]
&=&\frac{1}{2\pi i} \int_{\Gamma_{\tilde a_k-\tau
_2>}}d\omega\frac{e^{-\omega^3/3+\xi_2 \omega}}{\omega-\tilde
a_k+\tau
_2} \frac{1}{f(\tau_2,\xi_2)}\nonumber
\end{eqnarray}
and similarly,
%
%
\begin{eqnarray}\label{eq3.38}
\psi_k(\tau_1,\xi_1):\!&=&\lim_{n\to\infty}\psi_k^{(n)}(t_1,x_1)
\biggl(\frac
{n}{2e} \biggr)^{n/2}\nonumber\\[-8pt]\\[-8pt]
&=&\frac{1}{2\pi i} \int_{\Gamma_{<\tilde b_k-\tau
_1}}d\widetilde\omega\frac{e^{\widetilde\omega^3/3-\xi_1
\widetilde
\omega}}{\widetilde\omega-\tilde b_k+\tau_1} f(\tau_1,\xi_1).\nonumber
\end{eqnarray}

Now we can put together all the pieces, which make up the kernel (\ref
{2.Kext}), namely (\ref{eq3.26}), (\ref{eq3.37}), (\ref{eq3.38}) and
the asymptotics of the inverse of the matrix $B$ in Lem\-ma~\ref
{Lemma1.4}. Thus, we have
%
%
\begin{eqnarray}\label{eq3.39}\hspace*{35pt}
&&\lim_{n\to\infty} \frac{1}{\sqrt
{2}n^{1/6}}K_{n,m}(t_1,x_1;t_2,x_2)\frac
{f(\tau_2,\xi_2)}{f(\tau_1,\xi_1)} \nonumber\\[-8pt]\\[-8pt]
&&\qquad=K_\mathcal{A}(\tau_1,\xi
_1;\tau
_2,\xi_2)-\frac{f(\tau_2,\xi_2)}{f(\tau_1,\xi_1)}\sum_{i,j=1}^m \psi
_i(\tau_1,\xi
_1) [A^{-1}]_{i,j}\varphi_j(\tau_2,\xi_2).\nonumber
\end{eqnarray}
The last term in (\ref{eq3.39}) (including the minus sign) is equal to
%
%
\begin{eqnarray}
&&-\frac{1}{ (2\pi i )^2}\int_{\Gamma_{\tilde a-\tau_2>}} d\omega
\int_{\Gamma_{<\tilde b-\tau_1}} d\widetilde\omega\frac
{e^{-\omega
^3/3+\xi_2\omega}}{e^{-\widetilde\omega^{3}/3+\xi_1\widetilde
\omega}}\nonumber\\[-8pt]\\[-8pt]
&&\qquad{}\times
\sum_{i,j=1}^m \frac{[A^{-1}]_{i,j}}{(\widetilde\omega-\tilde
b_i+\tau
_1)(\omega-\tilde a_j+\tau_2)}.\nonumber
\end{eqnarray}
Applying the identity in Lemma \ref{LemmaIdentity}, we get as final
result the kernel $K_m^{\tilde a,\tilde b}(\tau_1,\xi_1;\break\tau_2,\xi_2)$
of Theorem \ref{th:2}, and this ends the proof of Proposition \ref
{Prop:kernel-ext}.
\end{pf}
\begin{pf*}{Proof of Theorem \protect\ref{th:2}}
For any bounded set $E$, the probability (\ref{0.Prob-ext}) is given by
the Fredholm determinant of the kernel, obtained in
Proposition \ref{Prop:kernel-ext}. Since this kernel is conjugate to
the one in Theorem \ref{th:2}, their Fredholm determinants are
identical.
\end{pf*}
\begin{pf*}{Proof of Theorem \protect\ref{th:2'}
\textup{(Universality)}}
The proof is a mild variation on the proof of Theorem \ref{th:2} and
Proposition \ref{Prop:kernel-ext}. Referring to the notation used in
the statement of Theorem \ref{th:2'}, one checks that formula (\ref
{ab}) for $a_kb_\ell/2$ remains the same, since $x_0^-x_0^+={2n}$ [see
(\ref{points})] and thus also asymptotic formula (\ref{eq2.35})
for~$\mu_{k\ell}$. Moreover, the scaling now reads
%
%
\begin{eqnarray}\label{x,t}\qquad
t &=& t_0+\frac{(1-t_0^2) \tau}{n^{1/3}},\nonumber\\
x &=& \sqrt{2n(1-t^2)} \biggl(1+\frac{\xi}{2n^{2/3}} \biggr)\\
&=& \sqrt{2n(1-t_0^2)}
\biggl(1-\frac{t_0\tau}{n^{1/3}}+\frac{\xi-\tau^2}{2n^{2/3}}-\frac
{t_0\tau}{2n}(\xi+\tau^2)+\mathcal{O}(n^{-4/3}) \biggr).\nonumber
\end{eqnarray}
Referring to the notation (\ref{x,a}), one checks that with $x$ and $t$
as in (\ref{x,t}) above, one has
%
%
\begin{eqnarray} \qquad
x'&=&\frac{x}{\sigma(t)^2\gamma(t)}=\frac{x}{\sqrt{1-t^2}}
=\sqrt{2n} \biggl(1+\frac{\xi}{2n^{2/3}} \biggr)+\mathcal{O}(n^{-5/6}),\nonumber
\nonumber\\
a_k'&=&a_k\gamma(t)
= a_k\sqrt{\frac{1-t}{1+t}}
=\sqrt{2n} \biggl(1-\frac{\tilde a_k-\tau
}{n^{1/3}} \biggr)+\mathcal{O}(n^{-1/6}),
\\
b_k'&=&b_k\gamma(-t)
=b_k\sqrt{\frac{1+t}{1-t}}
= \sqrt{2n} \biggl(1-\frac{\tilde b_k-\tau
}{n^{1/3}} \biggr)+\mathcal{O}(n^{-1/6}).\nonumber
\end{eqnarray}
With this information, one checks the following asymptotics, which is
the analogue of (\ref{3.35}), namely
%
%
\begin{eqnarray}
\hspace*{-10pt}&&\frac{e^{-x^2/2\sigma(\pm t)^2}}{\sigma(\pm t)}\gamma(\pm t)^n
e^{{x'}^2/2}\nonumber\\[-8pt]\\[-8pt]
&&\qquad=
\frac{e^{\pm{x'^2t}/{2}}}{\sqrt{1\pm t}} \biggl(\frac{1-t}{1+t}
\biggr)^{\pm n/2}
= \frac{e^{\pm nt_0}}{\sqrt{1\pm t_0}} \biggl(\frac{1- t_0}{1+ t_0}
\biggr)^{\pm n/2} f_n(\tau,
\xi)^{\mp1}e^{\mathcal{O}(n^{-1/3})}\nonumber\hspace*{-30pt}
\end{eqnarray}
with
%
%
\begin{equation}
f_n(\tau, \xi):=e^{t_0(n^{2/3}\tau t_0-(\xi-\tau^2)n^{1/3}+ t_0\tau
(\xi
+\tau^2))}e^{\tau^3/3-\xi\tau}.
\end{equation}
Here, $f_n(\tau, \xi)$ depends on $n$, besides $\tau$ and $\xi$.
Then we show
%
%
\begin{eqnarray}\label{eq3.37'}
&&\lim_{n\to\infty}\varphi_k^{(n)}(t_2,x_2) \biggl(\frac{n}{2e}
\biggr)^{n/2}e^{-nt_0} \biggl(\frac{1+t_0}{1-t_0} \biggr)^{n/2}\sqrt
{1+t_0} f_n(\tau_2,\xi_2) \nonumber\\[-8pt]\\[-8pt]
&&\qquad=\frac{1}{2\pi i} \int_{\Gamma_{\tilde a_k-\tau_2>}}d\omega\frac
{e^{-\omega^3/3+\xi_2 \omega}}{\omega-\tilde a_k+\tau_2},\nonumber
\end{eqnarray}
and similarly,
%
%
\begin{eqnarray}\label{eq3.38'}
&&\lim_{n\to\infty}\psi_k^{(n)}(t_1,x_1) \biggl(\frac{n}{2e}
\biggr)^{n/2}e^{nt_0}\nonumber
\biggl(\frac{1-t_0}{1+t_0} \biggr)^{n/2}\sqrt{1-t_0} \frac
{1}{f_n(\tau_1,\xi_1)} \\[-8pt]\\[-8pt]
&&\qquad=\frac{1}{2\pi i} \int_{<\Gamma_{\tilde b_k-\tau_1}}d\tilde\omega
\frac
{e^{\tilde\omega^3/3-\xi_1 \tilde\omega}}{\tilde\omega-\tilde
b_k+\tau_1}.\nonumber
\end{eqnarray}
Also, as before,
%
%
\begin{equation}
\lim_{n\rightarrow\infty}\frac{\mu^{-1}}{\sqrt2 n^{1/6}} \biggl(\frac
{2e}{n} \biggr)^n=-A^{-1}.
\end{equation}
Then, with $x_i,t_i$ as in (\ref{x,t}), the limit (\ref{eq3.26}) gets
replaced by
\begin{eqnarray*}
&&\lim_{n\to\infty} \frac{\sqrt{1-t_0^2}}{\sqrt
{2}n^{1/6}}K_{n,0}(t_1,x_1;t_2,x_2)\frac{f_n(\tau_2,\xi_2)}{f_n(\tau
_1,\xi_1)}
= K_\mathcal{A}(\tau_1,\xi_1;\tau_2,\xi_2)
\end{eqnarray*}
with very little change in the steepest descent argument.
So, putting all the pieces together, one checks
%
%
\begin{eqnarray}
\hspace*{25pt}&&\lim_{n\to\infty} \frac{\sqrt{1-t_0^2}}{\sqrt
{2}n^{1/6}}K_{n,m}(t_1,x_1;t_2,x_2)\frac{f_n(\tau_2,\xi_2)}{f_n(\tau
_1,\xi_1)}
\nonumber\\
&&\qquad
= K_\mathcal{A}(\tau_1,\xi_1;\tau_2,\xi_2)
-\frac{1}{ (2\pi i )^2}\int_{\Gamma_{\tilde a-\tau_2>}} d\omega
\int_{\Gamma_{<\tilde b-\tau_1}} d\widetilde\omega\frac
{e^{-\omega
^3/3+\xi_2\omega}}{e^{-\widetilde\omega^{3}/3+\xi_1\widetilde
\omega}}\\
&&\qquad\quad\hspace*{88.3pt}{}\times
\sum_{i,j=1}^m \frac{[A^{-1}]_{i,j}}{(\widetilde\omega-\tilde
b_i+\tau
_1)(\omega-\tilde a_j+\tau_2)}\nonumber
\end{eqnarray}
from which one proceeds in the same way as in the proof of
Proposition \ref{Prop:kernel-ext} and Theorem \ref{th:2}. This ends
the proof of Theorem \ref{th:2'}.
\end{pf*}

\section{Airy process with wanderers all leaving from point $a$ and all
going to point $b$}\label{SectTwoPackets}
In this section, we prove Theorem \ref{th:2} (and thus also Theorem
\ref{th:1}) for the case where $m$ wanderers all leave from one point and
all are forced to one point; that is,
%
%
\begin{equation}
\tilde a:=\tilde a_m=\cdots=\tilde a_1<\tilde b_1=\cdots=\tilde
b_m=:\tilde b.
\end{equation}
Thus, the $m$ top Brownian bridges start from $a$ and end at $b$ with
%
%
\begin{equation}\label{eq4.1}
a=\sqrt{2n}(1+\tilde a n^{-1/3}),\qquad b=\sqrt{2n}(1-\tilde b n^{-1/3}).
\end{equation}

The arguments presented in the previous sections break down. Therefore,
one should redo the proof, using an argument adapted to this case. It
is instructive to shortly present two different approaches. The first
follows the approach of the previous section, consisting in computing
the inverse of the $m\times m$ matrix $\mu$, and the second approach is
to perform the biorthogonalization. In principle, with some care
because of the $n\to\infty$ limit, one might also be able to do the
argument by analytic continuation, since the measure is analytic in the
$\tilde a_i$, $\tilde b_j$ as well as the final kernel (provided the
inequality $\tilde a_i<\tilde b_j$ for all $i,j$ is satisfied).

\subsection{Via the inversion of the moment matrix}
The start is almost the same as in the previous section. The only
difference is that the first and last
determinant in the measure,
instead of (\ref{eq3.6}) and (\ref{eq3.7}), are now
%
%
\begin{equation}\label{eq4.2}
\det\pmatrix{ \bigl((x_j^1)^{i-1} e^{a x_j^1/(1+T_1)}p(0,x_j^1,T_1+1) \bigr)_{
{\fontsize{8.36}{10}\selectfont{\begin{array}{l} 1\leq i\leq m\\ 1\leq j\leq m+n
\end{array}}}}
 \cr\bigl((x_j^1)^{i-1}p(0,x_j^1,T_1+1) \bigr)_{
{\fontsize{8.36}{10}\selectfont{\begin{array}{l} 1\leq i\leq n\\ 1\leq j\leq m+n
\end{array}}}}}
\end{equation}
and
%
%
\begin{equation}\label{eq4.3}
\det\pmatrix{
\bigl((x_j^\ell)^{i-1} e^{b x_j^\ell/(1-T_\ell)}p(x_j^\ell,0,1-T_\ell
) \bigr)_{{\fontsize{8.36}{10}\selectfont{
\begin{array}{l} 1\leq i\leq m\\ 1\leq j\leq m+n
\end{array}}}}
\cr\bigl((x_j^\ell)^{i-1}p(x_j^\ell,0,1-T_\ell) \bigr)_{
{\fontsize{8.36}{10}\selectfont{\begin{array}{l} 1\leq i\leq n\\ 1\leq j\leq m+n
\end{array}
}}}},
\end{equation}
respectively. The functions $\varphi^{(n)}_k$ and $\psi^{(n)}_k$, for
$1\leq k \leq m$, defined by
%
%
\begin{eqnarray}\label{eq4.4}
\varphi_k^{(n)}(t,x) &=& \frac{e^{-x^2/2\sigma(t)^2}}{\sigma(t)}
\frac
{1}{2\pi i}\oint_{\Gamma_{0,a/2}} dz\frac{e^{-z^2\gamma
(t)^2+2xz/\sigma
^2(t)}}{z^n (z-a/2 )^k},\nonumber\\[-8pt]\\[-8pt]
\psi_k^{(n)}(t,x) &=& \frac{e^{-x^2/2\sigma(-t)^2}}{\sigma(-t)}
\frac
{1}{2\pi i}\oint_{\Gamma_{0,b/2}} dz\frac{e^{-z^2\gamma
(-t)^2+2xz/\sigma
^2(-t)}}{z^n (z-b/2 )^k},\nonumber
\end{eqnarray}
replace those of (\ref{eqphipsim}), where we recall that
$\gamma(t)=\sqrt{\frac{1-t}{1+t}}, \sigma(t)=\sqrt{1+t}$.
Of course, since the last $n$ rows of the determinants (\ref{eq4.2})
and (\ref{eq4.3}) are exactly the same as in (\ref{eq3.6}) and (\ref
{eq3.7}), we keep the same choice for the functions $\varphi^{(n)}_{m+k}$
and $\psi^{(n)}_{m+k}$, $1\leq k\leq n$, as in (\ref
{eqphipsin}) and (\ref
{eqphipsinB}). Define the $m\times m$ matrix $\mu$ by $\mu
_{i,j}=\langle\varphi^{(n)}_{i},\psi^{(n)}_{j}\rangle$, $1\leq
i,j\leq m$.
Once again, this choice of $\varphi^{(n)}_{k}$ and $\psi^{(n)}_{k}$
generates the same vector space as the function in the above
determinants (\ref{eq4.2}) and~(\ref{eq4.3}).

Note that in this section we use the same notation as in the previous
section. However, the matrix $\mu$ and some of the functions are not
the same. What remains the same is the form of the kernel. Indeed,
since Proposition \ref{prop:Kext} holds exactly as before, one has
once again
%
%
\begin{eqnarray}\label{eq4.7}
K_{n,m}(t_1,x_1;t_2,x_2) &=&
K_{n,0}(t_1,x_1;t_2,x_2)\nonumber\\[-8pt]\\[-8pt]
&&{}+\sum_{i,j=1}^{m}
\psi^{(n)}_i(t_1,x_1)[\mu^{-1}]_{i,j}
\varphi^{(n)}_j(t_2,x_2),\nonumber
\end{eqnarray}
where $K_{n,0}$ is the kernel without wanderers and in the scaling
limit will converge to the extended Airy kernel. Thus, we only have to
deal with the double sum below.
\begin{lemma}\label{lemma4mu}
Under the scaling (\ref{eq4.1}), we have
%
%
\begin{equation}\label{eq4.7L}
\lim_{n\to\infty}\mu_{k,l} \biggl(\frac{n}{2e} \biggr)^{n}
\bigl(n^{1/6}/\sqrt{2} \bigr)^{k+l-1} = \frac{1}{2}\frac{1}{(\tilde b-\tilde
a)^{l+k-1}}\pmatrix{l+k-2\cr k-1}
\end{equation}
for $1\leq k,l\leq m$.
\end{lemma}
\begin{pf}
For convenience, in the proof we compute $\mu_{k+1,l+1}$ to avoid
$-1$'s in the formulas. Since, as before, $\mu_{k\ell}$ is
time-independent, we may set $t_1=t_2=0$ in the computation; so, as in
(\ref{eq2.21}) and after integrating over the $x$ variable, one finds
%
%
\begin{equation}
\mu_{k+1,l+1}=\frac{\sqrt{\pi}}{(2\pi i)^2}\oint_{\Gamma_{0,a/2}}dz
\frac{1}{z^n(z-a/2)^{k+1}}\oint_{\Gamma_{0,b/2}}dw \frac{e^{2wz}}{w^n
(w-b/2)^{l+1}}.\hspace*{-33pt}
\end{equation}
Then we apply twice the identity
%
%
\begin{equation}
\frac{1}{(z-a/2)^{k+1}}=\frac{2^k}{k!} \biggl(\frac{\partial}{\partial a
} \biggr)^k\frac{1}{z-a/2}
\end{equation}
and obtain
%
%
\begin{equation}
\mu_{k+1,l+1}=\frac{2^{k+l}}{k!l!} \biggl(\frac{\partial}{\partial a
} \biggr)^k \biggl(\frac{\partial}{\partial a } \biggr)^l \mu_{1,1}.
\end{equation}
But $\mu_{1,1}$ was already expressed as a single contour integral; see
(\ref{1.mukl}). Thus,
%
%
\begin{equation}
\mu_{k+1,l+1}=\frac{2^{k+l}}{k!l!} \biggl(\frac{\partial}{\partial a
} \biggr)^k \biggl(\frac{\partial}{\partial b} \biggr)^l \frac{\sqrt{\pi}
2^n}{2\pi i}\oint_{\Gamma_{0,ab/2}} dz \frac{e^z}{z^n(z-ab/2)}.
\end{equation}
We now compute the derivatives of $(z-ab/2)^{-1}$ and obtain
%
%
\begin{eqnarray}
&&\frac{2^{k+l}}{k!l!} \biggl(\frac{\partial}{\partial a } \biggr)^k
\biggl(\frac{\partial}{\partial a } \biggr)^l\frac{1}{z-{ab}/{2}}
\nonumber\\[-8pt]\\[-8pt]
&&\qquad=a^l b^k \sum_{j=0}^{\min(k,l)} \frac{(ab/2)^{-j}}{(z-
{ab}/{2})^{l+k-j+1}} \frac{(l+k-j)!}{(k-j)!(l-j)!j!},\nonumber
\end{eqnarray}
and so
%
%
\begin{eqnarray}\label{eq4.14}\quad
\mu_{k+1,l+1}&=&a^l b^k \sum_{j=0}^{\min(k,l)} \frac
{1}{(ab/2)^{j}}\frac{(l+k-j)!}{(k-j)!(l-j)!j!}\nonumber\\[-8pt]\\[-8pt]
&&\hspace*{49.3pt}{}\times\frac{\sqrt{\pi} 2^n}{2\pi i}\oint_{\Gamma_{0,ab/2}}dz
\frac
{e^z}{z^n (z-ab/2)^{l+k-j+1}}.\nonumber
\end{eqnarray}

Finally, we need to do the asymptotic analysis of the integral, which
essentially has already been made in
Lemma \ref{lemma_mubar}. Consider the following small change in~(\ref
{eq2.35}), with $a_k=a$, $b_\ell=b$: replace $(z-ab/2)^{-1}$ by
$(z-ab/2)^{-p-1}$ for any \textit{finite} $p=1,2,\ldots.$ Then the
steepest descent analysis is unchanged except for that finite power,
which would be present in (\ref{eq2.39}) too. This extra power gives a
factor $n^{-p/6}$, and we also have an extra factor coming from the
change of variables equal to $n^{-p/2}$. In the end, the result is
%
%
\begin{eqnarray}\label{eq4.15}
&&\frac{\sqrt{\pi} 2^n}{2\pi i}\oint_{\Gamma_{0,ab/2}}dz \frac
{e^z}{z^n (z-ab/2)^{p+1}} \nonumber\\[-8pt]\\[-8pt]
&&\qquad=
\biggl(\frac{2e}{n} \biggr)^n \frac{1}{(\tilde b-\tilde a)^{p+1}}\frac
{n^{-2p/3}}{\sqrt{2} n^{1/6}}\bigl(1+o(1)\bigr).\nonumber
\end{eqnarray}
We put (\ref{eq4.15}) into (\ref{eq4.14}) and compare the dependence in
$j$ of the different terms in the sum. Since $ab/2=n(1+\mathcal{O}(n^{-1/3}))$,
the $j$th term in the sum (\ref{eq4.14}) contains the following power
of $n$, namely
%
%
\begin{equation}
\frac{1}{n^{j}} n^{2j/3}=n^{-j/3}.
\end{equation}
Therefore, since the sum is finite, in the $n\to\infty$ limit, the
leading term is the one with $j=0$, the other ones being of smaller
order. Thus,
%
%
\begin{eqnarray}\label{4.17}\qquad
\mu_{k+1,l+1}&=&\pmatrix{k+l\cr k} \biggl(\frac{2e}{n} \biggr)^n
\frac{1}{(\tilde b-\tilde a)^{k+l+1}}\frac{a^l b^k
n^{-2(k+l)/3}}{\sqrt
{2} n^{1/6}}\bigl(1+o(1)\bigr)\nonumber\\[-8pt]\\[-8pt]
&=&\frac12 \pmatrix{k+l\cr k} \biggl(\frac{2e}{n} \biggr)^n
\biggl(\frac{\sqrt{2}}{n^{1/6}(\tilde b-\tilde a)}
\biggr)^{k+l+1}\bigl(1+o(1)\bigr),\nonumber
\end{eqnarray}
ending the proof of Lemma \ref{lemma4mu}.
\end{pf}
\begin{corollary}\label{CorMuInverse}
With the same scaling as in Lemma \ref{lemma4mu}, we have
%
%
\begin{equation}\label{4.asympt}
\qquad\lim_{n\to\infty}[\mu^{-1}]_{k,l} \biggl(\frac{2e}{n} \biggr)^{n}
\biggl(\frac{\sqrt{2}}{n^{1/6}} \biggr)^{k+l-1} = 2(\tilde b-\tilde
a)^{l+k-1}[(L^{-1})^{T} L^{-1}]_{k,l}
\end{equation}
for $1\leq k,l\leq m$, where
the $m\times m$ lower-triangular matrix $L^{-1}$ has binomial entries
%
\begin{equation}
(L^{-1})_{k,l}=(-1)^{k-l} \pmatrix{k-1\cr l-1}.
\end{equation}
\end{corollary}
\begin{pf}
From (\ref{eq4.7L}), it follows that computing the inverse of
$\mu$ reduces to computing the inverse of the $m\times m$ matrix $\nu$,
%
%
\begin{equation}
\nu= (\nu_{k,l})_{1\leq k,l\leq m}\qquad\mbox{with }\nu_{k,l}=\pmatrix
{k+l-2\cr k-1}.
\end{equation}
A convenient way of taking the inverse is to compute the $\nu=LL^T$
decomposition, where $L$ is lower-triangular,
yielding
%
%
\begin{equation}
\nu= L L^{T},\qquad L_{k,l}=\pmatrix{k-1\cr l-1},
\end{equation}
from which
%
%
\begin{equation}
\nu^{-1}=(L^{T})^{-1} L^{-1},\qquad (L^{-1})_{k,l}=(-1)^{k-l} \pmatrix
{k-1\cr l-1}.
\end{equation}
This establishes the asymptotics (\ref{4.asympt}).
\end{pf}

We now turn to the asymptotics of the functions $\varphi^{(n)}_k$ and
$\psi
^{(n)}_k$, defined in (\ref{eq4.4}).
\begin{lemma}\label{lemma4phipsi}
Under the scaling (\ref{eq4.1}) and\vspace*{-1pt}
%
%
\begin{equation}\label{eq4Scaling}
t_i=\tau_i n^{-1/3},\qquad x_i=\sqrt{2n}+\frac{\xi_i-\tau_i^2}{\sqrt
{2} n^{1/6}},
\end{equation}\vspace*{-1pt}
one has\vspace*{-1pt}
%
%
\begin{eqnarray}
\varphi_k(\tau_2,\xi_2):\!&=&\lim_{n\to\infty} \varphi^{(n)}_k(t_2,x_2)
\biggl(\frac{n}{2e} \biggr)^{n/2}
\bigl(n^{1/6}/\sqrt{2}\bigr)^{k-1}\nonumber\\[-9pt]\\[-9pt]
& =& \frac{1}{2\pi i} \int_{\Gamma_{\tilde a-\tau_2>}}d\omega\,
e^{-\omega^3/3+\xi_2 \omega}\frac{f(\tau_2,\xi_2)^{-1}}{(\omega
-\tilde
a+\tau_2)^k}\nonumber
\end{eqnarray}\vspace*{-1pt}
and\vspace*{-1pt}
%
%
\begin{eqnarray}
\psi_\ell(\tau_1,\xi_1):\!&=&\lim_{n\to\infty} \psi^{(n)}_\ell(t_1,x_1)
\biggl(\frac{n}{2e} \biggr)^{n/2} \bigl(n^{1/6}/\sqrt{2}\bigr)^{\ell-1}
\nonumber\\[-9pt]\\[-9pt]
&=& \frac{(-1)^{\ell-1}}{2\pi i} \int_{\Gamma_{<\tilde b-\tau
_1}}d\widetilde\omega\, e^{\widetilde\omega^3/3-\xi_1 \widetilde
\omega
}\frac{{f(\tau_1,\xi_1)}}{(\widetilde\omega-\tilde
b+\tau_1)^\ell}.\nonumber
\end{eqnarray}\vspace*{-1pt}
\end{lemma}
\begin{pf}
We must compute the asymptotics of\vspace*{-1pt}
%
%
\begin{equation}
\varphi^{(n)}_k(t,x)=\frac{e^{-x^2/2\sigma(t)^2}}{\sigma(t)} \frac
{1}{2\pi
i}\oint_{\Gamma_{0,a/2}} dz\frac{e^{-z^2\gamma(t)^2+2xz/\sigma
^2(t)}}{z^n (z-a/2 )^k}
\end{equation}\vspace*{-1pt}
and compare this expression with (\ref{eqphipsim}). One sees that the
only differences are that now $a$ replaces $a_k$ and the denominator in
$z-a/2$ has a power $k$ instead of power $1$. For any \textit{finite}
$k$, the asymptotic analysis for this case has only minor differences
with respect to the asymptotic for (\ref{eqphipsim}). Namely, one picks
up some extra factors by the changes of variables: setting $z=u\sqrt
{n/2}$, one gets a factor $(2/n)^{k-1}$ and then $\omega
=n^{1/3}(u-1)$ results in a factor $n^{(k-1)/3}$. In total, an extra
factor $(\sqrt{2}/n^{1/6})^{k-1}$ appears. A similar argument holds for
$\psi^{(n)}_\ell$.\vspace*{-1pt}
\end{pf}
\begin{pf*}{Proof of Theorems \protect\ref{th:1} and \protect\ref
{th:2} for Brownian bridges
starting from $a$ and ending up at $b$}
Putting together Corollary \ref{CorMuInverse} and Lemma \ref
{lemma4phipsi}, one obtains\vspace*{-1pt}
%
%
\begin{eqnarray}\label{eq4.26}
\hspace*{40pt}&&\lim_{n\to\infty} \frac{1}{\sqrt{2} n^{1/6}}\sum
_{i,j=1}^m\psi^{(n)}_i(t_1,x_1) [\mu^{-1}]_{i,j} \varphi^{(n)}_j(t_2,x_2)
\nonumber\\[-1pt]
&&\qquad = \sum_{i,j=1}^m\psi_i(\tau_1,\xi_1) (\tilde b-\tilde a)^{i+j-1}
[(L^{T})^{-1} L^{-1}]_{i,j} \varphi_j(\tau_2,\xi_2)
\frac{f(\tau_1,\xi_1)}{f(\tau_2,\xi_2)}
\nonumber\\[-9pt]\\[-9pt]
&&\qquad
= \frac{1}{(2\pi i)^2}\int_{\Gamma_{\tilde a-\tau_2>}}d\omega\int
_{\Gamma_{<\tilde b-\tau_1}}d\widetilde\omega\frac{e^{-\omega
^3/3+\xi
_2 \omega}}{e^{-\widetilde\omega^3/3+\xi_1\widetilde\omega}}
\frac
{f(\tau_1,\xi_1)}{f(\tau_2,\xi_2)}(\tilde a-\tilde b)
\nonumber\\[-1pt]
&&\qquad\quad{} \times\sum_{k=1}^m\sum_{i,j=1}^m
\pmatrix{k-1\cr i-1}\pmatrix{k-1\cr j-1}
\frac{(\tilde b-\tilde a)^{i+j-2}
(-1)^j}{(\widetilde\omega+\tau_1-\tilde b)^i (\omega+\tau_2-\tilde a)^j}.
\nonumber
\end{eqnarray}
Finally, using the fraction decomposition identity
%
%
\begin{eqnarray}
& &\frac{1}{V-U} \biggl( \biggl(\frac{U-a}{U-b} \biggr)^m \biggl(\frac
{V-b}{V-a} \biggr)^m-1 \biggr)\nonumber\\[-8pt]\\[-8pt]
& &\qquad= (a-b)\sum_{k=1}^m\sum_{i,j=1}^{m} \pmatrix{k-1\cr i-1}\pmatrix
{k-1\cr j-1}
\frac{(b-a)^{i+j-2}(-1)^{j}}{(U-b)^{i}(V-a)^{j}} \nonumber
\end{eqnarray}
with $U=\widetilde\omega+\tau_1$ and $V=\omega+\tau_2$, one gets the
final result
%
%
\begin{eqnarray}\label{eq4.29}
\hspace*{30pt}(\ref{eq4.26})&=& \frac{1}{(2\pi i)^2}\int_{\Gamma_{\tilde a-\tau
_2>}}d\omega\int_{\Gamma_{<\tilde b-\tau_1}}d\widetilde\omega
\frac
{e^{-\omega^3/3+\xi_2 \omega}}{e^{-\widetilde\omega^3/3+\xi
_1\widetilde
\omega}}
\frac{1}{(\omega+\tau_2)-(\widetilde\omega+\tau_1)}\nonumber\\[-8pt]\\[-8pt]
&&{}\times
\biggl[ \biggl(\frac{(\widetilde\omega-\tilde a+\tau_1)(\omega-\tilde
b+\tau_2)}{(\widetilde\omega-\tilde b+\tau_1)(\omega-\tilde a+\tau_2)}
\biggr)^m-1 \biggr]
 \frac{f(\tau_1,\xi_1)}{f(\tau_2,\xi_2)}.\nonumber
\end{eqnarray}
\upqed\end{pf*}

\subsection{Via biorthogonal functions}
Here, we present a slightly different approach, which consists in
using biorthogonal functions, instead of the functions defined in (\ref
{eq4.4}). We use a representation with determinants known from
classical orthogonal polynomial theory. Let us first define the
polynomials $\tilde\varphi^{(n)}_i$, $\tilde\psi^{(n)}_j$, and then show
they are actually biorthogonal.

Set
%
%
\begin{equation}
\mu_{i,j}(t):=\bigl\langle\varphi^{(n)}_i(0,\cdot),\psi^{(n)}_j(0,\cdot
)\bigr\rangle
\quad\mbox{and}\quad \Delta_k:=\det(\mu_{i,j})_{1\leq i,j\leq k}
\end{equation}
and define
%
%
\begin{equation}\label{eq4tildePhi}\hspace*{33pt}
\tilde\varphi^{(n)}_k(t,x):=\frac{1}{\sqrt{\Delta_k\Delta
_{k-1}}}\det
\pmatrix{
\mu_{1,1} & \mu_{1,2} & \cdots& \mu_{1,k-1} & \varphi^{(n)}_1(t,x)
\cr
\mu_{2,1} & \mu_{2,2} & \cdots& \mu_{2,k-1} & \varphi^{(n)}_2(t,x)
\cr
\vdots& \vdots& & \vdots& \vdots\cr
\mu_{k,1} & \mu_{k,2} & \cdots& \mu_{k,k-1} & \varphi^{(n)}_k(t,x)}
\end{equation}
and
%
%
\begin{equation}\label{eq4tildePsi}
\tilde\psi^{(n)}_l(t,x):=\frac{1}{\sqrt{\Delta_l\Delta
_{l-1}}}\det
\pmatrix{\mu_{1,1} & \mu_{1,2} & \cdots& \mu_{1,l} \cr
\mu_{2,1} & \mu_{2,2} & \cdots& \mu_{2,l} \cr
\vdots& \vdots& & \vdots\cr
\mu_{l-1,1} & \mu_{l-1,2} & \cdots& \mu_{l-1,l} \cr
\psi^{(n)}_1(t,x) & \psi^{(n)}_2(t,x) & \cdots&
\psi^{(n)}_l(t,x)}.\hspace*{-37pt}
\end{equation}
First of all, notice that $\tilde\varphi^{(n)}_k$ is linear
combination of
the $\varphi^{(n)}_\ell$ with $\ell=1,\ldots,k$, with a nonzero coefficient
in front of $\varphi^{(n)}_k$ (because $\Delta_k\neq0$, since both
$\{
\varphi
^{(n)}_\ell, 1\leq\ell\leq k\}$ and $\{\psi^{(n)}_\ell, 1\leq\ell
\leq
k\}$ form a basis of a $k$-dimensional vector space). The argument is
similar for $\tilde\psi^{(n)}_l$.
Observe, for $\ell<k$,
%
\begin{equation}
\bigl\langle\tilde\varphi^{(n)}_k(t,\cdot), \psi^{(n)}_\ell(t,\cdot
)\bigr\rangle= \left\{\begin{array}{l}\mbox{right-hand side of (\ref{eq4tildePhi}) with}
\\
\mbox{the last column replaced by}
\\
\mbox{the $\ell$th column of (\ref{eq4tildePhi})}
\end{array}
\right\}=0.
\end{equation}
Therefore, also $\langle\tilde\varphi^{(n)}_k(t,\cdot), \tilde\psi
^{(n)}_\ell(t,\cdot)\rangle=0$ for $\ell<k$ and thus also for $\ell
\neq k$, by merely interchanging the roles of $\tilde\varphi$ and
$\tilde\psi$. The above argument also shows
%
\begin{eqnarray}
\bigl\langle\tilde\varphi^{(n)}_k(t,\cdot), \tilde\psi^{(n)}_k(t,\cdot
)\bigr\rangle &=& \bigl\langle\tilde\varphi^{(n)}_k(t,\cdot), \psi
^{(n)}_k(t,\cdot
)\bigr\rangle
\frac{\Delta_{k-1}}{\sqrt{\Delta_k\Delta_{k-1}}^{}}
\nonumber\\[-8pt]\\[-8pt]
&=&\frac{\Delta_{k}\Delta_{k-1}}{\sqrt{\Delta_k\Delta_{k-1}}^{ 2}}
=1.\nonumber
\end{eqnarray}
The consequence is that now the kernel instead of (\ref{eq4.7}) reads
%
%
\begin{eqnarray}\label{eq4.36}
K_{n,m}(t_1,x_1;t_2,x_2) &=& K_{n,0}(t_1,x_1;t_2,x_2)
\nonumber\\[-8pt]\\[-8pt]
&&{} + \sum_{i,j=1}^m
\tilde\psi^{(n)}_i(t_1,x_1) [\tilde\mu^{-1}]_{i,j}\tilde\varphi
^{(n)}_j(t_2,x_2)\nonumber
\end{eqnarray}
with $\tilde\mu_{i,j}=\langle\tilde\varphi^{(n)}_i(t,\cdot),
\tilde
\psi
^{(n)}_j(t,\cdot)\rangle=\delta_{i,j}$.
Therefore, the double sum in (\ref{eq4.36}) becomes just
%
\begin{eqnarray}
K_m(t_1,x_1;t_2,x_2)&=&\sum_{i,j=1}^m \tilde\psi^{(n)}_i(t_1,x_1)
[\tilde\mu^{-1}]_{i,j}\tilde\varphi^{(n)}_j(t_2,x_2)
\nonumber\\[-8pt]\\[-8pt]
&=& \sum_{i=1}^m
\tilde\psi^{(n)}_i(t_1,x_1)\tilde\varphi{}^{(n)}_i(t_2,x_2).\nonumber
\end{eqnarray}
This last sum is of Darboux-type and can be rewritten as
%
%
\begin{eqnarray}\label{Darboux}
\hspace*{30pt}
&&K_m(t_1,x_1;t_2,x_2)\nonumber
\\[-8pt]\\[-8pt]
&&\qquad= -\frac{1}{\Delta_m}\det
{\fontsize{9.3}{9.3}\selectfont
\left(\matrix{0          & \psi^{(n)}_1(t_1,x_1)          & \psi^{(n)}_2(t_1,x_1) & \cdots & \psi^{(n)}_m(t_1,x_1) \cr
\varphi^{(n)}_1(t_2,x_2) & \mu_{1,1}                      & \mu_{1,2}             & \cdots & \mu_{1,m}\cr
\varphi^{(n)}_1(t_2,x_2) & \mu_{2,1}                      & \mu_{2,2}             & \cdots & \mu_{2,m}\cr
\vdots                   & \vdots                         & \vdots                &        & \vdots\cr
\varphi^{(n)}_m(t_2,x_2) & \mu_{m,1}                        & \mu_{m,2}             & \cdots & \mu_{m,m}}\right)},\nonumber
\end{eqnarray}
the latter follows from the fact that $K_m(t_1,x_1;t_2,x_2)$ is a
bilinear combination of $\psi^{(n)}_i(t_1,x_1) $ and $\varphi
^{(n)}_j(t_2,x_2)$, for $
1\leq i,j\leq m $ and is completely characterized
by $\langle K_m(t_1,x_1;t_2,\cdot),
\psi^{(n)}_i(t_2,\cdot)
\rangle= \psi^{(n)}_i(t_1,x_1)$ for $1\leq i\leq m $.

At this point, we have to determine the $n\to\infty$ limit of the
rescaled kernel, namely
%
%
\begin{equation}\label{eq4.39}
\lim_{n\to\infty}\frac{1}{\sqrt{2} n^{1/6}}K_m(t_1,x_1;t_2,x_2)
\end{equation}
with $x_i,t_i$ scaled as (\ref{eq4Scaling}). The asymptotics of $\mu
_{i,j}$ already appears in Lemma \ref{lemma4mu}, and for\vspace*{2pt} $\varphi^{(n)}_j$
and $\psi^{(n)}_i$ in Lemma \ref{lemma4phipsi}. Hence, Lemma \ref
{lemma4mu}, together with the fact that $m$ is finite, yields the
asymptotics of $\Delta_m$:
%
%
\begin{eqnarray}\label{eq4.40}
\hspace*{20pt}\lim_{n\to\infty} \Delta_m \biggl(\frac{n}{2e} \biggr)^{n m} \biggl(\frac
{n^{1/6}}{\sqrt{2}} \biggr)^{m^2}
&=& \frac{1}{2^m}\frac{1}{(\tilde b-\tilde a)^{m^2}}
\operatorname{det}\left(\pmatrix{k+l-2\cr k-1}\right)_{1\leq k,l\leq m}
\nonumber\\[-8pt]\\[-8pt]\nonumber\\
&=&\frac{1}{2^m}\frac{1}{(\tilde b-\tilde a)^{m^2}}.\nonumber\hspace*{-100pt}
\end{eqnarray}
This result (\ref{eq4.40}) and the linearity of the determinant,
together with the results of Lemmas \ref{lemma4mu} and \ref
{lemma4phipsi}, substituted in (\ref{Darboux}), lead to the limit in
(\ref{eq4.39}), namely\vspace*{-6pt}

\renewcommand{\theequation}{\fontsize{11}{11}\selectfont\arabic{section}.\arabic{equation}}
%
{{\fontsize{9.5}{9.5}\selectfont
\begin{eqnarray}
\hspace*{-26pt}&&\frac{1}{(2\pi i)^2} \int_{\Gamma_{\tilde a-\tau_2>}}d\omega
\int_{\Gamma_{<\tilde b-\tau_1}}d\widetilde\omega\frac{e^{-\omega^3/3+\xi_2\omega}}{e^{-\widetilde\omega^3/3+\xi_1\widetilde\omega}}
\frac{-(\tilde b-\tilde a)}{(\omega+\tau_2-\tilde a)(\widetilde\omega
+\tau_1-\tilde b)} \hspace*{-36pt}\nonumber\\
&&\qquad{}\times
\det{\fontsize{10.5}{10.5}\selectfont{\pmatrix{0 & 1 & \cdots& \biggl(\dfrac{\tilde b-\tilde a}{\omega+\tau_2-\tilde a} \biggr)^{m-1} \cr
1 & & & \cr
\vdots& &\left(\pmatrix{k+l-2\cr k-1}\right)_{1\leq k,l\leq m} & \cr
\biggl(\dfrac{-(\tilde b-\tilde a)}{\widetilde\omega+\tau_1-\tilde b} \biggr)^{m-1} & & &}}}\hspace*{-36pt}\\
&&\qquad{}\times
\frac{f(\tau_1,\xi_1)}{f(\tau_2,\xi_2)}.\nonumber\hspace*{-36pt}
\end{eqnarray}
}
The final\vspace*{1pt} step is to use the following identity, established by
observing that both sides are identical upon integrating against
$x^{i-1}$, $1\leq i\leq m$, from $x=0$ to $x=1$,
%
%
\begin{eqnarray}\label{4.40}
&&-\det
\pmatrix{0 & 1 & \cdots&(x+1)^{m-1} \cr
1 & & & \cr
\vdots& &\left(\pmatrix{k+l-2\cr k-1}\right)_{1\leq k,l\leq m} & \cr
(y+1)^{m-1} & & &}\nonumber\\[-8pt]\\[-8pt]
&&\qquad= \frac{(xy)^m-1}{xy-1}\nonumber
\end{eqnarray}
with $y+1=(\tilde a-\tilde b)/(\widetilde\omega+\tau_1-\tilde b)$ and
$x+1=(\tilde b-\tilde a)/(\omega+\tau_2-\tilde a)$.
Thus, one obtains for (\ref{eq4.39}):
%
%
\begin{eqnarray}\label{4.41}\qquad
&&\lim_{n\to\infty}\frac{1}{\sqrt{2} n^{1/6}}K_m(t_1,x_1;t_2,x_2)
\nonumber\\
&&\qquad= \frac{f(\tau_1,\xi_1)}{f(\tau_2,\xi_2)}\frac{1}{(2\pi i)^2}
\nonumber\\
&&\qquad\quad{}\times
\int
_{\Gamma_{\tilde a-\tau_2>}}d\omega\int_{\Gamma_{<\tilde b-\tau
_1}}d\widetilde\omega\frac{e^{-\omega^3/3+\xi_2\omega
}}{e^{-\widetilde\omega^3/3+\xi_1\widetilde\omega}}
\frac{(\tilde b-\tilde a)}{(\omega+\tau_2-\tilde
a)(\widetilde
\omega+\tau_1-\tilde b)} \\
&&\hspace*{-15pt}\hspace*{142.9pt}\times\biggl( \biggl(\frac{\widetilde\omega+\tau
_1-\tilde a}{\widetilde\omega+\tau_1-\tilde b}
\frac{\omega+\tau_2-\tilde b}{\omega+\tau_2-\tilde a} \biggr)^m-1\biggr)\nonumber\\
&&\hspace*{-15pt}\hspace*{142.9pt}{}\times\biggl(\frac
{\widetilde\omega+\tau_1-\tilde a}{\widetilde\omega+\tau_1-\tilde b}
\frac{\omega+\tau_2-\tilde b}{\omega+\tau_2-\tilde a}-1\biggr)^{-1},\nonumber
\end{eqnarray}
which is equal to the kernel (\ref{eq4.29}) obtained previously.
\begin{rem*}
Using (\ref{eq4tildePhi}), (\ref{eq4tildePsi}), a
1-border identity analogous to the $2$-border identity (\ref{4.40}),
Lemmas \ref{lemma4mu} and \ref{lemma4phipsi} and (\ref{eq4.40}), one
finds the limiting biorthogonal functions [which also yield (\ref{4.41})]:
\begin{eqnarray*}
\lim_{n\to\infty}2^{1/4}n^{1/12}\tilde\varphi_j^{n}(t,x)
&=&
\sqrt{\tilde b - \tilde a}\frac{(-1)^{j-1}}{2\pi i}
\int_{\Gamma_{\tilde a-\tau>}}\hspace*{-0.5pt}
d\omega\, e^{-{\omega^3}/{3}+\xi\omega}
\frac{(\omega+\tau-\tilde b_j)^{j-1}}{(\omega+\tau-\tilde a_j)^{j}},
\\
\lim_{n\to\infty}2^{1/4}n^{1/12}\tilde\psi_j^{n}(t,x)&=&
\sqrt{\tilde b - \tilde a}\frac{(-1)^{j-1}}{2\pi i}
\int_{\Gamma_{\tilde b-\tau>}}\hspace*{-0.5pt}
d\omega\, e^{ {\omega^3}/{3}-\xi\omega}
\frac{(\omega+\tau-\tilde a)^{j-1}}{(\omega+\tau-\tilde
b)^{j}}.
\end{eqnarray*}
\end{rem*}

\section{Limit to the Pearcey process}\label{SectPearcey}

In this section, we prove Theorem \ref{th:Pearcey}. To do this, we must
apply the scaling (\ref{scalingAB}) to the kernel with $m$ wanderers,
where all the $\tilde a_i=\tilde a$ and $\tilde b_i=\tilde b$ and where
one uses the shift $\omega\to\omega-\tau_2$ and $\widetilde\omega
\to
\widetilde\omega-\tau_1$, to yield
%
%
\begin{eqnarray}\label{eq5.1}\qquad
&&K_m^{\tilde a,\tilde b}(\tau_1,\xi_1;\tau_2,\xi_2)
\nonumber\\
&&\qquad=-\frac{\Id
(\tau
_2>\tau_1)}{\sqrt{4\pi(\tau_2-\tau_1)}}
e^{-{(\xi_2-\xi_1)^2}/({4(\tau_2-\tau_1)})-({1}/{2})(\tau
_2-\tau
_1)(\xi_2+\xi_1)+({1}/{12}) (\tau_2-\tau_1)^3}\nonumber\\[-8pt]\\[-8pt]
&&\qquad\quad{} + \frac{1}{ (2\pi i )^2}
\int_{\Gamma_{\tilde a>}}d\omega\int_{\Gamma_{<\tilde
b}}d\widetilde
\omega\frac{e^{-(\omega-\tau_2)^3/3+\xi_2(\omega-\tau
_2)}}{e^{-(\widetilde
\omega-\tau_1)^3/3+\xi_1 (\widetilde\omega-\tau_1)}}
\nonumber\\
&&\hspace*{-15.8pt}\hspace*{150.70pt}{}\times\frac{1}{\omega-\widetilde\omega} \biggl(\frac{\omega-\tilde b}{\omega
-\tilde a} \biggr)^m
\biggl(\frac{\widetilde\omega-\tilde a}{\widetilde\omega-\tilde b}
\biggr)^m.\nonumber
\end{eqnarray}
According to Theorem \ref{th:Pearcey}, we need to
take the following scaling limit:
%
%
\begin{eqnarray}
\tilde a&=&\alpha m^{1/3},\qquad \tilde b=\beta m^{1/3}, \nonumber\\
\tau_i&=&T_{\pm} m^{1/3}+\tfrac12\kappa^2 \theta_i m^{-1/6},\\
\xi_i&=&X m^{2/3}-\kappa^2\sigma_\pm\theta_i m^{1/6}-\kappa v_i
m^{-1/12}.\nonumber
\end{eqnarray}
Then, we have to compute the large $m$ limit of the rescaled kernel
(\ref{eq5.1}). We prove the following result, which implies Theorem
\ref{th:Pearcey}.
\begin{proposition}\label{Prop5.1}
Under the above scaling, for any fixed $\theta_1,\theta_2$, the limit
%
%
\begin{equation}\label{eq5.4}
\lim_{m\to\infty} \kappa m^{-1/12} K_m^{\tilde a,\tilde b}(\tau
_1,\xi
_1;\tau_2,\xi_2)\equiv K^\mathcal{P}(\theta_1,v_1;\theta_2,v_2)
\end{equation}
holds uniformly for $v_1,v_2$ in a bounded set.
\end{proposition}
\begin{pf}
The first term in (\ref{eq5.1}) is a straightforward limit. Indeed, for
$\theta_2>\theta_1$,
%
%
\begin{eqnarray}
&&-\frac{\kappa m^{-1/12} }{\sqrt{4\pi(\tau_2-\tau_1)}}e^{-
{(\xi
_2-\xi_1)^2}/({4(\tau_2-\tau_1)})-({1}/{2})(\tau_2-\tau_1)(\xi
_2+\xi
_1)+({1}/{12}) (\tau_2-\tau_1)^3} \nonumber\\[-8pt]\\[-8pt]
&&\qquad= \frac{-1}{\sqrt{2\pi(\theta_2-\theta_1)}} e^{ -
{(v_2-v_1)^2}/({2(\theta_2-\theta_1)}) } \frac{Q(1)}{Q(2)},\nonumber
\end{eqnarray}
where the conjugation terms $Q(i)$ are given by
%
%
\begin{equation}\label{eq5.5}
Q(i)=\exp\bigl(\tfrac12 \kappa^2 (\sigma^2+X) \theta_i m^{1/2}+\kappa
\sigma v_i m^{1/4}+\mathcal{O}(m^{-1/4}) \bigr).
\end{equation}

Next, one deals with the double integral, where it is natural to
introduce the change of integration variables:
%
%
\begin{equation}
\omega= w m^{1/3},\qquad \widetilde\omega=\widetilde w m^{1/3},
\end{equation}
leading to
%
%
\begin{eqnarray}\label{eq5.7}
&&\hspace*{25pt}\frac{\kappa m^{1/4} }{ (2\pi i )^2}
\int_{\Gamma_{\alpha>}}dw
\int_{\Gamma_{<\beta}}d\widetilde w \frac{1}{w-\widetilde w}\nonumber\\[-8pt]\\[-8pt]
&&\hspace*{25pt}\qquad
{}\times \frac{e^{m
F_0(w)+m^{1/2} F_2(w,\theta_2)+ m^{1/4} F_3(w,v_2)+ F_4(w,\theta
_2)+\mathcal{O}
(m^{-1/4})}}{e^{m F_0(\widetilde w)+m^{1/2} F_2(\widetilde w,\theta_1)+
m^{1/4} F_3(\widetilde w,v_1)+ F_4(\widetilde w,\theta_1)+\mathcal
{O}(m^{-1/4})}},\nonumber
\end{eqnarray}
where the functions $F_i$ are given by
%
%
\begin{eqnarray}\label{eq5.8}\qquad
F_0(w)&=&-\tfrac13 (w-T)^3+X(w-T)+\ln(w-\beta)-\ln(w-\alpha),\nonumber\\
F_2(w,\theta)&=&\tfrac12 \bigl((w-T)^2-X\bigr)\kappa^2 \theta-(w-T)\kappa
^2\sigma
\theta,\\
F_3(w,v)&=&-(w-T) \kappa v,\nonumber\\
F_4(w,\theta)&=&-\tfrac14(w-T-\sigma)\kappa^4 \theta^2.\nonumber
\end{eqnarray}
Setting
%
\begin{equation}\label{4.var}
w':=w-T,\qquad \alpha':=\alpha-T,\qquad \beta'=\beta-T,
\end{equation}
one defines
%
\begin{equation}\label{5.defF}\hspace*{33pt}
\tilde F_0(w'):=F_0(w'+T)=-\frac{w'^3}{3} +Xw'+\log(w'-\alpha')-\log
(w'-\beta').
\end{equation}
One now imposes the condition that $\tilde F'(w')$ experiences a
triple zero
at some critical point $w'_c$; this happens when the following
polynomial $P(w')$ is identically zero, with $w'_c\neq w'_0$:
%
%
\begin{eqnarray}\label{8}
\hspace*{35pt}0&\equiv& P(w')\nonumber
\\
&:=& -(w'-\alpha')(w'-\beta')\tilde F_0'(w')-(w'-w'_c)^3(w'-w'_1) \nonumber\\[-8pt]\\[-8pt]
&=&(3w'_c+w'_1-\alpha'-\beta')w'^3+(\alpha'\beta' -X-3{w'_c}^2-3w'_c w'_1){w'}^2\nonumber\\
&&{} + \bigl({w'_c}^3+3{w'_c}^2 w'_1+X(\alpha'+\beta')\bigr)w'-(\alpha'\beta'X-\alpha '+\beta'+{w'_c}^3 w'_1).\nonumber
\end{eqnarray}
Setting the coefficients of this cubic in $w'$ equal to $0$ amounts to 4
equations in 5~unknowns $\alpha', \beta', w'_c, w'_1, X$, thus yielding an algebraic
curve. At a first stage, let us look at it \textit{purely algebraically};
later we will have to take into account the real character of the
parameters, including various inequalities. Close inspection of the
four equations suggests the following birational map:
%
\begin{equation}\label{9}
\alpha'=\frac{ 2\sigma}{2-x }-r,\qquad \beta'=\frac{
2\sigma}{2-x}+r,\qquad
w'_1=\frac{ \sigma(3x -2)}{2-x}
\end{equation}
with inverse (assuming $\alpha'+\beta'\neq0$)
%
\begin{eqnarray}\label{10}
r&=&\frac{1}{2}(\beta'-\alpha'),\qquad x = \frac{2(2w'_1+\alpha'+\beta
')}{3(\alpha
'+\beta')},\nonumber\\[-8pt]\\[-8pt]
\sigma &=&-\frac{1}{3}(w'_1-\alpha'-\beta').\nonumber
\end{eqnarray}
Substituting this map in $P(w')$, one now solves the 4 equations (\ref
{8}) defined by $P(w')$ inductively, beginning with the highest degree
in $w'$. At first, one checks $P(w')=3(w'_c-\sigma)w'^3+\cdots,$ leading
to $w'_c=\sigma$, together with the value of $w'_1$, which we already
knew from (\ref{9}); thus
%
\begin{equation}\label{u0}
w'_c=\sigma\quad\mbox{and}\quad w'_1=\frac{ \sigma(3x -2)}{2-x}.
\end{equation}
Substituting this back into $P(w')$ yields a quadratic
polynomial in $w'$;
the vanishing of
the coefficient of $w'^2$ yields
%
\begin{equation}\label{11}
X=\frac{2\sigma^2(3x^2-6x +2)}{(x-2)^2}-r^2;
\end{equation}
$P(w')$ becomes thus linear, with vanishing linear and constant terms, yielding
%
\begin{equation}\label{12}
r^2=x^2\sigma^2\frac{(2x -3)}{(x -2)^2}
\quad\mbox{and}\quad r=\frac{2x^3\sigma^4}{2-x}.
\end{equation}
The compatibility between the $r$ and $r^2$ equations (\ref{12}) yields
a curve relating $x$ and $\sigma$,
%
\begin{equation}\label{14}
\mathcal{E}\dvtx\sigma^6=\frac{2x -3}{4x^4}.
\end{equation}
Incidentally, this curve is elliptic; indeed, viewed as a 6-fold cover
of the $x$-plane, the total ramification index equals $12$, with a
ramification of index $5$ above $x=3/2$, there are two ramification
points of index $2$ above $x=0$ and three simple branch points above
$x=\infty$; thus the genus $=1$. Then substituting the value (\ref
{12}) of
$r$ into (\ref{9}) and (\ref{11}), yields the following expressions for
$\alpha'$, $\beta'$ and $X$, all
defined on the algebraic curve (\ref{14}):
%
\begin{eqnarray}\label{4.0}
\alpha'&=&\frac{ 2\sigma}{2- x }(1-x^3\sigma^3)
,\nonumber\\[-8pt]\\[-8pt]
\beta'&=&\frac{ 2\sigma}{2- x }(1+x^3\sigma^3)\quad \mbox{and}\quad X =
{\sigma
^2}{(1-2x)}.\nonumber
\end{eqnarray}
Using these expressions, together with the value
of the critical point $w'_c=\sigma$, one checks from (\ref{5.defF}) that
%
\begin{equation} \label{4.F4}
\frac{1}{4!}\tilde F_0^{iv}(w'_c)=\frac{(x-1)}{2x^2\sigma},
\end{equation}
and thus
%
\begin{equation}\tilde F_0(w')= \tilde F_0(w'_c)+\frac
{(x-1)}{2x^2\sigma}(w'-w'_c)^4
+\mathcal{O}\bigl((w'-w'_c)^5 \bigr).
\end{equation}
One then \textit{requires the parameters} $\alpha',\beta',X, w'_c,w'_1$
\textit{to be
real with} $\alpha'<\beta'$, $w'_c\neq w'_1$, \textit{and} $\alpha'+\beta
'\neq0$.
This implies that $x$, $\sigma$ and $r$ must be real; the curve
relation (\ref{14}) yields two real solutions for $\sigma$, namely
$\sigma_+<0$ and $\sigma_-=-\sigma_+>0$. In
particular, from (\ref{14}), one must have $x > 3/2$, and since
$\alpha
'<\beta'$, one must have, according to (\ref{10}), that $2r=\beta
'-\alpha
'>0$ yielding from (\ref{12}) the inequality \mbox{$x <2$}. Thus, one has
$2>x >
3/2$. Moreover, (\ref{4.F4}) will be $<0$ for $\sigma=\sigma_+$ and
$>0$ for $\sigma=\sigma_-$; one then sets $\mp\frac{\kappa
^4}{4}:=\frac{(x-1)}{2x^2\sigma_\pm}$, with $\kappa>0$. Since
$\frac{2x
-3}{4x^4}$, the right-hand side of (\ref{14}), is an increasing
function of $3/2<x<2$, this function has its maximum at $x=2$, for
which $\sigma_\pm=\mp1/2$, according to the curve relation~(\ref{14}).
That $\alpha'+\beta'\neq0$ follows from adding the two first
equations in
(\ref{9}). Also one has $w'_c=\sigma_+>\beta'$ and $w'_c=\sigma
_-<\alpha
'$, since from (\ref{4.0}) and the curve (\ref{14}), one computes for
$\beta'-\sigma_+$ and similarly for $\alpha'-\sigma_-$
\[
\beta'-\sigma_+=\frac{\sigma_{_+} x}{2-x}(1+2x^2\sigma_+^3)=\frac
{\sigma
_+ x}{2-x}\bigl(1-\sqrt{2x-3}\bigr) <0
\quad\mbox{and}\quad
\alpha'-\sigma_- >0.
\]
Finally, from (\ref{u0}), it is clear that $w'_1\neq w'_c$, since
$\frac
{ \sigma(3x -2)}{2-x}\neq\sigma$ in the admissible range $x\in
(\frac32,2)$, with $w'_1\to w'_c$ for $x\to\frac32$.

To summarize, using the change of variables (\ref{4.var}), the
relations (\ref{4.0}) imply, for a given $\alpha<\beta$, two values
$T_\pm
$ of $T$ below, and thus
%
%
\begin{eqnarray}\label{0}
\alpha&=& T_\pm+\frac{2\sigma_\pm}{2-x}(1-x^3\sigma_\pm^3),\qquad
\beta=T_\pm+\frac{2\sigma_\pm}{2-x}(1+x^3\sigma_\pm^3),\nonumber
\nonumber\\[-8pt]\\[-8pt]
X &=& {\sigma_\pm^2}{(1-2x)},\qquad T_\pm=\frac{\alpha+\beta}{2}-\frac
{2\sigma_\pm}{2-x},\nonumber
\end{eqnarray}
from which (\ref{5.XT}) in the statement of Theorem \ref{th:Pearcey}
follows, with inequalities
%
%
\begin{equation}\hspace*{35pt}
3/2<x<2,\qquad 0<|\sigma_\pm|<\frac12,\qquad X<0,\qquad T_-
<\frac{\alpha+\beta}{2}<T_+.
\end{equation}
Also, the critical point $w_c$ of $F_0(w)$ and the extra-root $w_1$ of
$F'_0(w)$ occur at
%
%
\begin{eqnarray}\label{5.roots}
w_c:\!&=&\sigma_\pm+T_\pm\qquad\mbox{with } \matrix{{\alpha<\beta<w_c }\cr
{w_c<\alpha
<\beta}} \quad\mbox{and}\nonumber\\[-8pt]\\[-8pt]
w_1&=&T_\pm+\frac{ \sigma_\pm(3x -2)}{2-x}\neq
w_c .\nonumber
\end{eqnarray}

The statement about the \textit{uniqueness of the solution} $(x,\sigma)$
to the equations (\ref{0}) and (\ref{14}), given arbitrary $\alpha
<\beta$,
remains to be shown. Indeed, upon using the identities (\ref{0})
obtained for $\alpha$ and $\beta$, together with the curve equation
(\ref
{14}), it is easy to see that the right-hand side of the equation,
%
\begin{equation}\label{sol}
0<\beta-\alpha={\frac{4{ \sigma^4_\pm
}{x}^{3}}{2-x}}=\frac
{4x^{1/3} ( {2x-3} )^{2/3}}{2-x},
\end{equation}
is a monotonically increasing function in the range $3/2<x<2$;
therefore, the right-hand side of that equation takes on every value in
$(0,\infty)$ exactly once and thus given arbitrary $\alpha<\beta$, there
is a
unique value $x\in(3/2,2)$ satisfying the first equation in (\ref
{sol}). Substituting this value of $x$ in the $T_\pm$-equation of
(\ref
{0}), the value of
$T_\pm$ is specified unambiguously and thus $T_\pm$ can take on any
value in ${\mathbb R}$.
Therefore, only when $\alpha, \beta\to0$, do
%
\begin{eqnarray}\label{5.Quintic}\hspace*{45pt}
x&\to&3/2,\qquad \sigma_\pm\to0\quad\mbox{and}\nonumber\\[-8pt]\\[-8pt]
T_\pm &\to& 0\quad\mbox{and thus, by (\ref{5.roots})},\quad
w_c,w_1\to0,\nonumber
\end{eqnarray}
which proves the
remark at the end of Theorem \ref{th:Pearcey}.

Then the series expansions about the critical point $\sigma_\pm$ give
%
%
\begin{eqnarray}\label{eq5.13}
F_0(w)&=&F_0(w_c)\mp\tfrac14\kappa^4(w-w_c)^4+\mathcal{O}\bigl((w-w_c)^5\bigr),
\nonumber\\
F_2(w,\theta)&=&F_2(w_c,\theta)+\tfrac12\kappa^2 \theta(w-w_c)^2,
\nonumber\\[-8pt]\\[-8pt]
F_3(w,v)&=&F_3(w_c,v)-\kappa v (w-w_c), \nonumber\\
F_4(w,\theta)&=&\mathcal{O}(w-w_c).\nonumber
\end{eqnarray}
We now apply the steepest descent method, which we spell out for the
opening cusp; i.e., for $T_+$ and $\sigma=\sigma_+<0$. By Cauchy's
Residue theorem, one can deform the paths as indicated in Figure \ref
{FigPathPearcey}.
%
\begin{figure}

\includegraphics{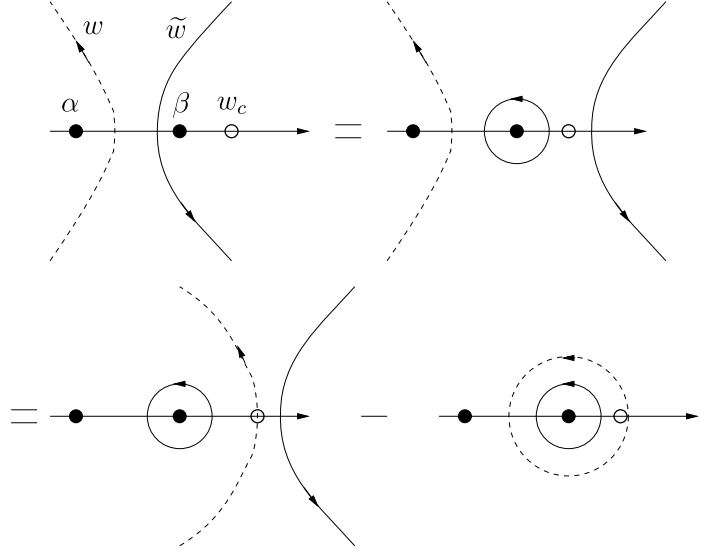}

\caption{First deformation of the paths which then pass close to the
critical point $w_c=\sigma_++T$. The solid contours are for
$\widetilde
w$, while the dashed ones for $w$. (For the case of $\sigma_-$, one has
$w_c<\alpha$ and the figures is essentially reflected.) The circles in the
first and second figures on the right-hand side become dashed and in
the second figure it sits on the right-hand side of the full curve. In
the third figure, the inner circle is dashed and the outer is full. In
effect, the roles of $w$ and $w'$ are interchanged.}
\label{FigPathPearcey}
\end{figure}
The contribution of the last contour is zero. Indeed, the integration
over $w$ is trivial, since the only pole is simple at $w=\widetilde w$.
All the factors involving $\alpha$ and $\beta$ cancel exactly. Thus, we
remain with a contour in $\widetilde w$ around $\beta$ of an analytic
function (no pole at $\beta$ anymore) which is zero. The deformation
also involves contributions which vanish at infinity.

The final and most important step is to deal with the previous last
contours of Figure \ref{FigPathPearcey}. First, a remark on the
integration paths in (\ref{eq5.7}). For large $w$ and $\widetilde w$,
the leading term is the cubic in $F_0$, which means that without any
error, we can let the the directions of the path $w$ go to infinity in
the cones of angles in $(\pi/2,5\pi/6)$ and $(-5\pi/6,-\pi/2)$ instead
of $2\pi/3$ and $-2\pi/3$. Similarly for $\widetilde w$, we can let it
go to infinity in the cones with angles in $(\pi/6,\pi/2)$ and $(-\pi
/2,-\pi/6)$ instead of $\pi/3$ and $-\pi/3$. Finally, the small contour
around $\beta$ can also be deformed to go to infinity as soon as it
does in directions $(5\pi/6,7\pi/6)$. Therefore, without errors, we
can deform the contours to become as in Figure \ref{FigSteepDescPearcey}.
%
%
\begin{figure}

\includegraphics{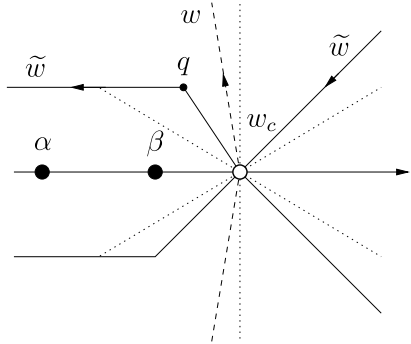}

\caption{Second deformation of the previous last contours in
Figure \protect\ref{FigPathPearcey}
for $\sigma=\sigma_+$. At the critical point $w_c$,
the path of $\widetilde w$ has angles $\pm\pi/4,\pm3\pi/4$, while the
path of $w$ leaves it with any angle in $(\pi/2,3\pi/4)$. The value of
$q$ will be chosen during the analysis. For $\sigma=\sigma_-$, one
picks the mirror image of the figure above about the vertical line
through $w_c$, with $w$ and $\tilde w$ also flipped;
that is, the dashed and solid lines are interchanged.}
\label{FigSteepDescPearcey}
\end{figure}
Let us verify that these paths satisfy the steepest descent property.
This will be done for the case $\sigma=\sigma_+$; the case $\sigma
=\sigma_-$ is essentially the same.

\textit{Slope of the function $F_0(w)$ starting from
$w_c=\sigma_+$.}\quad
Consider the curve given by
%
%
\begin{equation}
w=w_c+\zeta\frac{-\sigma x}{2-x}e^{\pm i(\pi/2+\delta)}\qquad \mbox{for
}0<\delta<\pi/3\mbox{ and }\zeta\geq0.
\end{equation}
Remember that $\sigma<0$. Then, at first, one verifies
%
\begin{equation}
\frac{\partial}{\partial\zeta}\operatorname{Re}F_0(w)
=\frac{x^2\sigma^3\zeta
^3P_3(\zeta
;x,\delta)}{ (2-x)^3 P_2(\zeta;x,-\sigma,\delta)P_2(\zeta;x,\sigma
,\delta)}<0\qquad
\mbox{for }\zeta>0,\hspace*{-37pt}
\end{equation}
with
%
%
\begin{eqnarray}
P_2(\zeta;x,-\sigma)&=&\zeta^2-2\zeta(1+2x^2\sigma^3)\sin\delta
+(1+2x^2\sigma^3)^2, \nonumber\\
P_3(\zeta;x,\delta)&=&\delta\bigl(3x +\mathcal{O}(\delta^2)\bigr) \zeta
^3+\bigl(2(2-x)+\mathcal{O}(\delta^2)\bigr)\zeta
^2\nonumber\\[-8pt]\\[-8pt]
&&{} +\delta\bigl((10x^2+4x -24)+\mathcal{O}(\delta^2)\bigr)
\zeta\nonumber\\
&&{} +\bigl(8(x-1)(2-x)+\mathcal{O}(\delta^2)\bigr).\nonumber
\end{eqnarray}
Indeed, $P_2(u;x,\pm\sigma)>0$, since its discriminant, as a quadric
in $u$, is $<0$. Moreover, $P_3(u;x,\delta)>0$, for
$0<\delta<\varepsilon(x)$ with $\varepsilon(x)$ sufficiently small,
since in that
case the
coefficients of $u^0,\ldots,u^3$ are positive [since $x\in(3/2,2)$] and
$\zeta\geq0$. Thus, the chosen path for $w$ is a path of steepest
descent.\vspace*{8pt}

\textit{Slope of the function $-F_0(\widetilde w)$ from
$w_c=\sigma_+$ to $q$ and $w_c$ to $e^{i\pi/4}\infty$.}\quad
Consider the curves parametrized by
%
%
\begin{equation}
\widetilde w=w_c+\frac{-\sigma x}{2-x}(\varepsilon\pm i)\zeta
\qquad\mbox{for
}\varepsilon
=\pm1\mbox{ and }\zeta\geq0.
\end{equation}
One verifies
%
%
\begin{equation}\label{V1}
\frac{\partial}{\partial\zeta}\operatorname{Re}(-F_0(\widetilde
w))=\frac{ 2x^3\sigma
^3\zeta
^3 P_3^{\varepsilon}(\zeta;x)}{ (2-x)^3 P_2^{\varepsilon}(\zeta
;x,-\sigma)P_2^{\varepsilon
}(\zeta;x,\sigma)},
\end{equation}
where, using the curve relation (\ref{14}),
%
%
\begin{equation}
P^{\varepsilon}_3(\zeta;x)=\varepsilon\zeta^3+2\zeta^2+\varepsilon
\zeta\frac
{-x^2+6x-4}{x}+\frac{4}{x}(x-1)(2-x)
\end{equation}
and
%
%
\begin{eqnarray}
|\widetilde w-\alpha|^2&=&2 \biggl(\frac{x\sigma}{x-2} \biggr)^2
P^{\varepsilon}_2(\zeta;x,-\sigma)>0,\nonumber\\[-8pt]\\[-8pt]
|\widetilde w-\beta|^2&=&2 \biggl(\frac{x\sigma}{x-2} \biggr)^2 P^{\varepsilon
}_2(\zeta;x,\sigma)>0\nonumber
\end{eqnarray}
with
%
\begin{equation}
P^{\varepsilon}_2(\zeta;x,-\sigma):= \zeta^2+\varepsilon\zeta
(1-2x^2\sigma
^3)-2x^2\sigma^3+x-1,
\end{equation}
showing at once the denominator of (\ref{V1}) is $>0$.

\textit{For $\varepsilon=1$}, the polynomial $P^{\varepsilon}_3$
in the
numerator of (\ref{V1}) is $>0$ for $\zeta>0$, because its coefficients
are all $>0$ in the range $2> x > 3/2$ and for $\zeta> 0$. Therefore,
the derivative (\ref{V1}) is $<0$ for $\zeta> 0$.

\textit{For $\varepsilon=-1$}, the polynomial $P^{\varepsilon
}_3(\zeta;x)$
will be $<0$ for large enough $\zeta>0$. However in the range $2> x > 3/2$,
%
\begin{equation}
-x^2+6x -4>0 \quad\mbox{and}\quad 0<\frac{4(x -1)(2-x)}{-x^2+6x
-4}<\frac{4}{11},
\end{equation}
and thus for $0<\zeta\leq\zeta_0$ with
%
\begin{equation}
\zeta_0:=\frac{4(x -1)(2-x)}{-x^2+6x-4},
\end{equation}
the cubic above is strictly positive:
%
\begin{eqnarray}\qquad
&&P_3^{\varepsilon}(\zeta;x) |_{\varepsilon=-1}
\nonumber\\[-8pt]\\[-8pt]
&&\qquad=
\zeta^2(2-\zeta)
-\frac
{1}{x}(-x^2+6x-4) \biggl(\zeta-\frac{4(x-1)(2-x)}{-x^2+6x-4} \biggr)
>0.\nonumber
\end{eqnarray}
This is the reason why for $\varepsilon=-1$ we bend the path at $q$ to be
horizontal, with $q$ set to be equal to
%
%
\begin{equation}
q:=w_c+\zeta_0 \frac{-\sigma x}{2-x}(\pm i-1).
\end{equation}

\textit{Slope of the function $-F_0(\widetilde w)$ from $q$ to
$q-\infty$.}\quad Consider the horizontal line given by
%
%
\begin{equation}
\widetilde w=q-\zeta\frac{-\sigma x}{2-x},\qquad \zeta\geq0.
\end{equation}
Then, in the range $-\sigma>0$ and $x\in(3/2,2)$,
%
\begin{eqnarray}
&&\frac{\partial}{\partial\zeta} \operatorname{Re}(-F(\widetilde w))
\nonumber\\[-8pt]\\[-8pt]
&&\qquad=\frac{ \sigma x}{2-x} \biggl[\frac{-\sigma x}{(2-x)^2(x^2-6x+4)} \biggr]^2
\frac{x\sigma^2 P_6(\zeta;x)}{|\widetilde w-\alpha|^2|\widetilde
w-\beta\nonumber
|^2}< 0
\end{eqnarray}
with
%
\begin{equation}
P_6(\zeta;x)=(2-x)^2Q_2(\zeta;x)+(2-x)\zeta^3 P_7(x)+\zeta^4\tilde
Q_2(\zeta;x).
\end{equation}
Indeed, $Q_2(\zeta;x)$ and $\tilde Q_2(\zeta;x)$ are quadratic
polynomials in
$\zeta$, with coefficients polynomial in $x$ and $P_7(x)$ is a seventh degree
polynomial in $x$. All three coefficients of $Q_2(\zeta;x)$ are $>0$
for $3/2\leq x\leq2$, while
$P_7(x)>0$ also as long as $2\geq x\geq3/2$. The coefficients of
$\zeta
^0$ and $\zeta^2$ of $\tilde
Q_2(\zeta;x)$ are $>0$ for $2\geq x\geq1.70$, which moreover has a
positive minimum in $\zeta$ for
$2\geq x \geq1.70$, thus proving the assertion for $2>x > 1.7$. A
little numerics in fact shows
we can remove the restriction $2>x >1.7$ and deduce the inequality for
$2>x >3/2$.

Thus, we have also shown that the chosen path for $\widetilde w$ is of
steepest descent. Thus, the steepest descent method can be applied along
these curves where the maximum of $\operatorname{Re}F_0(w)$,
$-\operatorname{Re}F_0(\widetilde
w)$ occur at the saddle point $w_c$. The main contribution comes from
the integration over a $\delta$-neighborhood of the critical point
$w_c$ for both $w$ and $\widetilde w$. For small $\delta$, the error
made is of order $e^{-\mu m}$ with $\mu\sim\delta^4$. Let us therefore
choose $\delta=\kappa^{-1} m^{-1/4} m^{\gamma}$ with any $\gamma\in
(0,1/20)$ fixed (i.e., $m^{-1/5}\gg\delta\gg m^{-1/4}$). Then, the
only nonvanishing contribution in the $m\to\infty$ limit is given by
the integrations with $|w-w_c|\leq\delta$, $|\widetilde w-w_c|\leq
\delta$. In these small neighborhoods, we can apply series expansions
(\ref{eq5.13}). After the change of variables
%
%
\begin{equation}
z:=\kappa m^{1/4} (w-w_c),\qquad \tilde z:=\kappa m^{1/4} (\widetilde w-w_c),
\end{equation}
we finally get for $\sigma=\sigma_\pm$,
%
%
\begin{equation}\label{eq5.30}\qquad
\frac{Q(2)}{Q(1)}(\ref{eq5.7}) =
\cases{\displaystyle{ \frac{1}{(2\pi i)^2}
\int_{{{\nwarrow\atop{\nearrow}}}
} dz
\int_{{\nwarrow\atop{\nearrow}} {\swarrow\atop{\searrow}}}
d\tilde z \frac{1}{z-\tilde z}\frac{e^{-z^4/4+\theta_2 z^2/2- v_2 z +
R_2}}{e^{-\tilde z^4/4+\theta_1 \tilde z^2/2- v_1 \tilde z + R_1}}
},\cr
\displaystyle{ \frac{1}{(2\pi i)^2}
\int_{{\nwarrow\atop{\nearrow}} {\swarrow\atop{\searrow}}}dz
\int_{{{\swarrow\atop{\searrow}}}
}
d\tilde z \frac{1}{z-\tilde z}\frac{e^{z^4/4+\theta_2 z^2/2- v_2 z +
R_2}}{e^{\tilde z^4/4+\theta_1 \tilde z^2/2- v_1 \tilde z + R_1}}
},}
\end{equation}
where $Q(i)=\exp(F_2(w_c,\theta
_i)m^{1/2}+F_3(w_c,v_i)m^{1/4}+\mathcal{O}
(m^{-1/4}) )$ is the conjugation given in (\ref{eq5.5}), and where
the $R_i$ are error terms, to be discussed later. Note the involution
$\theta_1\leftrightarrow-\theta_2, v_1\leftrightarrow
v_2, z\leftrightarrow-\tilde z$ between the two integrals on the
right-hand side of (\ref{eq5.30}), which also respect the integration
paths. The error terms $R_i$ include the following local
contributions:

(a) $\mathcal{O}(m^{-1/4})$ of (\ref{eq5.7}) (uniform for $v_i$ in a bounded
set),

(b) $\mathcal{O}(\delta)=\mathcal{O}(m^{-1/5})$ from $F_4(w)$ in
(\ref{eq5.13})
(uniform for $\theta_i$ in a bounded set),

(c) $\mathcal{O}(m \delta^5)=\mathcal{O}(m^{5(\gamma-1/20)})$,
which is the correction
in the series expansions of $F_0(w)$ of order higher than $4$, see
(\ref
{eq5.13}).

Indeed, (b) is immediate since $F_4$ is linear [see (\ref{eq5.8})]. To
see that (c) holds, we need to control the fifth derivative of $F_0$ at
$w_c$. We have
%
%
\begin{equation}\hspace*{33pt}
\max_{|w-w_c|\leq\delta} \bigl|F_0^{(v)}(w) \bigr| = 4 \max
_{|w-w_c|\leq\delta} \biggl|\frac{1}{(w-\alpha)^5}-\frac{1}{(w-\beta
)^5} \biggr| \leq\frac{16}{(w_c-\beta)^5}
\end{equation}
for $m$ large enough.

Finally, taking the $m\to\infty$ limit to (\ref{eq5.30}) the error
terms vanishes and at the same time the integrals extend to infinity.
The fact that $z$ is not exactly $i{\mathbb R}$ is irrelevant, since
the result
is identical as soon as the direction has an angle strictly smaller
than $\pi/4$ to the imaginary axis. Similarly, one can deform the
$\tilde z$-path as depicted in Figure \ref{Fig2}. This ends the proof of
Proposition \ref{Prop5.1} and thus also of Theorem~\ref{th:Pearcey}.
\end{pf}
\begin{rem}\label{Rem5}
For future use, we point out that the elliptic curve $\mathcal E$ has
three points above $x=\infty$, only one of which is real, namely
%
\begin{equation} \label{rem1}
\sigma= \frac{2^{-1/6}}{x^{1/2}}+\cdots\qquad\mbox{for } x\to\infty.
\end{equation}
At this point at infinity, one has, using the estimate (\ref{rem1}),
$
\beta-\alpha=\break{\lim_{x\to\infty}}\frac{4\sigma^4 x^3}{2-x}=
-2^{4/3}, $
and assuming $T=\frac{\alpha+\beta}2 -\frac{4\sigma}{2-x}=0$, also
$\alpha
+\beta=0$. This implies that
$\beta=-\alpha=-2^{1/3}$. Note how this contrasts with $(\sigma
,x)=(0,3/2), (X,T)=(0,0)$, in which case $\alpha=\beta=0$. To summarize,
near the real points on $\mathcal E$, namely near $x=3/2$ and $x=\infty$,
one has the following leading terms
(set
$\gamma:=\frac{2^{1/2}}{3^{2/3}}$):
%
%
\begin{eqnarray*}
(\sigma,x)&\sim&\biggl(\gamma\biggl(x - \frac32
\biggr)^{1/6},\frac32 \biggr),\qquad (\alpha,\beta)\sim0,\\
(X,T)&\sim&0,\qquad \frac{w_c-w_1}{-4\gamma}\sim\biggl(x-\frac32 \biggr)^{1/6},
\end{eqnarray*}
and
\begin{eqnarray}\label{estimate}
(\sigma,x)&\sim&\biggl(\frac{2^{-1/6}}{\sqrt{x}}, \infty\biggr) ,\qquad
(\alpha,\beta)\sim(2^{1/3},-2^{1/3}) , \nonumber\\[-8pt]\\[-8pt]
(X,T)&\sim& (-2^{2/3},0),\qquad \frac{w_c - w_1}{ 2^{-1/6} }
\sim\frac{4}{\sqrt{x}} .\nonumber
\end{eqnarray}
\end{rem}

\section{Limit to the quintic kernel}\label{SectConjecture}
In this section, we explain our guess concerning the process that
will occur in the situation illustrated in Figure~\ref{FigNewProcess}.
In Theorem~\ref{th:Pearcey} and, in particular, in formula (\ref
{5.Quintic}), it was observed that when $\alpha, \beta\to0$ (and only
then), the tips of the cusps $(\tau,\xi)\sim(\pm Tm^{1/3},Xm^{1/3})$
tend to the same point and that $ w_c-w_1=\frac{4\sigma_{\pm}
(1-x)}{2-x}\to0$, that is, the cube root of $F_0'(w)$ turns into a
quartic root. This also means that the starting and end points $a$ and
$b$ for the wanderers tend to coincide and that the line connecting
both points becomes vertical and tangent to the ellipse, as described
in Figure \ref{FigAsymmetric}. This corresponds to the first situation
in (\ref{estimate}). We now pick the second situation in (\ref
{estimate}), for which the cube root of $F_0'(w)$ also turns into a
quartic root. However, this forces the points $a$ and $b$ to be a bit
beyond $\sqrt{2n}$; this means in particular that $\tilde a> \tilde b$,
which actually violates the condition $\tilde a < \tilde b $ in
Theorem \ref{th:1}. One can think of the passage from $x=3/2$ to
$x=\infty
$ as a \textit{transition process}.

The most natural strategy would be to set $a=b=\sqrt{2n}+\sqrt{2m}$ and
take $m,n\to\infty$ together. Then, under an appropriate scaling limit,
we expect to get a process with a \textit{quintic kernel}. Of course,
there will be a parameter tuning regulating how close the two Airy
fields come together. For example, if $m=n$, then we have to choose
$a=2\sqrt{2n}+ \mathcal{O}(n^{-1/6})$, since the fluctuations of the
first $n$
Brownian bridges alone live on the $n^{-1/6}$ scale.

\textit{Evidence in favor of Conjecture \ref{conjecture}.}\quad
To give some evidence to this conjecture, we present two pieces of
rigorous mathematics, concerning the (one-time) kernel, with $\tilde
a<\tilde b$, with time $\tau$ absorbed into $\tilde a$,
%
%
\begin{eqnarray}\label{eqKernelExtended2}
&&K_m^{\tilde a,\tilde b}( \xi_1; \xi_2)\nonumber\\[-1pt]
&&\qquad= \frac{1}{ (2\pi i )^2}
\int_{\Gamma_{\tilde a >}}d\omega\int_{\Gamma_{<\tilde
b}}d\widetilde
\omega\frac{e^{-\omega^3/3+\xi_2\omega}}{e^{-\widetilde\omega
^{3}/3+\xi_1
\widetilde\omega}}
\\[-1pt]
&&\hspace*{-15.2pt}\hspace*{104.5pt}\hspace*{32.8pt}{}\times
\frac{ (({\widetilde\omega-\tilde a})/({\omega-\tilde a} ))^m
(({\omega-\tilde b})/({\widetilde\omega-\tilde b} ))^m}{ \omega
-\widetilde\omega}.\hspace*{-22pt}\nonumber
\end{eqnarray}
\begin{proposition}\label{prop:analyticcont}The kernel $K_m^{\tilde
a,\tilde b}$, as in (\ref{eqKernelExtended2}), can be continued
analytically to a new kernel $\widetilde K_m^{\tilde a,\tilde b}$, as
in (\ref{eqKernelExtended3}), with same integrand as kernel (\ref
{eqKernelExtended2}), by moving $\tilde a$ and $\tilde b$ in the
%
\begin{figure}

\includegraphics{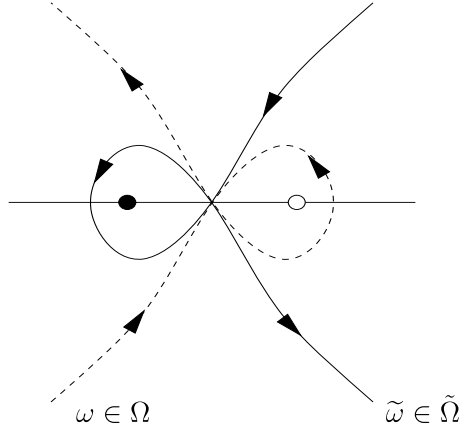}

\caption{New contour $\Omega$ and $\tilde\Omega$, with the \textit{black
dot} $= \tilde b$ and the \textit{white} $\mbox{\textit{dot}} = \tilde a$. The
solid line $\widetilde\Omega$ refers to the integration of the
$\widetilde\omega$-variable, while the dashed line $ \Omega$ refers to
the $ \omega$-integration.}
\label{newcontour}
\end{figure}
%
\begin{figure}[b]

\includegraphics{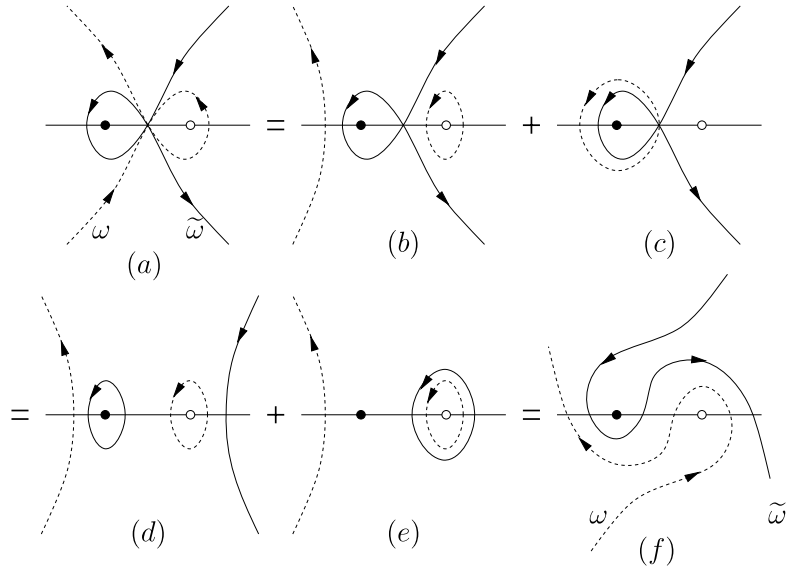}

\caption{Representation of the deformation of the integration variables
for the case $\tilde b<\tilde a$. All the contours are clockwise
oriented, the black dot is $\tilde b$, the white dot is $\tilde a$. The
solid line refers to the integration of the $\widetilde\omega
$-variable, while the dashed line for $ \omega$-integration. The
contributions of \textup{(c)} and \textup{(e)} are exactly zero.}
\label{FigDeformationQuintic}
\end{figure}
complex plane from their original position $\tilde a <\tilde b$ to a
new position $\tilde b <\tilde a$ on the real line:
%
%
\begin{eqnarray}\label{eqKernelExtended3}
&&\widetilde K_m^{\tilde a,\tilde b}( \xi_1; \xi_2)\nonumber\\[-1pt]
&&\qquad
= \frac{1}{ (2\pi i )^2}
\int_{\Omega}d\omega\int_{\tilde\Omega}d\widetilde\omega\frac
{e^{-\omega^3/3+\xi_2\omega}}{e^{-\widetilde\omega^{3}/3+\xi_1
\widetilde
\omega}}\\[-1pt]
&&\qquad\quad\hspace*{-15pt}\hspace*{156pt}\hspace*{-68.1pt}{}\times
\frac{ (({\widetilde\omega-\tilde a})/({\omega-\tilde a}) )^m
(({\omega-\tilde b})/({\widetilde\omega-\tilde b}) )^m}{ \omega
-\widetilde\omega}.\hspace*{-22pt}\nonumber
\end{eqnarray}
integrated over contours $\Omega$ and $\widetilde\Omega$ as in
Figure \ref{newcontour}.
\end{proposition}
\begin{pf}
$\tilde b$ corresponds to the black dot and $\tilde a$ to the
white dot in Figures \ref{newcontour} and \ref{FigDeformationQuintic};
the dashed line refers to the $\omega$-integration and the solid line to
the $\tilde\omega$-integration.

We noticed in Remark \ref{Rem5} that the elliptic curve $\mathcal{E}$,
introduced in (\ref{14}), contains another real point, namely one
covering $x=\infty$ for which $(\alpha,\beta)=(2^{1/3},-2^{-1/3})$. This
clearly violates the inequality $\tilde a=\alpha m^{1/3}<\tilde b=\beta
m^{1/3}$, crucial for the derivation of the kernel (\ref{eqKernelExtended2}).

Keeping $\tilde\omega$ fixed but arbitrary on the solid line,
one sees that the dashed line of Figure \ref{FigDeformationQuintic}(a)
can be deformed into the dashed lines of Figure \ref
{FigDeformationQuintic}[(b)${}+{}$(c)]. Then one notices that the
(c)-contribution vanishes. Indeed, (i) if $\tilde\omega$ belongs to the
solid line, outside the dashed circle, the $\omega$-integral vanishes, the
integrand being holomorphic; (ii)\vspace*{1pt} if $\tilde\omega$ belongs to the solid
line, inside the dashed circle, one picks up a residue and thus the
$\omega
$-integral equals $\frac{1}{2\pi i}e^{(\xi_2-\xi_1)\tilde\omega}$;
further integrated with regard to $\widetilde\omega$, one obtains
\[
\frac{1}{2\pi i} \int_{\mathrm{solid}\ \mathrm{circle}\ \mathrm{of}\ (\mathrm{c})} d\tilde
\omega\,
e^{(\xi_2-\xi_1)\tilde\omega}
=0
\]
and thus the only contribution comes from (b).

At the next stage, picking an arbitrary $\omega\in$ dashed
contour (b), one deforms the solid contour (b) into the solid contours
(d) + (e). In the same way, if $\omega\notin$ dashed circle, the
$\tilde
\omega$-integration contributes nothing, the integrand being holomorphic;
if $\omega\in$ dashed circle, the $\tilde\omega$-integration contributes
$\frac{1}{2\pi i}e^{(\xi_2-\xi_1) \omega}$; further integrated with
regard to $\omega$, one obtains
\[
-\frac{1}{2\pi i} \int_{\mathrm{dashed}\ \mathrm{circle}\ \mathrm{of}\ (\mathrm{e})} d \omega
\,e^{(\xi_2-\xi_1) \omega}
=0 ,
\]
and thus the integration over the (d)-contour is the only contribution.
Finally, the solid and dashed contours of (d) can further be deformed
into contours (f), thus, leading to the contours of Figure \ref{Fig1},
as the black dot $\tilde b$ migrates to the right of the white dot
$\tilde a$ through the ${\mathbb C}$-plane; this ends the proof of
Proposition \ref{prop:analyticcont}.
\end{pf}
\begin{proposition}\label{quintic}
Consider the kernel $\widetilde K_m^{\tilde a,\tilde b}( \xi_1; \xi
_2)$, as in (\ref{eqKernelExtended3}) with $\tilde a>\tilde b$.
Then defining the scaling
%
%
\begin{eqnarray}\label{eqScalingQuintic}
\tilde a&=&(2m)^{1/3} \bigl( 1+ \tfrac16 \theta m^{-2/5}+\tfrac{1}{2}\eta
m^{-3/5} \bigr),\nonumber\\
\tilde b&=&(2m)^{1/3} \bigl( -1- \tfrac16 \theta m^{-2/5}+\tfrac
{1}{2}\eta m^{-3/5} \bigr),\\
\xi_i&=&-(2m)^{2/3} \bigl(1-\tfrac16 \theta m^{-2/5}-\tfrac12
\bigl(v_i-\tfrac{1}{18}\theta^2 \bigr) m^{-4/5} \bigr),\nonumber
\end{eqnarray}
one obtains, in the $m\to\infty$ limit, the quintic kernel
$K^\mathcal
{Q}(\xi_1,\xi_2)$,
%
%
\begin{eqnarray}\label{eqQuinticKernel}
&&\lim_{m\to\infty}\frac{2^{-1/3} m^{-2/15}}{ (2\pi i )^2}
\int_{\Omega}d\omega\int_{\tilde\Omega}d\widetilde\omega
\frac{e^{-\omega^3/3+\xi_2\omega}}{e^{-\widetilde\omega^3/3+\xi_1\widetilde\omega}}
\frac{1}{\omega-\widetilde\omega} \biggl(\frac{\omega-\tilde b}{\omega-\tilde a} \biggr)^m
\biggl(\frac{\widetilde\omega-\tilde a}{\widetilde\omega-\tilde b}
\biggr)^m\hspace*{-28pt}\nonumber\\
&&\qquad=\frac{1}{(2\pi i)^2}\int_\mathcal{C}dz\int_{\widetilde{\mathcal{C}}}d\tilde z \frac{1}{z-\tilde z}
\frac{e^{(2/5) z^5-(1/3) \theta z^3-\eta z^2+v_2 z }}{e^{(2/5)\tilde z^5-
(1/3) \theta\tilde z^3-\eta\tilde z^2+v_1\tilde z }}\\
&&\qquad=: K^\mathcal{Q}(\theta,\eta;v_1,v_2),\nonumber
\end{eqnarray}
where $\mathcal C$ and $\widetilde{\mathcal{C}}$ are the paths
defined in
Figure \ref{QuinticProcess}.
The limit is uniform for $\theta,\eta, v_1,v_2$ in a bounded set.
\end{proposition}
\begin{pf}
We shall give the proof in the case of $\eta=0$; the case $\eta\neq0$
is easy to implement. As in the case of the Pearcey process (see
Theorem \ref{th:Pearcey}), consider the scaling $\xi_i=X m^{2/3}$,
$\tilde a=\alpha m^{1/3}$, $\tilde b=\beta m^{1/3}$ and the change of
integration variables $\omega= w m^{1/3}$, $\widetilde\omega=
\widetilde w m^{1/3}$. Then the kernel (\ref{eqKernelExtended3}) becomes
%
%
\begin{equation}\label{6.5}
(\ref{eqKernelExtended3})=\frac{m^{1/3}}{ (2\pi i )^2}
\int_{\Omega}dw \int_{\widetilde\Omega}d\widetilde w \frac{e^{m
F(w)-m F(\widetilde w)}}{w-\widetilde w}
\end{equation}
with
%
%
\begin{equation}\label{6.F}
F(w):=-w^3/3+Xw+\ln(w-\beta)-\ln(w-\alpha),
\end{equation}
where $\Omega$ and $\widetilde\Omega$ are the contours of Figure
\ref
{newcontour}, with the black dot being $\beta$ and
the white dot $\alpha$. Here, one imposes the property that $F'(w)$
experiences a 4-fold zero at some point $w_c$,
with $\alpha,\beta,X$ real and $\alpha\neq w_c$, $\beta\neq w_c$;
that is, one
requires all $A_i=0$:
%
%
\begin{eqnarray}
&&-(w-\alpha)(w-\beta)F'(w)-(w-w_c)^4\nonumber\\
&&\qquad=:A_0w^3+A_1w^2+A_2w+A_3
\nonumber\\[-8pt]\\[-8pt]
&&\qquad= (4w_c-\alpha-\beta)w^3+ (\alpha\beta-6w_c^2-X)w^2\nonumber\\
&&\qquad\quad{} +\bigl(4w_c^3+X(\alpha+\beta)\bigr)w-X\alpha\beta+\alpha-\beta
-w_c^4.\nonumber
\end{eqnarray}
The coefficients $A_0=A_1=0$ imply
$w_c=\frac{1}{4} (\alpha+\beta)\mbox{ and }
X=\alpha\beta-6w_c^2 $
and consequently
$A_2=-\frac{1}{16} (\alpha+\beta)(5\alpha^2-6\alpha\beta+5\beta^2)=0
$, whose only real solution is given by $\alpha=-\beta$ and thus $w_c=0$ and
$X=\alpha\beta=-\alpha^2$.
For these values, one has $A_3=-\alpha(\alpha^3-2)=0$, implying
$\alpha=-\beta=2^{1/3}$, $X=-2^{2/3}$ and $w_c=0$.
To summarize,
%
\begin{figure}[b]

\includegraphics{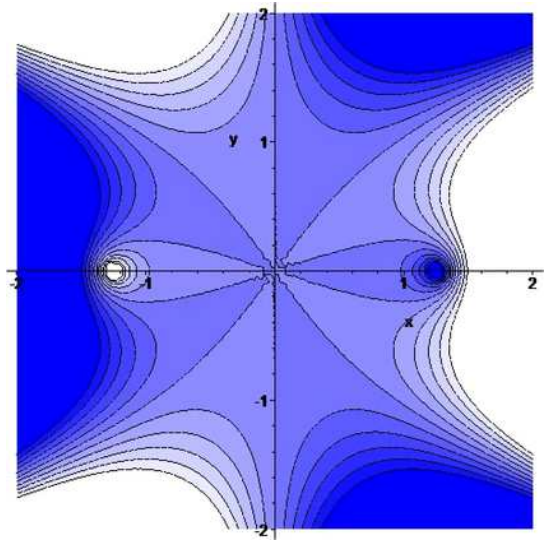}

\caption{Contourplot of the function $\operatorname{Re}(F(x+iy)-F(0))$.
The value is high in dark regions and low in light regions.}
\label{Contourplot}
\end{figure}
\[
\beta=-2^{1/3}<w_c=0<\alpha=2^{1/3} \quad\mbox{and}\quad X=-2^{2/3}.
\]
Note this solution corresponds precisely to the real point on the
elliptic curve~$\mathcal E$, covering $x=\infty$, as obtained on the
second line of (\ref{estimate}) (see Remark \ref{Rem5}).
Since we have a quintic leading term $\sim m w^5$, we make the change
of variables $w=m^{-1/5} z \alpha$ and $\widetilde w=m^{-1/5} \tilde z
\alpha$.
The precise coefficients are chosen in order to simplify the final
formula. Indeed, with (\ref{eqScalingQuintic}), we obtain
%
%
\begin{equation}\quad
m F(w) = m F(0)+v_2 z-\theta z^3/3 +2z^5/5+\mathcal{O}(z^7m^{-2/5},zm^{-2/5})
\end{equation}
with the error uniform for $\theta,v_i$ in a bounded set. The prefactor
in (\ref{6.5}), after the changes of variables, becomes $m^{1/3}
m^{-1/5} 2^{1/3}(1+\mathcal{O}(m^{-2/5}))$, which cancels with the
$2^{-1/3}m^{-2/15}$ in front of the left-hand side of (\ref{eqQuinticKernel})
(as $m\to\infty$). Except for the error terms, the result of the
theorem would follow.

What remains to be seen is that the higher order expansions in the
series do not contribute. We do it by the steepest descent method as
for the Pearcey case. Consider the curve parametrized by $w=e^{\pm3\pi
/5}x$. Then for the function $F$, as in~(\ref{6.F}), with $\alpha,
\beta
$ and $X$ substituted,
%
\begin{equation}
F(w)=-\tfrac13 w^3-{2}^{2/3}w+\ln( w+2^{1/3} ) -\ln
( w-2^{1/3} )
\end{equation}
one checks
%
\begin{equation}
\frac{\partial}{\partial x}\operatorname{Re}F(w)=-2{\frac{x^4
(x^2\cos(\pi/5)+1
)}{x^4+2x^2\cos(\pi/5 )+1}}<0 \qquad\mbox{for all }x>0.
\end{equation}
One then checks that along the dotted loop in Figure \ref{newcontour},
$\operatorname{Re}F(w)-\operatorname{Re}F(0)<0$ for $w\neq0$. It is
at once visible by
superimposing the dashed contour\vspace*{1pt} of Figure \ref{newcontour} onto the
contour plot, as in Figure \ref{Contourplot}.
Then along the curve given by $\widetilde w=e^{\pm2\pi/5}x$,
%
\begin{equation}\qquad
\frac{\partial}{\partial x}(-\operatorname{Re}F(\widetilde
w))=-2{\frac{x^4 (x^2\cos(\pi
/5)+1 )}{x^4+2x^2\cos(\pi/5 )+1}}<0 \qquad\mbox{for all }x>0,
\end{equation}
and along the solid loop in Figure \ref{newcontour}, $-\operatorname
{Re}F(\tilde
w)+\operatorname{Re}F(0)<0$ for $\tilde w\neq0$. This shows that the
curves have
the steepest descent property. Thus, if we integrate (in~$w$) around a
$\delta$-neighborhood of the origin, the error term will be only of
order $\mathcal{O}(e^{-\mu m})$ with $\mu\sim\delta^5$. We choose
$\delta
=m^{-1/5} m^{\gamma}$ for\vspace*{1pt} any $\gamma\in(0,2/35)$. Then, uniformly for
$\theta,v_1,v_2$ in a bounded set, the error term $\mathcal{O}(z^7
m^{-2/5})=\mathcal{O}
(m^{7\gamma-2/5})\to0$ as $m\to\infty$, as $z<\delta m^{1/5}/\alpha
<m^{\gamma
}/\alpha$. In the limit, the only part of the contour in Figure \ref
{newcontour}, which contributes in the end are the $8$ rays emanating
from the origin, which one deforms so as to form consecutive angles
$\pi
/5$. This then yields the quintic kernel $K^\mathcal{Q}(\theta,\eta
;v_1,v_2)$ with the integration paths $\mathcal C$ and $\widetilde
{\mathcal{C}}$ of Figure \ref{QuinticProcess}, thus establishing
Proposition \ref{quintic}.
\end{pf}
%

\section*{Acknowledgments}
The support of a National Science Foundation Grant DMS-07-04271 is
gratefully acknowledged (Mark Adler and Pierre van Moerbeke).
Also, a European Science Foundation
grant (MISGAM), a Marie Curie Grant (ENIGMA), a FNRS grant and a
Belgian ``Interuniversity Attraction Pole'' grant are gratefully acknowledged.
Patrik Ferrari started the work while being at the Weierstrass Institute for
Applied Analysis and Stochastics (WIAS), Berlin.

%

%
\printaddresses

\end{document}